\documentclass[12pt]{book}
\usepackage{amsmath,amsthm,amsfonts,amssymb,makeidx}

\def\myindex#1{{\em #1}\index{#1}}
\def\xL2{L^2[-1,1]}
\numberwithin{equation}{section}
\theoremstyle{plain}
\newtheorem{thm}{Theorem}[section]
\newtheorem{lemma}[thm]{Lemma}
\newtheorem{cor}[thm]{Corollary}
\newtheorem{prop}[thm]{Proposition}
\newtheorem{defn}[thm]{Definition}

\newtheorem*{completeness_axiom}{The Completeness Axiom}
\newtheorem*{axiom_of_choice}{The Axiom of Choice}
\newtheorem*{notn}{Notation}

\theoremstyle{definition}
\newtheorem{exer}[thm]{Exercise}

\newtheorem{ex}[thm]{Example}
\newtheorem{exs}[thm]{Examples}

\DeclareMathOperator{\dist}{dist}

\DeclareMathOperator{\len}{len}
\newcommand{\la}{\langle}
\newcommand{\ra}{\rangle}
\newcommand{\ms}{{\mu^{*}}}
\newcommand{\A}{{\mathcal A}}
\newcommand{\B}{{\mathcal B}}
\newcommand{\C}{{\mathcal C}}
\newcommand{\D}{{\mathcal D}}
\newcommand{\F}{{\mathcal F}}
\newcommand{\HH}{{\mathcal H}}
\newcommand{\M}{{\mathcal M}}
\newcommand{\PP}{{\mathcal P}}
\newcommand{\V}{{\mathcal V}}
\newcommand{\U}{{\mathcal U}}
\newcommand{\yL}{{\mathcal L}}
\newcommand{\I}{{I}}
\newcommand{\RR}{{\mathcal R}}
\newcommand{\R}{\mathbb R}
\newcommand{\N}{\mathbb N}
\newcommand{\Q}{\mathbb Q}

\newcommand{\X}{\mathfrak{X}}
\newcommand{\Z}{\mathbb Z}
\newcommand{\eps}{\epsilon}

\title{Notes on Measure and Integration}
\author{John Franks\\
Department of Mathematics\\
Northwestern University}
\date{}
\makeindex
\begin{document}
\maketitle
\pagenumbering{roman}
\chapter*{Preface}
This text grew out of notes I have used in teaching a one quarter course
on integration at the advanced undergraduate level.
My intent is to introduce the Lebesgue integral in a quick, and
hopefully painless, way and then go on to investigate the standard
convergence theorems and a brief introduction to the Hilbert space
of $L^2$ functions on the interval.  

The actual construction of Lebesgue measure and proofs of its key
properties are relegated to an appendix.  Instead the text introduces
Lebesgue measure as a generalization of the concept of length and
motivates its key properties: monotonicity, countable additivity, and
translation invariance.  This also motivates the concept of
$\sigma$-algebra. If a generalization of length has these
properties then to make sense it should be defined on a
$\sigma$-algebra.

The text introduces null sets (sets of measure zero) and 
shows that any generalization of length satisfying 
monotonicity must assign zero to them.
We then {\em define} Lebesgue measurable sets to be sets in 
the $\sigma$-algebra generated by Borel sets and null sets.

At this point we state a theorem which asserts that Lebesgue
measure exists and is unique, i.e. there is a function $\mu$
defined for measurable subsets of a closed interval which satisfies
monotonicity, countable additivity, and translation invariance.

The proof of this theorem (Theorem~(\ref{thm: lebesgue measure})) is
included in an appendix where it is also shown that the more common
definition of measurable sets using outer measure is equivalent to
being in the $\sigma$-algebra generated by Borel sets and null sets.

The text presupposes a background which a student obtain from
an undergraduate course in real analysis.
Chapter~0 summarizes these prerequisites with many proofs and
some references.  Chapter~1 gives a brief treatment of 
the ``regulated integral'' (as found in
Dieudonn\'e \cite{D}) and the Riemann integral in a way that
permits drawing parallels with the presentation of the Lebesgue integral in
subsequent chapters.  Chapter~2 introduces Lebesgue measure in
the way described above.  

Chapter~3 discusses bounded Lebesgue measurable functions and
their Lebesgue integral, while Chapter~4 considers unbounded
functions and some of the standard convergence theorems.
In Chapter~5 we consider the Hilbert space of $L^2$ functions
on $[-1,1]$ and show several elementary properties leading up
to a definition of Fourier series.

In Appendix~A we construct Lebesgue measure and prove
it has the properties cited in Chapter 2.  Finally in Appendix B
we construct a non-measurable set.

\medskip

\noindent
This work is licensed under a 
Creative Commons Attribution-Noncommercial-Share Alike 3.0 
United States License.
\medskip

\noindent
See http://creativecommons.org/licenses/by-nc-sa/3.0/us/

\tableofcontents

\setcounter{chapter}{-1}
\chapter{Background and Foundations}
\setcounter{page}{1}
\pagenumbering{arabic}
\section{The Completeness of $\R$}

This chapter gives a very terse summary of the properties of the real
numbers which we will use throughout the text.  It is intended as
a review and reference for standard facts about the real numbers
rather than an introduction to these concepts.

\begin{notn}
We will denote the set of \myindex{real numbers} by $\R$, the 
\myindex{rational numbers} by $\Q$, the \myindex{integers} by $\Z$ and the 
\myindex{natural numbers} by $\N.$
\index{$\R$ real numbers}
\index{$\Q$ rational numbers}
\index{$\Z$ integers}
\index{$\N$ natural numbers}
\end{notn}

In addition to the standard properties
of being an ordered field (i.e. the properties of arithmetic)
 the real numbers $\R$ satisfy a property
which makes analysis as opposed to algebra possible.

\begin{completeness_axiom}
Suppose $A$ and $B$ are non-empty subsets of $\R$ such that $x \le y$
for every $x \in A$ and every $y \in B$.  Then there exists
at least one real number $z$ such that $x \le z$ for all
$x \in A$ and $z \le y$ for all $y \in B.$
\end{completeness_axiom}

\begin{ex}
The rational numbers, $\Q,$ fail to satisfy this property.  If
$A = \{x\ |\ x^2 < 2\}$ and $B = \{y\ |\ y >0 \text{ and } y^2 > 2\},$
then there is no $z \in \Q$ such that $x \le z$ for all
$x \in A$ and $z \le y$ for all $y \in B.$
\end{ex}

\begin{defn}[Infimum, Supremum]
 If $A \subset \R,$ then $b \in \R$ is called an 
\myindex{upper bound} for $A$ if $b \ge x$ for all $x \in A$.
The number $\beta$ is called the \myindex{least upper bound} or
\myindex{supremum} of the set $A$ if $\beta$ is an upper bound
and $\beta \le b$ for every upper
bound $b$ of $A$.  A number $a \in \R$ is called a
\myindex{lower bound} for $A$ if $a \le x$ for all $x \in A$.
The number $\alpha$ is called the \myindex{greatest lower bound} or
\myindex{infimum} of the set $A$ if $\alpha$ is a lower bound and
$\alpha \ge a$ for every lower
bound $a$ of $A$.
\end{defn}

\begin{thm}\label{thm: sup exists}
If a non-empty set $A \subset \R$ has an upper bound, then it has a unique 
supremum $\beta$.  If $A$  has a lower bound, then it has a unique 
infimum $\alpha$.
\end{thm}

\begin{proof}
Let $B$ denote the non-empty set of upper bounds for $A$.  Then
$x \le y$ for every $x \in A$ and every $y \in B$.  The Completeness
Axiom tells us there is a $\beta$ such that 
$x \le \beta \le y$ for every $x \in A$ and every $y \in B$.
This implies that $\beta$ is an upper bound of $A$ and that
$\beta \le y$ for every upper bound $y$.  Hence $\beta$ is a supremum
or least upper bound of $A$.  It is unique, because any $\beta'$
with the same properties must satisfy $\beta \le \beta'$ (since
$\beta$ is a least upper bound) and $\beta' \le \beta$ (since
$\beta'$ is a least upper bound).  This, of course implies 
$\beta = \beta'.$

The proof for the {\it infimum} is similar.
\end{proof}

We will denote the {\it supremum} of a set $A$ by $\sup A$ and the
{\it infimum}  by $\inf A$.

\begin{prop}\label{prop: sup}
If $A$ has an upper bound and $\beta = \sup A$, then for any $\eps >0$
there is an $x \in A$ with $\beta - \eps < x \le \beta$.  Moreover $\beta$ is
the only upper bound for $A$ with this property.  If $A$ has a lower
bound its {\it infimum} satisfies the analogous property.
\end{prop}

\begin{proof}
If $\beta = \sup A$ and there is no $x \in (\beta - \eps, \beta),$
then every $x \in A$ satisfies $x \le \beta -\eps.$  It follows
that $\beta -\eps$ is an upper bound for $A$ and is smaller than
$\beta$ contradicting the definition of $\beta$ as the least upper
bound.  Hence there must be an $x \in A$ with $x \in (\beta - \eps, \beta)$.

If $\beta' \ne \beta$ is another upper bound for $A,$ then $\beta' > \beta.$
There is no $x \in A$ with 
$x \in(\beta, \beta']$, since such an $x$ would be greater than
$\beta$ and hence $\beta$ would not be an upper bound for $A$.

The proof for the {\it infimum} is similar.
\end{proof}

\section{Sequences in $\R$}

There are a number of equivalent formulations we could have chosen for
the Completeness Axiom.  For example, we could have take Theorem
(\ref{thm: sup exists}) as an axiom and then proved the Completeness
Axiom as a theorem following from this axiom.  In this section we
prove several more theorems which we will derive from the Completeness
Axiom, but which are in fact equivalent to it in the sense that if we
assumed any one as an axiom we could prove the others as consequences.
Results of this type include Theorem (\ref{thm: monotone limits}),
Corollary (\ref{cor: limsup}), and Theorem (\ref{thm: cauchy}).

We recall the definition of limit of a sequence.

\begin{defn}
Suppose $\{x_n\}_{n=1}^\infty$ is a sequence in $\R$ and $L \in \R$.  We say
\[
\lim_{n \to \infty} x_n = L
\]
provided for every $\eps > 0$ there exists $N \in \N$
such that
\[
|x_m - L| < \eps
\]
for all $m \ge N.$
\end{defn}

Let $\{x_n\}_{n=1}^\infty$ be a sequence in $\R$.  We will 
say it is \myindex{monotone increasing} if $x_{n+1} \ge x_n$ for 
all $n$ and \myindex{monotone decreasing} if $x_{n+1} \le x_n$ for 
all $n.$  
\begin{thm}\label{thm: monotone limits}
If $\{x_n\}_{n=1}^\infty$ is a bounded \myindex{monotone} sequence
then $\displaystyle{\lim_{n \to \infty} x_n}$ exists.  
\end{thm}
\begin{proof}
If $\{x_n\}_{n=1}^\infty$ is a bounded monotone increasing
sequence, let $L = \sup\{x_n\}_{n=1}^\infty$.  Given any $\eps >0$
there is an $N$ such that $L - \eps < x_N \le  L$ by Proposition
(\ref{prop: sup}).  For any $n > N$ we have $x_N \le x_n \le L$
and hence $|L - x_n| <\eps$.  Thus $\displaystyle{ \lim_{n \to \infty}
x_n = L.}$  

If $\{x_n\}_{n=1}^\infty$ is a monotone decreasing sequence, then
$\{-x_n\}_{n=1}^\infty$ is increasing and 
$\displaystyle{ \lim_{n \to \infty} x_n =-\lim_{n \to \infty} -x_n.}$  
\end{proof}

\begin{cor}\label{cor: limsup}
If $\{x_n\}_{n=1}^\infty$ is a bounded sequence
then 
\[
\lim_{m \to \infty} \sup \{x_n\}_{n=m}^\infty \text{\ \ \ and\ \ \ }
\lim_{m \to \infty} \inf \{x_n\}_{n=m}^\infty
\]
both exist.  We will denote them by $\displaystyle{ \limsup_{n \to \infty} x_n 
\text{ and } \liminf_{n \to \infty} x_n }$ respectively.
The sequence  $\{x_n\}_{n=1}^\infty$ has limit $L$, i.e., $\lim x_n = L,$
if and only if 
\[
\displaystyle{\liminf_{n \to \infty} x_n = \limsup_{n \to \infty} x_n} = L.
\]
\end{cor}

\begin{proof}
If $y_m = \sup \{x_n\}_{n=m}^\infty$, then $\{y_m\}_{m=1}^\infty$ is
a monotone decreasing sequence, so $\displaystyle{\lim_{m \to \infty} y_m}$
exists.  The proof that $\liminf x_n$ exists is similar.

The fact that $\inf \{x_n\}_{n=m}^\infty \le x_m \le \sup \{x_n\}_{n=m}^\infty$
implies that if 
\[
\displaystyle{\liminf_{n \to \infty} x_n = \limsup_{n \to \infty} x_n} = L
\]
then $\displaystyle{\lim_{n \to \infty} x_n}$ exists and equals $L$.
\end{proof}

\begin{defn}[Cauchy Sequence]
A sequence $\{x_n\}_{n=1}^\infty$ is called a \myindex{Cauchy sequence}
if for every $\eps >0$ there is an $N >0$ (depending on $\eps$) such
that $|x_n - x_m| < \eps$ for all $n,m \ge N.$
\end{defn}

\begin{thm}[Cauchy Sequences Have Limits]\label{thm: cauchy}
If $\{x_n\}_{n=1}^\infty$ is a Cauchy sequence, then 
$\displaystyle{\lim_{n\to \infty} x_n}$ exists.
\end{thm}

\begin{proof}
First we show that if 
$\{x_n\}_{n=1}^\infty$ is a Cauchy sequence, then it is 
bounded.  For $\eps =1$  there is an $N_1$ such that
$|x_n - x_m| < 1$ for all $n,m \ge N_1.$  Hence for any
$n \ge N_1$ we have $|x_n| \le |x_n - x_{N_1}| + |x_{N_1}| \le  |x_{N_1}| + 1.$
It follows that if $M = 1+ max \{x_n\}_{n=1}^{N_1}$, then
$|x_n| \le M$ for all $n.$
Hence $\displaystyle{\limsup_{n \to \infty} x_n}$ exists. 

Since the sequence is Cauchy, given $\eps >0$ there is an $N$ such that
that $|x_n - x_m| < \eps/2$ for all $n,m \ge N.$
Let 
\[
L = \limsup_{n\to \infty} x_n = \lim_{n \to \infty} \sup \{x_m\}_{m=n}^\infty.
\]
Hence by Proposition (\ref{prop: sup}) there is an $M \ge N$ such that
$|x_M - L| < \eps/2$.  It follows that for any $n > M$ we have
$|x_n - L| \le |x_n - x_M| + |x_M - L| < \eps/2 + \eps/2 = \eps.$
So $\displaystyle{\lim_{n \to \infty} x_n = L.}$
\end{proof}

\begin{defn}[Convergent and Absolutely Convergent]
\label{thm: absolute convergence}
\index{absolute convergence}
An infinite series $\sum_{n=1}^\infty x_n$ of real numbers is said 
to \myindex{converge} provided the sequence $\{S_m\}_{m=1}^\infty$ 
converges where $S_m = \sum_{n=1}^m x_n.$ It is said
to \myindex{converge absolutely} provided the series
$\sum_{n=1}^\infty |x_n|$ converges.
\end{defn}

\begin{thm}[Absolutely Convergent Series]
\label{thm: abs conv series}
If the series $\sum_{n=1}^\infty x_n$ converges absolutely, then it
converges.
\end{thm}
\begin{proof}
Let $S_m = \sum_{i=1}^m x_i$ be the partial sum.  We must show that
$\displaystyle{\lim_{m\to \infty} S_m}$ exists.  We will do this by
showing it is a Cauchy sequence.  Since the series 
$\sum_{i=1}^\infty |x_i|$ converges, given $\eps >0,$
there is an $N>0$ such that $\sum_{i=N}^\infty |x_i| < \eps.$
Hence if $m > n \ge N$ 
\[
|S_m - S_n| = \Big | \sum_{i=n+1}^m x_i \Big | \le
\sum_{i=n+1}^m |x_i| \le \sum_{i=N}^\infty |x_i| < \eps.
\]
Hence $\{S_n\}$ is a Cauchy sequence and converges.
\end{proof}

\section{Set Theory and Countability}\label{sec: sets}

\begin{prop}[Distributivity of $\cap$ and $\cup$]\label{prop: set distrib}
\index{distributivity, of $\cup$ and $\cap$}
If for each $j$ in some index set $J$ there is a set $B_j$  and 
$A$ is an arbitrary set, then
\[
A \cap \bigcup_{j \in J} B_j = \bigcup_{j \in J} (A \cap  B_j)
\text{ and } A \cup \bigcap_{j \in J} B_j = \bigcap_{j \in J} (A \cup  B_j).
\]
\end{prop}

The proof which is straightforward is left to the reader.

\begin{defn}[Set Difference, Complement]
We define the \myindex{set difference} of sets $A$ and $B$ by
\[
A \setminus B = \{ x\ |\ x \in A \text{ and } x \notin B\}.
\]
If all the sets under discussion are subsets of some fixed
larger set $E$, then we can define the \myindex{complement} of $A$ with
respect to $E$ to be $A^c = E \setminus A.$
\end{defn}

We will normally just speak of the complement $A^c$ of $A$ when
it is clear what the larger set $E$ is.  Note the obvious facts
that $(A^c)^c = A$ and that $A \setminus B = A \cap B^c$.

\begin{prop}\label{prop: comp intersection}
If for each $j$ in some index set $J$ there is a set $B_j \subset E$, then
\[
\bigcap_{j \in J} B_j^c = \big (\bigcup_{j \in J} B_j \big)^c
\text{ and } \bigcup_{j \in J} B_j^c = \big (\bigcap_{j \in J} B_j \big)^c.
\]
\end{prop}

Again the elementary proof is left to the reader.

\begin{prop}[Well Ordering of $\N$]\label{prop: well-order}
Every non-empty subset $A$ of $\N$ has a least element which
we will denote $min(A).$
\end{prop}
\begin{proof}  Every finite subset of $\N$ clearly has a greatest
element and a least element.  Suppose $A \subset \N$ is non-empty.
Let $B = \{ n \in \N\ |\ n < a \text{ for all } a \in A\}.$
If $1 \notin B$, then $1 \in A$ and it is the least element.
Otherwise $1 \in B$ so $B \ne \emptyset$. 
Let $b$ be the greatest element of the
finite set $B$.  The element $a_0 = b + 1$ is in $A$ and is its
least element.
\end{proof}

\begin{defn}[Injection, Surjection]\label{def: jection}
\index{injective} \index{surjective} \index{bijective}
Suppose $A$ and $B$ are sets and $\phi: A \to B$ is a
function. Then 
\begin{enumerate}
\item[(1)] The function $\phi$ is called {\em injective}
(or {\em one-to-one})
if $\phi(x) = \phi(y)$ implies $x=y.$
\item[(2)] The function $\phi$ is called {\em surjective}
(or {\em onto})
if for every $b \in B$ there exists $x \in A$ such that
$\phi(x) = b.$
\item[(3)] The function $\phi$ is called {\em bijective}
if it is both injective and surjective.
\item[(4)] If $C \subset B$ the \myindex{set inverse}
$\phi^{-1}(C)$ is defined to be 
$\{a\ |\  a \in A \text{ and } \phi(a) \in C\}$.
If $C$ consists of a single element $c$ we write
$\phi^{-1}(c)$ instead of the more cumbersome
$\phi^{-1}(\{c\})$.
\end{enumerate}
\end{defn}

The notion of countability, which we now define, turns out
to be a crucial ingredient in the concept of measure which
is the main focus of this text.

\begin{defn}[Countable]
A set $A$ is called \myindex{countable} if it is finite
or there is a bijection from $A$ to the natural numbers $\N$,
(i.e. a one-to-one correspondence between elements of $A$ and 
elements of $\N).$  A set which is not countable is
called \myindex{uncountable}.
\end{defn}

The following are standard properties of countable sets which 
we will need.

\begin{prop}[Countable Sets]\label{prop: countable}\ 
\begin{enumerate}
\item[(1)] If $A$ is countable, then any non-empty subset of $A$ is countable.
\item[(2)] A set $A$ is countable if and only if there is a 
surjective function $f: \N \to A$.
\end{enumerate}
\end{prop}

\begin{proof}
Item (1) is trivial if $A$ is finite.  Hence in proving
it we may assume
there is a bijection from $A$ to $\N$, and indeed, without loss of
generality, we may assume $A$ in fact equals $\N.$

To prove (1) suppose $B$ is a non-empty subset of $A = \N.$ 
If $B$ is finite it is countable so assume it is infinite.
Define $\phi : \N \to B$ by $\phi(1) = min(B),$ and
\[
\phi(k) = min(B \setminus \{\phi(1), \dots, \phi(k-1)\}).
\] 
The function $\phi$ is injective and defined for all $k
\in \N$.
Suppose $m \in B$ and let $c$ be the number of elements in the finite
set $\{ n \in B\ |\ n \le m\}$.  Then $\phi(c) = m$ and hence
$\phi$ is surjective.

To prove (2) suppose $f: \N \to A$ is surjective. Define
$\psi: A \to \N$ by $\psi(x) = min (f^{-1}(x)).$ This
is a bijection from $A$ to $\psi(A).$
Since $\psi(A)$ is a subset of $\N$ it is countable by (1).
This proves one direction of (2).  The converse is nearly
obvious.  If $A$ is countably infinite, then there is a 
bijection (and hence a surjection) $f: \N \to A$.  But if
$A$ is finite  one can easily define a
surjection $f: \N \to A$. 
\end{proof}

\begin{prop}[Products and Unions of Countable Sets]\label{prop: cartesian}
If $A$ and $B$ are countable, then their Cartesian product
$A \times B = \{(a,b)\ |\ a \in A,\ b \in B\}$ is a countable set.
If $A_n$ is countable for each $n \in \N$ 
then $\bigcup_{n=1}^\infty A_n$ is countable.
\end{prop}
\begin{proof}
We first observe that from part (1) of the previous proposition
a set $A$ is countable if there is an injective function 
$\phi: A \to \N$.  The function $\phi:\N \times \N \to \N$
given by $\phi(m,n) = 2^m 3^n$ is easily seen to be injective.
This is because $2^m 3^n = 2^r 3^s$ only if $2^{m-r} =  3^{s-n}.$
This is only possible if $m-r = s-n =0.$  Hence $\N \times \N$
is countable.  To show $A \times B$ is countable when $A$ and $B$
are, we note that there are surjective functions $f: \N \to A$ and
$g: \N \to B$ so
\[
f \times g: \N \times \N \to A \times B
\]
is surjective.  Since $\N \times \N$ is countable it follows
from part (2) of the previous proposition that $A \times B$
is countable.

To prove that a countable union of countable sets is countable
note that if $A_n$ is countable there is
a surjection $\psi_n : \N \to A_n$.
The function 
\[
\Psi: \N \times \N \to \bigcup_{n=1}^\infty A_n
\]
given by
\[
\Psi(n,m) = \psi_n(m)
\]
is a surjection.  Since $\N \times \N$ is countable it follows that
$\bigcup_{n=1}^\infty A_n$ is countable by part (2) of the previous
proposition.
\end{proof}

\begin{cor}[$\Q$ is countable]
The rational numbers $\Q$ are countable.
\end{cor}
\begin{proof}
The set $\Z$ is countable (see exercises below)
so $\Z \times \N$ is countable and the function
$\phi: \Z \times \N \to \Q$ given by $\phi(n,m) = n/m$
is surjective so the set of rationals $\Q$ is countable.
\end{proof}

For an arbitrary set $A$ we will denote by $\PP(A)$ its 
\myindex{power set}, which is the set of all subsets of $A$.  

\begin{prop}\label{prop: power} 
Suppose $A$ is a non-empty set and
$f: A \to \PP(A).$  Then $f$ is not surjective.
\end{prop}

\begin{proof}  This proof is short and elegant, but slightly
tricky.  For $a \in A$ either $a \in f(a)$ or 
$a \notin f(a)$.  Let $B = \{ a \in A\ |\ a \notin f(a)\}.$

Let $x$ be any element of $A$.
If $x \in B$, then, by the definition of $B$, we know
$x \notin f(x)$ and $x \in B$ so $f(x) \ne B.$
On the other hand if $x \notin B$, then by the definition
of $B$ we know $x \in f(x)$
and since $x \notin B$ we again conclude $f(x) \ne B$.
Thus in either case $f(x) \ne B$, i.e.
there is no $x$ with $f(x) = B$, so $f$ is not 
surjective.
\end{proof}

As an immediate consequence we have the existence of an
uncountable set. 

\begin{cor} 
The set $\PP(\N)$ is uncountable.
\end{cor}

\begin{proof}
This is an immediate consequence of Proposition~(\ref{prop: power})
and part (2) of Proposition~(\ref{prop: countable}).
\end{proof}

\begin{cor} \label{cor: uncountable}
If $f: A \to B$ is surjective and $B$ is uncountable, then
$A$ is uncountable.
\end{cor}

\begin{proof}
This is an immediate consequence of 
part (2) of Proposition~(\ref{prop: countable}), since if
$A$ were countable the set $B$ would also have to be 
countable.
\end{proof}

Later we will give an easy proof using measure theory
that the set of irrationals is not countable
(see Corollary (\ref{cor: reals uncountable})).
But an elementary proof of this fact is outlined in the exercises
below.

The next axiom asserts
that there is a way to pick an element from each non-empty subset of $A.$

\begin{axiom_of_choice}\label{axiom_of_choice}
\index{Axiom of Choice}
For any non-empty set $A$ there is a choice function 
\[
\phi: \PP(A) \setminus \{\emptyset\} \to A,
\]
i.e. a function such that for every non-empty subset $B \subset A$
we have $\phi(B) \in B.$
\end{axiom_of_choice}

\begin{exer}\label{exer: uncountable}\
\begin{enumerate}
\item Prove Propositions (\ref{prop: set distrib}) 
and (\ref{prop: comp intersection}).
\item ({\em Inverse Function})\\ 
If $f : A \to B$, then $g: B \to A$ is called the 
\myindex{inverse function} of $f$ provided $g(f(a)) = a$ for all $a\in A$ and
$f(g(b)) = b$ for all $b \in B.$ \subitem{(a)} Prove that if the
inverse function exists it is unique (and hence it can be referred to
as {\em the} inverse).  \subitem{(b)} Prove that $f$ has an inverse if
and only if $f$ is a bijection.  \subitem{(c)} If it exists we denote
the inverse function of $f$ by $f^{-1}$.  This is a slight abuse of
notation since we denote the set inverse (see part (4) of
Definition~(\ref{def: jection})) the same way. To justify this abuse
somewhat prove that if $f$ has an inverse $g$, then for each $b \in B$
the set inverse $f^{-1}(\{b\})$ is the set consisting of the single
element $g(b)$.  Conversely show that if for every $b \in B$ the
set inverse $f^{-1}(\{b\})$ contains a single element, then $f$ has
an inverse $g$ defined by letting $g(b)$ be that single element.

\item Prove that any subset of $\Z$ is countable by finding an
explicit bijection $f: \Z \to \N.$
\item ({\em Uncountabilitity of $\R$})\\
Let $\D$ be the set of all infinite sequences 
$d_1 d_2 d_3 \dots d_n \dots$ where each $d_n$ is either
$0$ or $1$.  
\subitem{(a)} Prove that $\D$ is uncountable.  {\em Hint: } Consider
the function $f: \PP(\N) \to \D$ defined as follows.  If
$A \subset \N$, then $f(A) = d_1 d_2 d_3 \dots d_n \dots$ 
where $d_n = 1$ if $n\in A$ and $0$ otherwise.
\subitem{(b)} Define $h: \D \to [0,1]$ by
letting $h(d_1 d_2 d_3 \dots d_n \dots)$ be the real
number whose decimal expansion is $0.d_1 d_2 d_3 \dots d_n \dots$ 
Prove that $h$ is injective.
\subitem{(c)} Prove that the closed interval $[0,1]$ is uncountable.  {\em Hint: }
Show there is a surjective function $\phi: [0,1] \to \D$
defined by $\phi(x) = h^{-1}(x)$ if $x \in h(\D)$ and
$\phi(x) = 0$ otherwise.
\subitem{(d)} Prove that if $a < b$, the closed interval
$\{x \ |\ a \le x \le b\}$, the open interval
$\{x \ |\ a < x < b\}$,
the ray $\{x \ |\ a \le x <\infty\}$, and $\R$
are all uncountable.
\end{enumerate}
\end{exer}

\section{Open and Closed Sets}

We will denote the \myindex{closed interval} $\{x \ |\ a \le x \le b\}$
by $[a,b]$ and the \myindex{open interval} $\{x \ |\ a < x < b\}$
by $(a,b).$  We will also have occasion to refer to 
the \myindex{half open} intervals $(a,b] = \{x \ |\ a < x \le b\}$
and $[a,b) = \{x \ |\ a \le  x < b\}$.
Note that the interval $[a,a]$ is the set consisting of the 
single point $a$ and $(a,a)$ is the empty set.

\begin{defn}[Open, Closed, Dense]
A subset $A \subset \R$ is called \myindex{open} if for every
$x \in A$ there is an open interval $(a,b) \subset A$ such
that $x \in (a,b).$  A subset $B \subset \R$ is called
\myindex{closed} if $\R \setminus B$ is open.  A set $A \subset \R$ is 
said to be \myindex{dense} in $\R$ if every open subset contains a
point of $A$.
\end{defn}

\begin{prop}[$\Q$ is dense in $\R$]
The rational numbers $\Q$ are a dense subset of $\R.$
\end{prop}
\begin{proof}
Let $U$ be an open subset of $\R.$ By the definition of open set there
is a non-empty interval $(a,b) \subset U$.
Choose an integer $n$ such
that $\frac{1}{n} < b-a.$ Then every point of $\R$ is in one of the
intervals $[\frac{i-1}{n}, \frac{i}{n}).$ In particular,
for some integer $i_0,$\ $\frac{i_0-1}{n} \le a < \frac{i_0}{n}.$  Since
$\frac{1}{n} < b-a$ it follows that 
\[
\frac{i_0-1}{n} \le a <
\frac{i_0}{n} \le a + \frac{1}{n} < b.
\]
Hence the rational number $i_0/n$ is in $(a,b)$ and hence in $U.$
\end{proof}

\begin{thm}\label{thm: characterize open}
An open set $U \subset \R$ is a countable union of pairwise
disjoint open intervals $\bigcup_{n=1}^\infty (a_n,b_n)$.
\end{thm}
\begin{proof}
Let $x \in U$. Define $a_x = \inf \{y\ |\ [y,x] \subset U\}$ and $b_x
= \sup \{y\ |\ [x,y] \subset U\}$ and let $U_x = (a_x,b_x).$ Then $U_x
\subset U$ but $a_x \notin U$ since otherwise for some $\eps >0,\
[a_x-\eps ,a_x+\eps] \subset U$ and hence $[a_x-\eps ,x] \subset
[a_x-\eps ,a_x+\eps] \cup [a_x+\eps ,x]\subset U$ and this would
contradict the definition of $a_x$.  Similarly $b_x \notin U$.  It
follows that if $z \in U_x$, then $a_z = a_x$ and $b_z = b_x$.
Hence if $U_z \cap U_x \ne \emptyset$, then $U_z = U_x$ or 
equivalently, if $U_z \ne U_x$, then they are disjoint.

Thus $U$ is a union of open intervals, namely the
set of all the open intervals $U_x$ for $x \in U.$  Any two
such intervals are either equal or disjoint, so the collection of
distinct intervals is pairwise disjoint.

To see that this is a countable collection observe that the
rationals $\Q$ are countable so $U \cap \Q$ is countable and
the function $\phi$ which assigns to each $r \in U \cap \Q$
the interval $U_r$ is a surjective map onto this collection.  
By Proposition (\ref{prop: countable}) this collection must
be countable.
\end{proof}

\begin{exer}\ 
\begin{enumerate}
\item Prove that the complement of a closed subset of $\R$ is open.
\item Prove that an arbitrary union of open sets is open
and an arbitrary intersection of closed sets is closed.
\item A point $x$ is called a \myindex{limit point} of
a set $S$ if every open interval containing $x$ contains
points of $S$ other than $x$.  Prove that a set $S \subset \R$
is closed if and only if it contains all its limit points.
\end{enumerate}
\end{exer}

\section{Compact Subsets of $\R$}

One of the most important concepts for analysis 
is the notion of compactness.

\begin{defn}
A closed set $X \subset \R$ is called
\myindex{compact} provided every open cover
of $X$ has a finite subcover.
\end{defn}

Less tersely, $X$ is compact if for every collection
$\V$ of open sets with the property that
\[
X \subset \bigcup_{U \in \V} U
\]
there is a finite collection $U_1,U_2,\dots U_n$ of
open sets in $\V$ such that
\[
X \subset \bigcup_{k=1}^n U_k.
\]

For our purposes the key property is that closed and
bounded subsets of $\R$ are compact.

\begin{thm}[The Heine-Borel Theorem]\label{thm: compactness}
\index{Heine-Borel Theorem}
A subset $X$ of $\R$ compact if and only if it is closed and bounded.
\end{thm}
\begin{proof}
To see that a compact set is bounded observe that $U_n = (-n,n)$
defines an open cover of any subset $X$ of $\R$.  If this cover has
a finite subcover for a set $X$, then $X \subset U_m$ for some
$m$ and hence $X$ is bounded.  To show a compact set $X$
is closed observe that if $y \notin X$, then 
$U_n = (-\infty, y - \frac{1}{n})\ \cup\ (y + \frac{1}{n}, \infty)$ defines an
open cover of $\R \setminus \{y\}$ and hence of $X$.
Since this cover of $X$ has a finite subcover  there is $m >0$
such that $X \subset U_m$.  It follows that $(y -1/m, y+1/m)$
is in the complement of $X$.  Since $y$ was an arbitrary point
of the complement of $X$, this complement is open and $X$ is
closed.

To show the converse we first consider the special case 
that $X = [a,b]$ is a closed
interval.  Let $\V$ be an open cover of $X$ and define
\[
z = \sup \{ x\in [a,b] \ |\ \text{The cover $\V$ of $[a,x]$ has a
finite subcover}\}.
\]

Our aim is to prove that $z = b$ which we do by showing that the
assumption that $z < b$ leads to a contradiction. There is an
open set $U_0 \in \V$ with $z \in U_0.$  From the definition
of open sets we know there are
points $z_0,z_1 \in U_0$ satisfying $z_0 < z < z_1$.
From the definition of $z$ the cover $\V$ of $[a,z_0]$
has a finite subcover $U_1,U_2,\dots U_n$.  Then
the finite subcover $U_0,U_1,U_2,\dots U_n$ of $\V$ is
a cover of $[a,z_1]$.  Since $z < z_1$ this is contradiction
arising from the assumption $z < b.$

For an arbitrary closed bounded set $X$ we choose $a, b \in \R$
such that $X \subset [a,b]$.  If $\V$ is any open cover of
$X$ and we define $U_0 = \R \setminus X$, then 
$\V \cup \{U_0\}$ is an open cover of $[a,b]$ which must
have a finite subcover, say $U_0,U_1,U_2,\dots U_n$.  Then
$U_1,U_2,\dots U_n$ must be a cover of $X$.
\end{proof}

There is a very important property of nested families of bounded
closed sets which we will use.

\begin{thm}[Nested Families of Compact Sets]\label{thm: nested}
If $\{A_n\}_{n=1}^\infty$ is a nested family of closed bounded subsets
of $\R$, i.e.  $A_n \subset A_{n-1}$, then 
$\cap_{n=1}^\infty A_n$ is non-empty.
\end{thm}
\begin{proof}
Let $x_n = \inf A_n$. Then $\{ x_n \}$ is a bounded
monotonic sequence so the limit $z = \lim x_n$ exists
by Theorem (\ref{thm: monotone limits})
Since $A_n$ is closed $x_n \in A_n$ and
hence $x_n \in A_m$ for all $m \le n.$ It follows that for
any $m >0$ we have $z \in A_m$, i.e. $z \in \cap_{n=1}^\infty A_n$.
\end{proof}

\begin{exer}\ 
\begin{enumerate}
\item Prove that the set $\D = \{ m/2^n\ |\ m \in \Z,\ n\in \N\}$
is dense in $\R.$ 
\item Give an example of a nested family of non-empty open intervals
$U_1 \supset U_2 \dots \supset U_n \dots$ such that 
$\cap U_n = \emptyset.$
\end{enumerate}
\end{exer}

\section{Continuous and Differentiable Functions}

\begin{defn}[Continuous and Uniformly Continuous Functions]
A function $f:\R \to \R$ is \myindex{continuous} if for every $x$ and
every $\eps >0$ there is a $\delta(x)$ (depending on $x$) such
that $|f(y) - f(x)| < \eps$ whenever $|y -x| < \delta(x).$
A function $f:\R \to \R$ is \myindex{uniformly continuous} if for 
every $\eps >0$ there is a $\delta$ (independent of $x$ and $y$) such
that $|f(y) - f(x)| < \eps$ whenever $|y -x| < \delta.$
\end{defn}

\begin{thm}\label{thm: unif continuity}
If $f$ is defined and continuous on a closed interval $[a,b]$ 
then it is uniformly continuous on that interval.
\end{thm}

\begin{proof}
Suppose $\eps >0$ is given. For any $x \in [a,b]$ 
and any positive number $\delta$ let 
$U(x,\delta) = (x - \delta, x + \delta)$
From the definition of continuity 
it follows that for each $x$ there is a $\delta(x) >0$ such that for every
$y \in U(x,\delta(x))$ we have $|f(x)-f(y)| < \eps/2.$
Therefore if $y_1$ and $y_2$ are both in $U(x, \delta(x))$ we note
\[
|f(y_1)-f(y_2)| \le |f(y_1)-f(x)| + |f(x)-f(y_2)| < 
 \frac{\eps}{2} +  \frac{\eps}{2} = \eps.
\]

The collection $\{ U(x, \delta(x)/2)\ |\ x \in [a,b]\}$ is
an open cover of the compact set $[a,b]$ so it has a
finite subcover $\{ U(x_i, \delta(x_i)/2)\ |\ 1\le i \le n\}$.
Let 
\[
\delta = \frac{1}{2} \min\{ \delta(x_i)\ |\ 1\le i \le n \}.
\]

Suppose now $y_1,y_2 \in [a,b]$ and $|y_1 - y_2| < \delta.$
Then $y_1$ is in $U(x_j, \delta(x_j)/2)$ for some $1 \le j \le n$
and 
\[
|y_2 - x_j| \le |y_2 - y_1| + |y_1 - x_j| < \delta + \frac{\delta(x_j)}{2}
\le \delta(x_j).
\]
So both $y_1$ and $y_2$ are in $U(x_j,\delta(x_j))$ and hence 
$|f(y_1)-f(y_2)|  < \eps.$
\end{proof}

We will also make use of the following result from elementary
calculus.  

\begin{thm}[Mean Value Theorem]\label{thm: mean value thm}
\index{Mean Value Theorem}
If $f$ is is differentiable on the interval $[a,b]$ 
then there is $c \in (a,b)$ such that
\[
f'(c) = \frac{f(b) - f(a)}{b-a}.
\]
\end{thm}

\begin{cor}\label{cor: $f+C$}
If $f$ and $g$ are differentiable functions on $[a,b]$ and
$f'(x) = g'(x)$ for all $x$, then there is a constant
$C$ such that $f(x) = g(x) + C.$  
\end{cor}
\begin{proof}
Let $h(x) = f(x) - g(x)$, then $h'(x) = 0$ for all $x$ and we
wish to show $h$ is constant.  But if $a_0,b_0 \in [a,b]$, then
the Mean Value Theorem says
$h(b_0) - h(a_0) = h'(c)(b_0 -a_0) = 0$ since $h'(c) = 0.$
Thus for arbitrary $a_0,b_0 \in [a,b]$ we have 
$h(b_0) = h(a_0)$ so $h$ is constant.
\end{proof}
\begin{exer}\ 
\begin{enumerate}
\item ({\em Characterization of continuity})\\
Suppose $f$ is a function $f: \R \to \R$.
\subitem{(a)} Prove that $f$ is continuous
if and only if the set inverse $f^{-1}(U)$ is open for
every open set $U \subset \R.$
\subitem{(b)} Prove that $f$ is continuous
if and only if the set inverse $f^{-1}((a,b))$ is open for
every open interval $(a,b).$
\subitem{(c)} Prove that $f$ is continuous
if and only if the set inverse $f^{-1}(C)$ is closed for
every closed set $C \subset \R.$
\end{enumerate}
\end{exer}

\section{Real Vector Spaces}
\begin{defn}[Inner Product Space]\label{def: inner prod space}
\index{inner product space}
A real vector space $\V$ is called an \myindex{inner product space} if there
is a function $\la\ ,\ \ra : \V \times \V \to \R$ 
which for any $v_1, v_2, w \in \V$ and any $a, c_1, c_2 \in \R$ satisfies:
\begin{enumerate}
\item {\bf Commutativity}: $\la v_1,v_2 \ra = \la v_2,v \ra.$
\item {\bf Bi-linearity}: $\la c_1v_1 +c_2 v_2,w \ra = c_1 \la v_1,w \ra + 
c_2 \la v_2,w \ra.$
\item {\bf Positive Definiteness}: 
$\la w,w \ra \ge 0$ with equality only if $w = 0.$
\end{enumerate}
\end{defn}

\begin{defn}[Norm]
If\ $\V$ is a real vector space with inner product $\la\ ,\ \ra$, we define
the associated \myindex{norm} $\|\  \|$ by $\|v\| = \sqrt{\la v,v \ra}.$
\end{defn}

\begin{prop}[Cauchy-Schwarz Inequality]\label{prop: cauchy schwartz}
\index{Cauchy-Schwarz Inequality}
If $(\V, \la\ ,\ \ra)$ is an inner product space and $v,w \in \V$, then 
\[
|\la v, w \ra| \le \|v\|\ \|w\|,
\]
with equality if and only if $v$ and $w$ are multiples of a single
vector.
\end{prop}
\begin{proof}
First assume $\|v\| = \|w\| =1.$ Then
\begin{align*}
\|\la v,w\ra w\|^2 + \|v - \la v,\ w\ra w\|^2  
&= \la v,w\ra^2 \|w\|^2 + \la v - \la v, w\ra w,\ v - \la v, w\ra w\ra\\
&= \la v,w\ra^2 \|w\|^2 + \|v\|^2 - 2\la v,w\ra^2 + \la v,w\ra^2\|w\|^2\\
&= \|v\|^2 = 1,
\end{align*}
since $\|v\|^2 = \|w\|^2 = 1.$  
Hence
\[
\la v,w\ra^2 = \|\la v,w\ra w\|^2 \le \|\la v,w\ra w\|^2 
+ \|v - \la v,w\ra w\|^2  = 1
\]
with equality only if $\|v - \la v,w\ra w\| = 0$ or $v = \la v,w\ra
w$.  This implies the inequality $|\la v, w \ra| \le 1 = \|v\|\
\|w\|$, when $v$ and $w$ are unit vectors.  The result is trivial if
either $v$ or $w$ is $0$.  Hence we may assume the vectors are
non-zero multiples $v = av_0$ and $w = bw_0$ of unit vectors $v_0$ and
$w_0$.  In this case we have $|\la v, w \ra| 
= |\la av_0, bw_0 \ra| = |ab| |\la v_0, w_0 \ra|
\le |ab| = \|av_0\|\ \|bw_0\| = \|v\|\ \|w\|$,

Observe that we have equality only if $v = \la v,w\ra w$,
i.e. only if one of the vectors is a multiple of the other.
\end{proof}

\begin{prop}[Normed Linear Space]\label{prop: norm props}
\index{normed linear space}
If $\V$ is an  inner product space and $\|\ \|$
is the norm defined by $\|v \| = \sqrt{\la v,v \ra},$
then
\begin{enumerate}
\item[(1)] For all $a \in \R$ and $v \in \V,\ \| av\| = |a|\| v\|.$
\item[(2)] For all $v \in \V,\ \| v\| \ge 0$ with equality only if
$v = 0.$
\item[(3)] {\em Triangle Inequality}: \index{triangle inequality} For all 
$v,w \in \V,\ \| v+w \| \le \| v\| + \| w\|.$
\item[(4)] {\em Parallelogram Law}: \index{parallelogram law} For all 
$v,w \in \V,$
\[
\| v - w\|^2 + \| v + w\|^2 = 2\| v\|^2 + 2\| w\|^2.
\]
\end{enumerate}
\end{prop}

\begin{proof}
The first two of these properties follow immediately from
the definition of inner product.
To prove item (3), the triangle inequality, observe
\begin{align*}
\|v+w\|^2 &= \la v+w, v+w \ra \\
&= \la v, v \ra + 2 \la v,w \ra + \la w,w \ra\\
&= \|v\|^2 + 2 \la v,w \ra + \|w\|^2\\
&\le \|v\|^2 + 2 |\la v,w\ra| + \|w\|^2\\
&\le \|v\|^2 + 2 \|v\|\ \|w\| + \|w\|^2 \hfil \text{ by Cauchy-Schwarz,}\\
&= (\|v\|+ \|w\|)^2
\end{align*}

To prove item (4), the parallelogram law, note 
$\| v - w\|^2 = \la v -w, v-w\ra =  \| v\|^2 - 2\la v,w\ra +\| w\|^2.$
Likewise $\| v+w\|^2 = \la v+w, v+w\ra =  \| v\|^2 + 2\la v,w\ra +\| w\|^2.$
Hence the sum $\| v - w\|^2 + \| v + w\|^2$ equals
$2\| v\|^2 + 2\| w\|^2.$
\end{proof}

\vfill\eject
\chapter{The Regulated and Riemann Integrals}
\section{Introduction}
We will consider several different approaches to defining the 
definite integral
\[
\int_a^b f(x) \ dx
\]
of a function $f(x)$.  These definitions will all assign the same
value to the definite integral, but they differ in the size of the
collection of functions for which they are defined.  For example,
we might try to evaluate the Riemann integral 
(the ordinary integral of beginning
calculus) of the function
\[
f(x) = 
\begin{cases}
	0, &\text{if $x$ is rational;} \\
	1, &\text{otherwise.} 
\end{cases}
\]
The Riemann integral $\int_0^1 f(x) \ dx$ is, as we will see, undefined.
But the Lebesgue integral, which we will develop, has no difficulty
with $f(x)$ and indeed $\int_0^1 f(x) \ dx = 1.$

There are several properties which we want an integral to satisfy no
matter how we define it. It is worth enumerating them at the
beginning.  We will need to check them for our different definitions.

\section{Basic Properties of an Integral}\label{sec: basic properties}

We will consider the value of the integral of functions in various
collections.  These collections all have a common domain which, for
our purposes, is a closed interval.  They also
are closed under the operations
of addition and scalar multiplication.  We will call such a collection a
\myindex{vector space} of functions.  More precisely a
non-empty set of real valued functions $\V$ defined on a fixed closed 
interval will be called a {\bf vector
space of functions} provided:
\begin{enumerate}
\item If $f,g \in \V$, then $f + g \in \V.$
\item If $f\in \V$ and $r \in \R$, then $rf \in \V.$
\end{enumerate}
Notice that this implies that the constant function $0$ is in $\V$.
All of the vector spaces we consider will contain all the constant
functions.

Three simple examples of vector spaces of functions defined on some
closed interval $I$ are the constant functions, the polynomial
functions, and the continuous functions.

An ``integral'' defined on a vector space of functions $\V$ is a 
way to assign a real number to each function in $\V$ and each subinterval
of $\I.$  For the function $f \in \V$ and the subinterval $[a,b]$
we denote this value by $\int_a^b f(x) \ dx$ and call it 
``the integral of $f$ from $a$ to $b$.''

All the integrals we consider will satisfy five basic properties
which we now enumerate.

\begin{description}
\item[I. Linearity:]
For any functions $f,g \in \V$, any $a,b \in \I$,  
and any real numbers $c_1, c_2,$
\[
\int_a^b c_1f(x) +  c_2 g(x) \ dx 
= c_1 \int_a^b f(x)\ dx +  c_2 \int_a^b g(x)\ dx.
\]
\item[II. Monotonicity:]
If functions $f,g \in \V$ satisfy $f(x) \ge g(x)$ for all $x$ and
$a,b \in \I$ satisfy $a \le b$, then
\[
\int_a^b f(x) \ dx \ge \int_a^b g(x) \ dx.
\]
In particular if $f(x) \ge 0$ for all $x$ and $a \le b$ 
then $\int_a^b f(x) \ dx \ge 0.$
\item[III. Additivity:]
For any function $f \in \V,$ and any $a,b,c \in \I,$
\[
\int_a^c f(x) \ dx = \int_a^b f(x)\ dx +  \int_b^c f(x)\ dx .
\]
In particular we allow $a, b$ and $c$ to occur in any order on the
line and we note that two easy consequences
of additivity are 
\[
\int_a^a f(x) \ dx = 0 \text{ and }
\int_a^b f(x) \ dx = -\int_b^a f(x) \ dx.
\]

\item[IV. Constant functions:]
The integral of a constant function $f(x) = C$ should 
be given by 
\[
\int_a^b C \ dx = C(b-a).
\]
If $C > 0$ and $a < b$ this just says the integral of 
$f$ is the area of the rectangle under its graph.

\item[V. Finite Sets Don't Matter:]  If $f$ and $g$ are functions
in $\V$ with $f(x) = g(x)$ for all $x$ except possibly a
finite set, then for all $a,b \in \I$
\[
\int_a^b f(x) \ dx = \int_a^b g(x) \ dx.
\]
\end{description}

Properties III, IV and V are not valid for all
mathematically interesting theories of integration.  Nevertheless,
they hold for all the integrals we will consider so we include them
in our list of basic properties.  It is important to note that these
are {\it assumptions}, however, and there are many mathematically
interesting theories where they do not hold.

There is one additional property which we will need.  It differs
from the earlier ones in that we can {\em prove} that it holds
whenever the properties above are satisfied.

\begin{prop}[Absolute Value]\label{prop: abs value}
Suppose the integral $\int_a^b f(x) \ dx$  has been defined for all
$f$ in some vector space of functions $\V$  and for all 
subintervals $[a,b]$ of $I$.
And suppose this integral satisfies properties I and II above.  Then
for any function $f \in \V$ for which $|f| \in \V$
\[
\Big | \int_a^b f(x) \ dx \Big | \le \int_a^b |f(x)| \ dx,
\]
for all $a < b$ in $\I.$
\end{prop}

\begin{proof}
This follows from monotonicity and linearity.
Since $f(x) \le |f(x)|$ for all $x$ we know
$\int_a^b f(x) \ dx \le \int_a^b |f(x)| \ dx.$
Likewise $-f(x) \le |f(x)|$ so 
$-\int_a^b f(x) \ dx = \int_a^b -f(x) \ dx \le \int_a^b |f(x)| \ dx.$
But $|\int_a^b f(x) \ dx|$ is either equal to 
$\int_a^b f(x) \ dx$ or to $-\int_a^b f(x) \ dx.$  In either
case $\int_a^b |f(x)| \ dx$ is greater so
$|\int_a^b f(x) \ dx| \le \int_a^b |f(x)| \ dx.$
\end{proof}

\section{Step Functions and the Regulated Integral}

The easiest functions to integrate are {\it step functions} which we
now define.

\begin{defn}[Step Function]
A function $f: [a,b] \to \R$ is called a \myindex{step function}
provided there numbers $x_0=a < x_1 < x_2 < \dots < x_{n-1}< x_n = b$ such that
$f(x)$ is constant on each of the open intervals $(x_{i-1}, x_{i}).$
\end{defn}

We will say that the points $x_0=a < x_1 < x_2 < \dots < x_{n-1}< x_n
= b$ define an \myindex{interval partition} 
for the step function $f$.  Note that
the definition says that on the open intervals $(x_{i-1}, x_{i})$ of
the partition $f$ has a constant value, say $c_i$, but it says nothing
about the values at the endpoints. The value of $f$ at the points
$x_{i-1}$ and $x_{i}$ may or may or may not be equal to $c_i.$ Of
course when we define the integral this won't matter because the
endpoints form a finite set.

Since the area under the graph of a positive step function is a finite
union of rectangles, it is pretty obvious what the integral should be.
The $i^{th}$ of these rectangles has width $(x_{i} - x_{i-1})$ and
height $c_i$ so we should sum up the areas $c_i (x_{i} - x_{i-1}).$ Of
course if some of the $c_i$
are negative, then the corresponding $c_i (x_{i} - x_{i-1})$ are also
negative, but that is appropriate since the area between the graph and
the $x$-axis is below the $x$-axis on the interval $(x_{i-1}, x_{i}).$

\begin{defn}[Integral of a step function] \label{def: stepfcn}
Suppose $f(x)$ is a step function with partition 
$x_0=a < x_1 < x_2 < \dots < x_{n-1}< x_n = b$
and suppose
$f(x) = c_i$ for $x_{i-1} < x < x_i.$  Then we define
\[
\int_a^b f(x) \ dx = \sum_{i=1}^n c_i (x_{i} - x_{i-1})
\]
\end{defn}

\begin{exer}\ 
\begin{enumerate}
\item Prove that the collection of all step functions on a closed
interval $[a,b]$ is a vector space which contains the constant functions.

\item Prove that if $x_0=a < x_1 < x_2 < \dots < x_{n-1}< x_n = b$
is a partition for a step function $f$ with value $c_i$ on 
$(x_{i-1},x_i)$ and 
$y_0=a < y_1 < y_2 < \dots < y_{n-1}< y_m = b$ is another partition
for the same step function with value $d_j$ on  $(y_{j-1},y_i)$, then
\[
\sum_{i=1}^n c_i (x_{i} - x_{i-1}) = \sum_{j=1}^m d_i (y_{j} - y_{j-1}).
\]
In other words the value of the integral of a step function depends
only on the function, not on the choice of partition.  {\em Hint:}
the union of the sets of points defining the two partitions defines a
third partition and the integral using this partition is equal to
the integral using each of the partitions.

\item Prove that the integral of step functions as given in Definition
\ref{def: stepfcn} satisfies properties I-V of \S\ref{sec: basic properties}.
\end{enumerate}
\end{exer}

We made the ``obvious'' definition for the integral of a step
function, but in fact, we had absolutely no choice in the matter if we
want the integral to satisfy properties I-V above.

\begin{thm}
The integral as given in Definition \ref{def: stepfcn}
is the unique real valued function defined on step functions 
which satisfies properties I-V of \S\ref{sec: basic properties}.
\end{thm}

\begin{proof}
Suppose that there is another ``integral'' defined on step functions
and satisfying I-V.  We will denote this alternate integral as 
\[
\oint_a^b f(x) \ dx.
\]
What we must show is that for every step function $f(x),$
\[
\oint_a^b f(x) \ dx = \int_a^b f(x) \ dx.
\]

Suppose that $f$ has partition 
$x_0=a < x_1 < x_2 < \dots < x_{n-1}< x_n = b$ and satisfies
$f(x) = c_i$ for $x_{i-1} < x < x_{i}.$  

Then from the additivity property 
\begin{equation}\label{eqn:1}
\oint_a^b f(x)\ dx = \sum_{i=1}^n\oint_{x_{i-1}}^{x_i} f(x)\ dx.
\end{equation}
But on the interval $[x_{i-1}, x_i]$ the function $f(x)$ is 
equal to the constant function with
value $c_i$ except at the endpoints.  Since functions
which are equal except at a finite set of points have the 
same integral, the integral of $f$ is the same as the integral
of $c_i$ on $[x_{i-1}, x_i]$.  Combining this with the
constant function property we get
\[
\oint_{x_{i-1}}^{x_i} f(x) \ dx = \oint_{x_{i-1}}^{x_i} c_i  \ dx 
= c_i (x_{i} - x_{i-1}).
\]

If we plug this value into equation (\ref{eqn:1}) we obtain
\[
\oint_a^b f(x)\ dx  = \sum_{i=1}^n c_i (x_{i} - x_{i-1}) = \int_a^b f(x) \ dx.
\]
\end{proof}

Recall the definition of a uniformly converging sequence of functions.

\begin{defn}[Uniform Convergence]
A sequence of functions $\{f_m\}$ is said to \myindex{converge uniformly}
on $[a,b]$ to a function $f$ if for every $\eps >0$ there is an
$M$ ({\em independent of $x$}) such that for all $x \in [a,b]$
\[
|f(x) - f_m(x)| < \eps \text{ whenever } m \ge M.
\]
\end{defn}
Contrast this with the following.
\begin{defn}[Pointwise Convergence]
A sequence of functions $\{f_m\}$ is said to \myindex{converge pointwise}
on $[a,b]$ to a function $f$ if for each $\eps >0$ and each
$x \in [a,b]$ there is an
$M_x$ ({\em depending on $x$}) such that
\[
|f(x) - f_m(x)| < \eps \text{ whenever } m \ge M_x.
\]
\end{defn}

\begin{defn}[Regulated Function]
\index{regulated function}
A function $f: [a,b] \to \R$ is called {\em regulated} provided there
is a sequence $\{f_m\}$ of step functions which converges {\em uniformly}
to $f$.
\end{defn}

\begin{exer}\ 
\begin{enumerate}
\item Prove that the collection of all regulated functions on a closed
interval $\I$ is a vector space which contains the constant functions.

\item Give an example of a sequence of step functions which
converge uniformly to $f(x) = x$ on $[0,1].$  Give an example of a
sequence of step functions which converge {\em pointwise} to $0$
on $[0,1]$, but which do not converge uniformly. 
\end{enumerate}
\end{exer}

Every regulated function can be uniformly approximated  as closely
as we wish by a step function.  Since we know how to integrate
step functions it is natural to take a sequence of better and better
step function approximations to a regulated function $f(x)$ and 
define the integral of $f$ to be the limit of the integrals of the
approximating step functions.  For this to work we need to know
that the limit exists and that it does not depend on the choice of
approximating step functions.

\begin{thm}
Suppose $\{f_m\}$ is a sequence of step functions on $[a,b]$ 
converging uniformly to a regulated function $f.$  Then
the sequence of numbers $\{ \int_a^b f_m(x)\ dx\}$ converges.
Moreover if $\{g_m\}$ is another sequence of step functions
which also converges uniformly to $f$, then
\[
\lim_{m \to \infty} \int_a^b f_m(x) \ dx 
= \lim_{m \to \infty} \int_a^b g_m(x) \ dx.
\]
\end{thm}

\begin{proof}
Let $z_m = \int_a^b f_m(x)\ dx$.  We will show that the sequence
$\{z_m\}$ is a Cauchy sequence and hence has a limit.
To show this sequence is Cauchy we must show that for any
$\eps > 0$ there is an $M$ such that 
$|z_p - z_q| \le \eps$ whenever $p,q \ge M.$

If we are given $\eps >0$, since $\{f_m\}$ is a sequence of
step functions on $[a,b]$ converging uniformly to $f,$ there is an $M$
such that for all $x$ 
\[
|f(x) - f_m(x)| < \frac{\eps}{2(b-a)} \text{ whenever } m \ge M.
\]
Hence whenever, $p,q \ge M$
\begin{equation}\label{eqn:2}
|f_p(x) - f_q(x)| < |f_p(x) - f(x)| + |f(x) - f_q(x)| < 
2 \frac{\eps}{2(b-a)} = \frac{\eps}{b-a}.
\end{equation}
Therefore, whenever $p,q \ge M$
\[
|z_p - z_q| = \Big |\int_a^b f_p(x) - f_q(x)\ dx \Big | 
\le \int_a^b |f_p(x) - f_q(x)|\ dx
\le \int_a^b \frac{\eps}{b-a}\ dx = \eps,
\]
where the first inequality comes from the absolute value property
of Proposition \ref{prop: abs value} and the 
second follows from the monotonicity property and equation (\ref{eqn:2}).
This shows that the sequence $\{z_m\}$ is Cauchy and hence converges.

Now suppose that $\{g_m\}$ is another sequence of step functions
which also converges uniformly to $f$, then for any $\eps >0$
there is an $M$ such that for all $x$ 
\[
|f(x) - f_m(x)| < \eps \text{ and } |f(x) - g_m(x)|< \eps
\]
whenever $m \ge M.$  It follows that 
\[
|f_m(x) - g_m(x)| \le |f_m(x) - f(x)| + |f(x) - g_m(x)|< 2\eps.
\]
Hence, using the absolute value and monotonicity properties, we see
\[
\Big |\int_a^b f_m(x) - g_m(x)\ dx \Big | 
\le \int_a^b |f_m(x) - g_m(x)|\ dx
\le \int_a^b 2\eps\ dx = 2\eps (b-a),
\]
for all $m \ge M.$  Since $\eps$ is arbitrarily small we may
conclude that 
\[
\lim_{m \to \infty} \Big |\int_a^b f_m(x)\  dx - \int_a^b g_m(x)\ dx \Big |
= \lim_{m \to \infty} \Big |\int_a^b f_m(x) - g_m(x)\ dx \Big | = 0.
\]
This implies 
\[
\lim_{m \to \infty} \int_a^b f_m(x) = 
\lim_{m \to \infty} \int_a^b g_m(x)\ dx.
\]
\end{proof}

This result enables us to define the regulated integral.

\begin{defn}[The Regulated Integral] \label{def: reg int}
\index{regulated integral}
If $f$ is a regulated
function on $[a,b]$ we define the \myindex{regulated integral} by
\[
\int_a^b f(x) \ dx = \lim_{n \to \infty} \int_a^b f_n(x) \ dx
\]
where $\{f_n\}$ is any sequence of step functions converging
uniformly to $f.$
\end{defn}

We next need to see that the regulated functions form a large
class including all continuous functions.

\begin{thm}[Continuous functions are regulated] 
Every continuous function $f: [a,b] \to \R$ is a regulated function.
\end{thm}

\begin{proof}
By Theorem~\ref{thm: unif continuity} 
a continuous function $f(x)$ defined on a closed interval $[a,b]$ is
\myindex{uniformly continuous}.  That is, given $\eps > 0$ there is
a corresponding $\delta >0$ such that
$|f(x) - f(y)| < \eps$ whenever $|x - y| < \delta.$
Let $\eps_n = 1/2^n$ and let $\delta_n$ be the corresponding
$\delta$ guaranteed by uniform continuity.

Fix a value of $n$ and choose a partition $x_0 =a < x_1 < x_2 < \dots
< x_m = b$ with $x_{i} - x_{i-1} < \delta_n.$  For example, we could 
choose $m$ so large that if we define $\Delta x = (b-a)/m$
then $\Delta x < \delta_n$ and, then we could define 
$x_i$ to be $a + i \Delta x$.  Next we define a step function
$f_n$ by
\[
f_n(x) = f(x_{i}) \text{ for all $x \in [x_{i-1}, x_{i}[.$}
\]
That is, on each half open interval $[x_{-1}i, x_{i}[$ we define
$f_n$ to be the constant function whose value is the value of
$f$ at the left endpoint of the interval.  The value of $f_n(b)$
is defined to be $f(b).$

Clearly $f_n(x)$ is a step function with the given partition.  We must
estimate its distance from $f$.  Let $x$ be an arbitrary point of
$[a,b].$ It must lie in one of the open intervals of the partition or
be an endpoint of one of them;
say $x \in [x_{i-1}, x_{i}[.$ Then since $f_n(x) = f_n(x_{i-1}) =
f(x_{i-1})$ we may conclude
\[
|f(x) - f_n(x)| \le |f(x) - f(x_{i-1})| < \eps_n
\]
because of the uniform continuity of $f$ and the fact that
$|x - x_{i-1}| < \delta_n.$

Thus we have constructed a step function $f_n$ with the property
that for all $x \in [a,b]$ 
\[
|f(x) - f_n(x)|  < \eps_n.
\]
So the sequence $\{ f_n \}$ converges uniformly to $f$ and $f$
is a regulated function. 
\end{proof}

\begin{exer}\ 
\begin{enumerate}
\item Give an example of a continuous function on the {\em open}
interval $(0,1)$ which is not regulated, i.e. which cannot be 
uniformly approximated by step functions.

\item Prove that the regulated integral, as given in
(\ref{def: reg int}), satisfies properties I-V
of \S\ref{sec: basic properties}.

\item Prove that $f$ is a regulated function on $\I = [a,b]$ if and
only if both the limits
\[
\lim_{x \to c+} f(x) \text{ and } \lim_{x \to c-} f(x)
\]
exist for every $c \in (a,b).$ (See section VII.6 of Dieudonn\'e \cite{D}).
\end{enumerate}
\end{exer}

\section{The Fundamental Theorem of Calculus}

The most important theorem of elementary calculus asserts that if
$f$ is a continuous function on $[a,b]$, then its integral 
$\int_a^b f(x) \ dx$ can be evaluated by finding an anti-derivative.
More precisely, if $F(x)$ is an anti-derivative of $f$, then
\[
\int_a^b f(x) \ dx = F(b) - F(a).
\]
We now can present a rigorous proof of this result.  We will actually
formulate the result slightly differently and show that the result
above follows easily from that formulation.

\begin{thm}\label{thm: fund}
If $f$ is a continuous function and we define
\[
F(x) = \int_a^x f(t) \ dt
\]
then $F$ is a differentiable function and $F'(x) = f(x).$
\end{thm}
\begin{proof}
By definition 
\[
F'(x_0) = \lim_{h \to 0} \frac{F(x_0 + h) - F(x_0)}{h};
\]
so we need to show that 
\[
\lim_{h \to 0} \frac{F(x_0 + h) - F(x_0)}{h} =  f(x_0).
\]
or equivalently
\[
\lim_{h \to 0} \Big |\frac{F(x_0 + h) - F(x_0)}{h} - f(x_0) \Big | = 0.
\]
To do this we note that
\begin{align}
\Big | \frac{F(x_0 + h) - F(x_0)}{h} - f(x_0) \Big |
&= \Big |\frac{\int_{x_0}^{x_0 + h} f(t)\ dt}{h} - f(x_0)\Big | \notag \\
& = \Big |\frac{\int_{x_0}^{x_0 + h} f(t)\ dt - f(x_0)h}{h} \Big |\notag \\
& = \frac{\big |\int_{x_0}^{x_0 + h} (f(t) - f(x_0))\ dt\big |}{|h|}
\label{eqn:fund 1}
\end{align}

Monotonicity tells us that when $h$ is positive
\[
\Big |\int_{x_0}^{x_0 + h} (f(t) - f(x_0))\ dt \Big | \le \int_{x_0}^{x_0 + h} 
|f(t) - f(x_0)|\ dt
\] and if $h$  is negative 
\[
\Big |\int_{x_0}^{x_0 + h} (f(t) - f(x_0))\ dt \Big | 
\le \int_{x_0+h}^{x_0} |f(t) - f(x_0)|\ dt 
= - \int_{x_0}^{x_0+h} |f(t) - f(x_0)|\ dt.
\]
In either case we see
\[
\Big |\int_{x_0}^{x_0 + h} (f(t) - f(x_0))\ dt \Big |
\le \Big |\int_{x_0}^{x_0 + h} |f(t) - f(x_0)|\ dt \Big |
\]
Combining this with inequality (\ref{eqn:fund 1}) above we obtain
\begin{equation}\label{eqn:fund 2}
\Big | \frac{F(x_0 + h) - F(x_0)}{h} - f(x_0) \Big |
\le \frac{\big |\int_{x_0}^{x_0 + h} |f(t) - f(x_0)|\ dt \big |}{|h|}.
\end{equation}

But the continuity of $f$ implies that given $x_0$ and any $\eps >0$
there is a $\delta >0$ such that whenever $|t - x_0| < \delta$
we have $|f(t) - f(x_0)| < \eps$.  Thus if $|h| < \delta$
then $|f(t) - f(x_0)| < \eps$ for 
all $t$ between $x_0$ and $x_0 + h$.  It follows that
$\big |\int_{x_0}^{x_0 + h} |f(t) - f(x_0)|\ dt\big | < \eps |h|$ and hence
that 
\[
\frac{\big |\int_{x_0}^{x_0 + h} |f(t) - f(x_0)|\ dt\big |}{|h|} < \eps.
\]
Putting this together with equation (\ref{eqn:fund 2}) above we have
that 
\[
\Big |\frac{F(x_0 + h) - F(x_0)}{h} - f(x_0)\Big | < \eps
\]
whenever $|h| < \delta$ which was exactly what we needed to show.
\end{proof}

\begin{cor}
{\bf Fundamental Theorem of Calculus.} 
If $f$ is a continuous function on $[a,b]$ and $F$ is any anti-derivative of 
$f$, then 
\[
\int_a^b f(x) \ dx = F(b) - F(a)
\]
\end{cor}

\begin{proof}
Define the function $G(x) = \int_a^x f(t) \ dt.$  By Theorem
(\ref{thm: fund}) the derivative of $G(x)$ is $f(x)$ which is 
also the derivative of $F$.  Hence
$F$ and $G$ differ by a constant, say $F(x) = G(x) + C$
(see Corollary (\ref{cor: $f+C$})).

Then 
\begin{align*}
F(b) - F(a) &= (G(b) + C)  - (G(a) + C)\\
& = G(b) - G(a) \\
&= \int_a^b f(x) \ dx - \int_a^a f(x) \ dx \\
&=  \int_a^b f(x) \ dx.
\end{align*}
\end{proof}

\begin{exer}\label{exer: integral continuous}\ 
\begin{enumerate}
\item Prove that if $f : [a,b] \to \R$ is a regulated function
and $F:[a,b] \to \R$ is defined to by $F(x) = \int_a^x f(t) \ dt$
then $F$ is continuous.
\end{enumerate}
\end{exer}
\section{The Riemann Integral}

We can obtain a larger class of functions for which a good integral 
can be defined by using a different method of comparing with step
functions.

Suppose that $f(x)$ is a bounded function on the interval ${\mathcal
I} = [a,b]$ and that it is an element of a vector space of functions
which contains the step functions and for which there is an integral
defined satisfying properties I-V of \S\ref{sec: basic properties}.
If $u(x)$ is a step function satisfying $f(x) \le u(x)$ for all $x \in
\I$, then monotonicity implies that if we can define $\int_a^b
f(x)\ dx$ it must satisfy $\int_a^b f(x)\ dx \le \int_a^b u(x)\ dx.$

This is true for every step function $u$ satisfying $f(x) \le u(x)$
for all $x \in \I.$
Let $\U(f)$ denote the set of all step functions with this property.
Then if we can define
$\int_a^b f(x)\ dx$ in a way that satisfies monotonicity it must 
also satisfy
\begin{equation}\label{eqn:inf}
\int_a^b f(x)\ dx \le  \inf \Big \{ \int_a^b u(x)\ dx\ 
\Big |\   u \in \U(f) \Big \}. 
\end{equation}
The {\it infimum} exists because all the step functions in $\U(f)$
are bounded below by a lower bound for the function $f.$

Similarly we define $\yL(f)$ to be the set of all step functions
$v(x)$ such that $v(x) \le f(x)$ for all $x \in \I.$  Again 
if we can define
$\int_a^b f(x)\ dx$ in such a way that it satisfies monotonicity
it must also satisfy
\begin{equation}\label{eqn:sup}
\sup \Big \{ \int_a^b v(x)\ dx \ \Big |\   v \in \yL(f)\ \Big \} 
\le \int_a^b f(x)\ dx.
\end{equation}
The {\it supremum} exists because all the step functions in $\U(f)$
are bounded above by an upper bound for the function $f.$

Putting inequalities (\ref{eqn:inf}) and (\ref{eqn:sup}) together,
we see if $\V$ is any vector space of
bounded functions which contains the step functions and we manage to
define the integral of functions in $\V$ in a way that satisfies
monotonicity, then this integral must satisfy
\begin{equation}\label{eqn:Riem}
\sup \Big \{ \int_a^b v(x)\ dx \ \Big |\   v \in \yL(f)\Big \} 
\le \int_a^b f(x)\ dx 
\le \inf \Big \{ \int_a^b u(x)\ dx\ \Big |\   u \in \U(f)\Big \}
\end{equation}
for every $f \in \V.$  Even if we can't define an integral for $f$,
however, we still have the inequalities of the ends.

\begin{prop}\label{prop: ineq}
Let $f$ be any bounded function on the interval $\I = [a.b].$ 
Let $\U(f)$ denote the set of all step functions $u(x)$ on $\I$ such
that $f(x) \le u(x)$ for all $x$ and let $\yL(f)$ denote the set of all step
functions $v(x)$ such that $v(x) \le f(x)$ for all $x.$
Then
\[
\sup \Big \{ \int_a^b v(x)\ dx \ \Big |\   v \in \yL(f)\Big \} 
\le \inf \Big \{ \int_a^b u(x)\ dx\ \Big |\   u \in \U(f)\Big \}.
\]
\end{prop}
\begin{proof}
If $v \in \yL(f)$ and $u \in \U(f)$, then $v(x) \le f(x) \le u(x)$ for
all $x \in \I$ so monotonicity implies that 
$\int_a^b v(x)\ dx \le \int_a^b u(x)\ dx$. 
Hence if 
\[
V =  \Big \{ \int_a^b v(x)\ dx \ \Big |\   v \in \yL(f)\Big \} 
\text{ and\ \ \  }
U =  \Big \{ \int_a^b u(x)\ dx\ \Big |\   u \in \U(f)\Big \}
\]
then every number in the set $V$ is less than or equal to every number
in the set $U$.  Thus $\sup V \le \inf U$ as claimed
\end{proof}

It is not difficult to see that sometimes the two ends of
this inequality are not equal (see Exercise (\ref{exer: non-integrable}) below), 
but if it should happen that
\[
\sup \Big \{ \int_a^b v(x)\ dx \ \Big |\   v \in \yL(f)\Big \} 
= \inf \Big \{ \int_a^b u(x)\ dx\ \Big |\   u \in \U(f)\Big \}.
\]
then we have only one choice for $\int_a^b f(x)\ dx$; it must be this
common value.

This motivates the definition of the next vector space of functions we can
integrate.  Henceforth we will use the more compact notation
\[
\sup_{v \in \yL(f)} \Big \{\int_a^b v(x)\ dx \Big\} \text{\ \ instead of\ \ } 
\sup \Big \{ \int_a^b v(x)\ dx \ \Big |\   v \in \yL(f)\Big \}
\] and 
\[
\inf_{u \in \U(f)} \Big \{\int_a^b u(x)\ dx \Big\} \text{\ \ instead of\ \ } 
\inf \Big \{ \int_a^b u(x)\ dx\ \Big |\   u \in \U(f)\Big \}.
\]

\begin{defn} {\bf The Riemann Integral.}  
\index{Riemann integral}
Suppose $f$ is a bounded function on the interval $\I = [a,b].$ 
Let $\U(f)$ denote the set of all step functions $u(x)$ on $\I$ such
that $f(x) \le u(x)$ for all $x$ and let $\yL(f)$ denote the set of all step
functions $v(x)$ such that $v(x) \le f(x)$ for all $x.$
The function $f$ is said to be Riemann integrable provided
\[
\sup_{v \in \yL(f)} \Big \{\int_a^b v(x)\ dx \Big\}
= \inf_{u \in \U(f)} \Big \{\int_a^b u(x)\ dx \Big\}.
\]
In this case its Riemann integral $\int_a^b f(x) \ dx$ is defined
to be this common value.
\end{defn}

\begin{thm}\label{thm: riem criterion}
A bounded function $f:[a,b] \to \R$ is Riemann integrable if and
only if, for every $\eps >0$ there are step functions $v_0$ and
$u_0$ such that $v_0(x) \le f(x) \le u_0(x)$ for all $x \in [a,b]$
and 
\[
\int_a^b u_0(x) \ dx - \int_a^b v_0(x) \ dx \le \eps.
\]
\end{thm}
\begin{proof}
Suppose the functions $v_0 \in \yL(f)$ and $u_0 \in
\U(f)$ have integrals within $\eps$ of each other. Then
\[
\int_a^b v_0(x) \ dx \le
\sup_{v \in \yL(f)} \Big \{\int_a^b v(x)\ dx \Big\}\le 
\inf_{u \in \U(f)} \Big \{\int_a^b u(x)\ dx \Big\}
\le \int_a^b u_0(x) \ dx.
\]
This implies 
\[
\inf_{u \in \U(f)} \Big \{\int_a^b u(x)\ dx \Big\}-
\sup_{v \in \yL(f)} \Big \{\int_a^b v(x)\ dx \Big\}\le \eps.
\]
Since this is true for all $\eps >0$ we conclude that $f$ is Riemann
integrable.

Conversely if $f$ is Riemann integrable, then 
by Proposition \ref{prop: sup} there exists a step function
$u_0 \in \U(f)$ such that
\[
\int_a^b u_0(x) \ dx - \int_a^b f(x)\ dx
= \int_a^b u_0(x) \ dx - \inf_{u \in \U(f)} \Big \{\int_a^b u(x)\ dx \Big\}< \eps/2.
\]
Similarly there exists a step function
$v_0 \in \yL(f)$ such that
\[
\int_a^b f(x)\ dx - \int_a^b v_0(x) \ dx < \eps/2.
\]
Hence 
\[
\int_a^b u_0(x) \ dx - \int_a^b v_0(x) \ dx < \eps/2 + \eps/2 = \eps.
\]
and $u_0$ and $v_0$ are the desired functions.
\end{proof}

\begin{exer}\label{exer: non-integrable}\ 
\begin{enumerate}
\item At the beginning of these notes we mentioned the function
$f: [0,1] \to \R$ which has the value $f(x) = 0$ if $x$ is rational
and $1$ otherwise.  Prove that for this function
\[
\sup_{v \in \yL(f)} \Big \{\int_0^1 v(x)\ dx \Big\}= 0  \text{ and }
\inf_{u \in \U(f)} \Big \{ \int_0^1 u(x)\ dx \Big\}= 1.
\]
Hence $f$ is not Riemann integrable.
\end{enumerate}
\end{exer}

There are several facts about the relation with the regulated integral
we must establish.  Every regulated function is Riemann integrable,
but there are Riemann integrable functions which have no regulated
integral.  Whenever a function has both types of integral the values
agree.  We start by giving an example of a function which is Riemann
integrable, but not regulated.

\begin{ex}
Define the function $f: [0,1] \to \R$ by
\[
f(x) = 
\begin{cases}
	1, &\text{ if } x = \frac{1}{n} \text{ for  $n \in \Z^+$;} \\
	0, &\text{otherwise.} 
\end{cases}
\]
Then $f(x)$ is Riemann integrable and $\int_0^1 f(x)\ dx = 0$ but it is
not regulated.
\end{ex}

\begin{proof}
We define a step function $u_m(x)$ by
\[
u_m(x) = 
\begin{cases}
	1, &\text{ if } 0 \le x \le \frac{1}{m} \text{;} \\
	f(x), &\text{otherwise.} 
\end{cases}
\]
A partition for this step function is given by 
\[
x_0 = 0 < x_1 = \frac{1}{m} < x_2 = \frac{1}{m-1} < \dots <
x_{m-1} = \frac{1}{2}< x_m = 1.
\]
Note that $u_m(x) \ge f(x).$ Also $\int_0^1 u_m(x) \ dx = 1/m.$
This is because it is constant and equal to $1$ on the interval
$[0,1/m]$ and except for a finite number of points it is constant
and equal to $0$ on the interval $[1/m, 1].$
Hence
\[
\inf_{u \in \U(f)} \Big \{\int_0^1 u(x)\ dx \Big\}
\le \inf_{m \in \Z^+} \Big \{\int_0^1 u_m(x)\ dx \Big\}
= \inf_{m \in \Z^+} \{\frac{1}{m}\} = 0.
\]
Also the constant function $0$ is $\le f(x)$ and its integral is $0$,
so
\[
0 \le \sup_{v \in \yL(f)} \Big\{\int_0^1 v(x)\ dx \Big\}.
\]
Putting together the last two inequalities with
Proposition (\ref{prop: ineq}) we obtain
\[
0 \le \sup_{v \in \yL(f)} \Big \{\int_0^1 v(x)\ dx \Big\}
 \le \inf_{u \in \U(f)} \Big \{\int_0^1 u(x)\ dx \Big\}\le 0.
\]
So all of these inequalities are equalities and
by definition, $f$ is Riemann integrable with integral $0.$

To see that $f$ is not regulated suppose that $g$ is an 
approximating step function with partition
$x_0 =0 < x_1 < \dots < x_m = 1$  and satisfying
$|f(x) - g(x)| \le \eps$ for some $\eps > 0.$
Then $g$ is constant,
say with value $c_1$ on the open interval $(0,x_1).$ 

Now there are points $a_1, a_2 \in (0,x_1)$ with $f(a_1) =0$
and $f(a_2) = 1.$  Then
$|c_1| = |c_1 -0| = |g(a_1) - f(a_1)| \le \eps$ and
$|1 - c_1| = |f(a_2) - g(a_2)| \le \eps.$  But 
$|c_1| + |1 - c_1| \ge |c_1 + 1 -c_1| =1$ so at
least one of $|c_1|$ and $|1 - c_1|$ must be $\ge 1/2.$
This implies that $\eps \ge 1/2.$ That is, $f$ cannot
be uniformly approximated by any step function to within $\eps$
if $\eps < 1/2.$  So $f$ is not regulated.
\end{proof}

\begin{thm}[Regulated functions are Riemann integrable]
Every regulated function $f$ is Riemann integrable and the 
regulated integral of $f$ is equal to its Riemann integral.
\end{thm}

\begin{proof}
If $f$ is a regulated function on the interval $\I = [a,b]$, then,
for any $\eps > 0$, it can be uniformly approximated within
$\eps$ by  a step function.  In particular, if 
$\eps_n = 1/2^n$ there is a step function $g_n(x)$ such that
$|f(x) - g_n(x)| < \eps_n$ for all $x \in \I.$
The regulated integral $\int_a^b f(x) \ dx$ was defined to be
$\lim \int_a^b g_n(x) \ dx.$

We define two other approximating sequences of step functions for $f.$
Let $u_n(x) = g_n(x) + 1/2^n$ and $v_n(x) = g_n(x) - 1/2^n.$ Then
$u_n(x) \ge f(x)$ for all $x \in \I$ because $u_n(x) - f(x) = 1/2^n +
g_n(x) - f(x) \ge 0$ since $|g_n(x) - f(x)| < 1/2^n.$ Similarly
$v_n(x) \le f(x)$ for all $x \in \I$ because $f(x) - v_n(x) = 1/2^n +
f(x) - g_n(x) \ge 0$ since $|f(x) - g_n(x)| < 1/2^n.$

Since $u_n(x) - v_n(x) = g_n(x) + 1/2^n  - (g_n(x) - 1/2^n) = 1/2^{n-1},$
\[
\int_a^b u_n(x)\ dx - \int_a^b v_n(x)\ dx 
=  \int_a^b u_n(x) -  v_n(x)\ dx =  \int_a^b \frac{1}{2^{n-1}}\ dx =  
\frac{b-a}{2^{n-1}}.
\]
Hence we may apply Theorem(\ref{thm: riem criterion}) to conclude that
$f$ is Riemann integrable.

%Clearly, since $\lim 1/2^n = 0,$
Also
\begin{align*}
\lim_{n \to \infty} \int_a^b g_n(x) \ dx &=
\lim_{n \to \infty} \int_a^b v_n(x) + \frac{1}{2^n} \ dx
= \lim_{n \to \infty} \int_a^b v_n(x) \ dx, \text{ and }\\
\lim_{n \to \infty} \int_a^b g_n(x) \ dx &=
\lim_{n \to \infty} \int_a^b u_n(x) - \frac{1}{2^n} \ dx
= \lim_{n \to \infty} \int_a^b u_n(x) \ dx. 
\end{align*}
Since for all $n$
\[
\int_a^b v_n(x) \ dx \le \int_a^b f(x) \ dx \le \int_a^b u_n(x) \ dx
\]
we conclude that
\[
\lim_{n \to \infty} \int_a^b g_n(x) \ dx = \int_a^b f(x) \ dx.
\]
That is, the regulated integral equals the Riemann integral.
\end{proof}

\begin{thm}
The set $\RR$ of Riemann integrable functions on an interval $\I=[a,b]$
is a vector space containing the vector space of regulated functions.
\end{thm}
\begin{proof}
We have already shown that every regulated function is Riemann 
integrable.  Hence we need only show
that whenever $f,g \in \RR$ and $r \in \R$ we also have $(f + g) \in \RR$ 
and $rf \in \RR.$
We will do only the sum and leave the product as an exercise.

Suppose $\eps > 0$ is given.
Since $f$ is Riemann integrable there are step functions $u_f$ and
$v_f$ such that $v_f(x) \le f(x) \le u_f(x)$ for $x \in \I$ 
(i.e. $u_f \in \U(f)$ and $v_f \in \yL(f)$) and with the property
that
\[
\int_a^b u_f(x)\ dx - \int_a^b v_f(x)\ dx < \eps.
\]
Similarly there are
$u_g \in \U(g)$ and $v_g \in \yL(g)$ with the property
that
\[
\int_a^b u_g(x)\ dx - \int_a^b v_g(x)\ dx < \eps.
\]
This implies that 
\[
\int_a^b (u_f + u_g)(x)\ dx - \int_a^b (v_f + v_g)(x)\ dx < 2\eps.
\]

Since $(u_f + u_g) \in \U(f + g)$ and $(v_f + v_g) \in \yL(f + g)$ 
we may conclude that
\[
\inf_{u \in \U(f + g)} \Big \{\int_a^b u(x)\ dx \Big\}- 
\sup_{v \in \yL(f + g)} \Big \{\int_a^b v(x)\ dx \Big\}< 2\eps.
\]
As $\eps > 0$ is arbitrary we conclude that
\[
\inf_{u \in \U(f + g)} \Big \{\int_a^b u(x)\ dx \Big\}=
\sup_{v \in \yL(f + g)} \Big \{\int_a^b v(x)\ dx\Big\}
\]
and hence $(f+g) \in \RR.$
\end{proof}

\begin{exer}\ 
\begin{enumerate}
\item Prove that if $f$ and $g$ are Riemann integrable functions on
an interval $[a,b]$, then
so is $fg.$  In particular if $r \in \R$, then $rf$ is a Riemann integrable
function on $[a,b]$.
\end{enumerate}
\end{exer}

\vfill\eject

\chapter{Lebesgue Measure}\label{chap: leb meas}
\section{Introduction}
In the previous section we studied two definitions of integral that
were based on two important facts: (1) There is only the one obvious
way to define the integral of step assuming we want it to satisfy
certain basic properties, and (2) these properties force the
definition for the integral for more general functions which are
uniformly approximated by step functions (regulated integral) or
squeezed between step functions whose integrals are arbitrarily close
(Riemann integral).

To move to a more general class of functions we first find a 
more general notion to replace step functions.  For a step 
function $f$ there is a partition of $I = [0,1]$ into intervals on each of
which $f$ is constant. We now would like to allow functions
for which there is a finite partition of $I$ into sets on each which
$f$ is constant, but with the sets {\em not necessarily intervals.} 
For example we will consider functions like
\begin{equation}\label{eqn: rational-irrational}
f(x) = 
\begin{cases}
	3, &\text{if $x$ is rational;} \\
	2, &\text{otherwise.} 
\end{cases}
\end{equation}

The interval $I$ is partitioned into two sets $A = I \cap \Q$ and
$B = I \cap \Q^c$, i.e. the rational points of $I$ and the
irrational points.  Clearly the integral of this function
should be $3 \len(A) + 2  \len(B)$, but only if we can make sense of
$\len(A)$ and $\len(B)$.  That is the problem to which this
chapter is devoted.  We want to generalize the concept of
length to include as many subsets of $\R$ as we can.
We proceed in much the same way we
did in previous chapters.  We first decide what are the
``obvious'' properties this generalized length must
satisfy to be of any use, and, then try to define it by
approximating with simpler sets where the definition is
clear, namely sets of intervals.

The generalization of length we want is called {\em
Lebesgue measure}.  Ideally we would like it to work
for {\it any} subset of the
interval $I = [0,1]$, but it turns out that it is not
possible to achieve that.  

There are several properties which we want any notion
of ``generalized length'' to satisfy.
For each bounded subset $A$ of $\R$ we would like to be able
to assign a non-negative
real number $\mu(A)$ that satisfies the following:

\begin{description}
\item[I. Length.] If $A = (a,b)$ or $[a,b]$ 
then $\mu(A) = \len (A) = b-a,$ i.e. the measure of an open or closed interval
is its length
\item [II. Translation Invariance.]
If $A \subset \R$ is a bounded subset of $\R$ and $c\in \R.$, then 
$\mu( A + c ) = \mu(A),$ where $A + c$ denotes the set
$\{ x +c\ |\ x \in A\}.$

\item [III. Countable Additivity.]
If $\{A_n\}_{n=1}^\infty$ is a countable collection of
bounded subsets of $\R$, then
\[
\mu( \bigcup_{n=1}^\infty A_n) \le \sum_{n=1}^\infty \mu(A_n)
\]
and if the sets are  {\em pairwise disjoint}, then
\[
\mu( \bigcup_{n=1}^\infty A_n) = \sum_{n=1}^\infty \mu(A_n)
\]
Note the same conclusion applies to finite collections
$\{A_n\}_{n=1}^m$ of bounded sets (just let $A_i = \emptyset$  for $i > m$).

\item [IV. Monotonicity] If $A \subset B$, then $\mu(A) \le \mu(B).$
Actually, this property is a consequence of additivity since $A$ and
$B \setminus A$ are disjoint and their union is $B$.
\end{description}
\bigskip

It should be fairly clear why we most of these properties are
absolutely necessary for any sensible notion of length.  The only
exception is property III, which deserves some comment.  We might ask
that additivity only hold for finite collections of sets, but that is
too weak.  For example, if we had a collection of pairwise disjoint
intervals of length $1/2, 1/4, 1/8,\dots 1/2^n,\dots,$ etc., then we
would certainly like to be able say that the measure of their union is
the sum $\sum 1/2^{n} = 1$ which would not follow from finite
additivity.  Alternatively, one might wonder why
additivity is only for {\em countable} collections of pairwise
disjoint sets. But it is easy to see why it would lead to problems if
we allowed uncountable collections.  If $\{x\}$ is the set consisting
of a single point $x\in [0,1]$, then $\mu(\{x\}) = 0$ by property I and
$[a,b]$ is an uncountable union of pairwise disjoint sets, namely each
of the sets consisting of one point of $[a,b]$.  Hence we would have
$\mu([a,b]) = b-a$ is an {\em uncountable} sum of zeroes.  This is one
reason the concept of uncountable sums isn't very useful.  Indeed, we
will see that the concept of countability is intimately related to the
concept of measure.

Unfortunately, as mentioned above, it turns out that it is impossible
to find a $\mu$ which satisfies I--IV and which is defined for {\em
all} bounded subsets of the reals.  But we can do it for a very large
collection which includes all the open sets and all the closed sets.
The measure we are interested in using is called {\em Lebesgue
measure} Its actual construction is slightly technical and we have
relegated that to an appendix.  Instead we will focus the properties
of Lebesgue measure and how to use it.

\section{Null Sets}

One of our axioms for the regulated integral was, ``Finite
sets don't matter.'' Now we want to generalize that to 
say that sets whose ``generalized length,'' or measure,
is zero don't matter.  It is a somewhat surprising fact
that even without defining Lebesgue measure in general
we can easily define those sets whose 
measure must be $0$ and investigate the properties of these
sets.

\begin{defn}[Null Set]\index{null set}\label{def: null}
A set $X \subset \R$ is called a {\em null set} provided
for every $\eps > 0$ there is a collection of
open intervals $\{U_n\}_{n=1}^\infty$ such that
\[
\sum_{n=1}^\infty \len(U_n) < \eps
\text{ and }
X \subset \bigcup_{n=1}^\infty U_n.
\]
\end{defn}

Perhaps surprisingly this definition makes no use of the measure 
$\mu$.  Indeed, we have not yet defined the measure $\mu$
of a set $X$ for any choice of the set $X$!
However, it is clear that if we can do so in a
way that satisfies properties I-IV above and the
Hence to denote $X$ has measure zero, we will write $\mu(X) = 0$ even
though we have not yet defined $\mu.$ 

If $X$ is a null set in $I = [0,1]$ we will say that its complement
$X^c$ has \index{full measure} {\em full measure} in $I$.

\begin{exer}\
\begin{enumerate}
\item Prove that a finite set is a null set.
\item Prove that a countable union of null sets
is a null set (and hence, in particular, countable
sets are null sets).
\item Assuming that a measure $\mu$ has been defined
and satisfies properties I-IV above, find the numerical
value of the integral of the function $f(x)$ defined
in Equation (\ref{eqn: rational-irrational}).  Prove
that the Riemann integral of this function does not
exist.
\end{enumerate}
\end{exer}

It is not true that countable sets are the only sets which are null
sets.  We give an example in Exercise~(\ref{exer: measure}) below,
namely, the Cantor middle third set, which is an uncountable null set.

\section{Sigma algebras}

As mentioned before there does not exist function $\mu$ satisfying
properties I-IV and defined for every subset of $I = [0,1]$.
In this section we want to consider what is the
best collection of
subsets of $I$ for which we can define a ``generalized length''
or {\em measure} $\mu$.
Suppose we have somehow defined 
$\mu$ for all the sets in some collection
$\A$ of subsets of $I$ and it satisfies properties I--IV.
Property I only makes sense if $\mu$ is defined
for open and closed intervals, i.e. we need open and closed intervals
to be in $\A$.  For property III to make sense we will need
that any countable union of sets in $\A$ is also in $\A.$  Finally it seems
reasonable that if $A$ is a set in the collection $\A$, then
the set $Z^c$, its complement in $I,$ should also be in $\A$.

All this motivates the following definition.

\begin{defn}[Sigma algebra]\index{Sigma algebra}
Suppose $X$ is a set and $\A$ is a collection of subsets of $X$.
$\A$ is called a
\myindex{$\sigma$-algebra} of subsets of $X$ 
provided it contains the set $X$ and is
closed under taking complements (with respect to $X$), 
countable unions, and countable intersections.
\end{defn}

In other words if $\A$ is a $\sigma$-algebra of subsets of 
$X$, then any complement (with respect to $X$) of a set in $\A$ is also in $\A$,
any countable union of sets in $\A$ is in $\A$, and
any countable intersection of sets in $\A$ is in $\A$.
In fact the property about countable intersections follows
from the other two and Proposition (\ref{prop: comp intersection})
which says the intersection of a family of sets is the complement
of the union of the complements of the sets. Also notice that if 
$A,B \in \A$, then their set difference
$A \setminus B = \{x \in A\ |\ x \notin B\}$ is in $\A$
because $A \setminus B = A \cap B^c.$

Since $X$ is in any $\sigma$-algebra of subsets of $X$ (by
definition), so is its complement, the empty set.  A trivial example
of a $\sigma$-algebra of subsets of $X$
is $\A = \{ X, \emptyset\}$, i.e. it consists of
only the whole set $X$ and the empty set.  Another example is
$\A = \PP(X)$, the power set of $X$, i.e. the
collection of all subsets of $X$.  Several more interesting examples
are given in the exercises below.  Also in these exercises
we ask you to show that any intersection of $\sigma$-algebras
is a $\sigma$-algebra.  Thus for any collection $\C$ of subsets
of $\R$ there is a smallest  $\sigma$-algebra
of subsets of $\R$ which contains all sets in $\C$, 
namely the intersection of all $\sigma$-algebras containing
$\C$.

\begin{defn}[Borel Sets]
If $\C$ is a collection of subsets of $\R$ and $\A$ is
the the smallest $\sigma$-algebra
of subsets of $\R$ which contains all the sets of $\C$
then $\A$ is called the 
$\sigma$-algebra generated by $\C$.
\index{$\sigma$-algebra generated by a family
of sets}
Let $\B$ be the $\sigma$-algebra
of subsets of $\R$ generated by the collection of all open intervals.
$\B$ is called the \myindex{Borel $\sigma$-algebra}
and elements of $\B$ are called \myindex{Borel sets}.
\end{defn}

In other words $\B$ is the collection of subsets of $\R$ which can be formed
from open intervals by any finite sequence of countable unions,
countable intersections or complements.

\begin{exer}\label{exer: Borel}\
\begin{enumerate}
\item Let $\A = \{ X \subset I\ |\ X \text{ is countable, or }
X^c\text{ is countable}\}.$ Prove that $\A$ is a $\sigma$-algebra.

\item Let $\A = \{ X \subset I\ |\ X \text{ is a null set, or }
X^c\text{ is a null set}\}.$ Prove that $\A$ is a $\sigma$-algebra.

\item Suppose $\A_\lambda$ is a $\sigma$-algebra of subsets of $X$ for
each $\lambda$ in some indexing set $\Lambda.$ Prove that
\[
\A = \bigcap_{\lambda \in \Lambda}\A_\lambda
\]
is a $\sigma$-algebra of subsets of $X$.

\item Let $\A$ be a $\sigma$-algebra of subsets of $\R$ and
suppose $I$ is a closed interval which is in $\A$.
Let $\A(I)$ denote the collection of all subsets of $I$ which
are in $\A$.  Prove that $\A(I)$ is 
a $\sigma$-algebra of subsets of $I$.

\
\item Suppose $\C_1$ is the collection of closed intervals in $\R$,

$\C_2$ is the collection of all open subsets of $\R$, and

$\C_3$ is the collection of all closed subsets of $\R$.  

Let 
$\A_i$ be the $\sigma$-algebra generated by
$\C_i$.  Prove that
$\B_1, \B_2,$ and $\B_3$ are all equal to the
Borel $\sigma$-algebra $\B.$
\end{enumerate}
\end{exer}

\section{Lebesgue Measure}
The $\sigma$-algebra of primary interest to us is the one
generated by Borel sets and null sets.  Alternatively,
as a consequence of part 5. of Exercise~(\ref{exer: Borel}),
it is the $\sigma$-algebra of subsets of $\R$ generated
by open intervals, and null sets, or the one
generated by closed intervals and null sets.

\begin{defn}\label{def: leb sigma}
\index{Lebesgue measurable} \index{measurable}
The $\sigma$-algebra of subsets of $\R$
generated by open intervals and null sets will
be denoted $\M$.  Sets in $\M$ will be called Lebesgue measurable,
or measurable for short.  If $I$ is a closed interval, then
$\M(I)$ will denote the Lebesgue measurable subsets of $I$.
\end{defn}

For simplicity we will focus on subsets of $I = [0,1]$ though
we could use any other interval just as well.
Notice that it is a consequence of part 4. of 
Exercise~(\ref{exer: Borel}) that
$\M(I)$ is a $\sigma$-algebra of
subsets of $I$.
It is by no means obvious that $\M$ is not the
$\sigma$-algebra of all subsets of $\R$.  However, in section
(\ref{sec: non-measurable}) of the appendix we will construct a subset
of $I$ which is not in $\M$.

We are now ready to state the main theorem of this Chapter.

\begin{thm}[Existence of Lebesgue Measure]\label{thm: lebesgue measure}
\index{Lebesgue measure}
There exists a unique  function $\mu$, called {\em Lebesgue measure},
from $\M(I)$ to the non-negative real numbers satisfying:
\begin{description}
\item[I. Length.] If $A = (a,b)$ 
then $\mu(A) = \len (A) = b-a,$ i.e. the measure of an open interval
is its length
\item [II. Translation Invariance.]
Suppose $A \subset I,\ c\in \R$ and $A +c \subset I$
where $A + c$ denotes the set $\{ x +c\ |\ x \in A\}.$
Then  $\mu( A + c ) = \mu(A)$

\item [III. Countable Additivity.]
If $\{A_n\}_{n=1}^\infty$ is a countable collection of
subsets of $I$, then
\[
\mu( \bigcup_{n=1}^\infty A_n) \le \sum_{n=1}^\infty \mu(A_n)
\]
and if the sets are  {\em pairwise disjoint}, then
\[
\mu( \bigcup_{n=1}^\infty A_n) = \sum_{n=1}^\infty \mu(A_n)
\]
\item [IV. Monotonicity] If $A,B \in \M(I)$ and $A \subset B$
then $\mu(A) \le \mu(B)$
\item [V. Null Sets] A subset $A \subset I$ is a null set
set if and only if $A \in \M(I)$ and $\mu(A) = 0.$
\end{description}
\end{thm}

Note that the countable additivity of property III implies
the analogous statements about finite additivity.
Given a finite collection $\{A_n\}_{n=1}^m$ of sets 
just let $A_i = \emptyset$  for $i > m$ and the
analogous conclusions follow.

We have relegated the proof of most of this theorem to Appendix~A,
because it is somewhat technical and is a diversion from our main task
of developing a theory of integration. However there are some
properties of Lebesgue measure we can easily derive; so we do so now.
For example, we will use properties I-III of Theorem~(\ref{thm:
lebesgue measure}) to prove property IV.

\begin{prop}[Monotonicity] 
If $A,B \in \M(I)$ and
$A \subset B$, then $\mu(A) \le \mu(B).$
\end{prop}
\begin{proof}
Since $A \subset B$ we have $B = A \cup (B \setminus A)$.
Also $A$ and $B \setminus A$ are disjoint so
by property III we know $\mu(A) + \mu(B\setminus A) = \mu(B).$
But $\mu(B \setminus A) \ge 0$ so 
$\mu(A) \le \mu(A) + \mu(B\setminus A) = \mu(B).$
\end{proof}

\begin{prop}
If $X \subset I$ is a null set, then $X \in\M(I)$ and $\mu(X) = 0.$
\end{prop}

\begin{proof}
If $X \subset I$ is a null set, then by the definition of
$\M(I)$ we know $X \in \M(I)$.  If $\eps >0$
there is a collection of
open intervals $\{U_n\}_{n=1}^\infty$ such that
\[
\sum_{n=1}^\infty \len(U_n) < \eps
\text{ and }
X \subset \bigcup_{n=1}^\infty U_n.
\]

Property I says $\len(U_n) = \mu(U_n)$, so
\[
\sum_{n=1}^\infty \mu(U_n) = \sum_{n=1}^\infty \len(U_n) < \eps.
\]
Since $X \subset \cup U_n$ properties II and III imply
\[
\mu(X) \le \mu\big ( \bigcup_{n=1}^\infty U_n\big ) \le
\sum_{n=1}^\infty \mu(U_n) < \eps,
\]
This is true for any $\eps >$ so the only possible
value for $\mu(X)$ is zero.
\end{proof}

Recall that \myindex{set difference}
$A \setminus B = \{x \in A\ |\ x \notin B\}$.  Since
we are focusing on subsets
of $I$ complements are with respect to $I$ so 
$A^c = I \setminus A.$ 

\begin{prop}\ 
\begin{enumerate}
\item[(1)] The Lebesgue measure of $I$,\ $\mu(I)$, is $1$ and
hence $\mu(A^c) = 1 - \mu(A).$
\item[(2)] If $A$ and $ B$ are in $\M(I)$, then $A \setminus B$ is in $\M$
and $\mu( A \cup B) = \mu(A \setminus B) + \mu(B).$
\end{enumerate}
\end{prop}
\begin{proof}
To see (1) observe that $A$ and $A^c$ are disjoint
and $A \cup A^c = I$, so additivity
implies $\mu(A) + \mu(A^c) = \mu(A \cup A^c) = \mu(I) = 1.$

For (2) note that $A \setminus B = A \cap B^c$ which is in $\M$.
Also $A\setminus B$ and $B$ are disjoint and their union is $A \cup B$.
So once again additivity implies 
implies $\mu(A\setminus B) + \mu(B) = \mu(A \cup B).$
\end{proof}

If we have a countable increasing family of measurable sets
then the measure of the union can be expressed as a limit.

\begin{prop}\label{prop: increasing union}
If $A_1 \subset A_2 \subset \dots \subset A_n \dots$ is an increasing 
sequence of measurable subsets of $I$, then
\[
\mu(\bigcup_{n=1}^\infty A_n) = \lim_{n \to \infty} \mu(A_n).
\]
If $B_1 \supset B_2 \supset \dots \supset B_n \dots$ is a decreasing
sequence of measurable subsets of $I$, then
\[
\mu(\bigcap_{n=1}^\infty B_n) = \lim_{n \to \infty} \mu(B_n).
\]
\end{prop}
\begin{proof}
Let $F_1 = A_1$ and $F_n = A_n \setminus A_{n-1}$ for $n > 1.$
Then $\{F_n\}_{n=1}^\infty$ are pairwise disjoint measurable
sets, $A_n = \cup_{i=1}^n F_i$ and 
$\cup_{i=1}^\infty A_i = \cup_{i=1}^\infty F_i.$
Hence by countable additivity we have
\begin{align*}
\mu \big( \bigcup_{i=1}^\infty A_i \big)  
&= \mu\big( \bigcup_{i=1}^\infty F_i\big)
= \sum_{i=1}^\infty \mu(F_i)\\
&= \lim_{n \to \infty} \sum_{i=1}^n \mu(F_i)
= \lim_{n \to \infty} \mu\big (\bigcup_{i=1}^n F_i \big )\\
&= \lim_{n \to \infty} \mu(A_n).
\end{align*}

For the decreasing sequence we define
$E_n = B_n^c.$  Then $\{E_n\}_{n=1}^\infty$ is an increasing
sequence of measurable functions and 
\[
\big ( \bigcap_{n=1}^\infty B_n \big)^c =  \bigcup_{n=1}^\infty E_n.
\]
Hence 
\[
\mu\big( \bigcap_{n=1}^\infty B_n \big ) 
= 1 - \mu\big( \bigcup_{i=1}^\infty E_i \big )
= 1 - \lim_{n \to \infty}\mu( E_n)
= \lim_{n \to \infty} (1 - \mu( E_n))
= \lim_{n \to \infty} \mu( B_n).
\]
\end{proof}

\section{The Lebesgue Density Theorem}

The following theorem asserts that if a subset of an interval $I$
is ``equally distributed'' throughout the interval
then it must be a null set or the complement of a null set.
For example it is not possible to have a set
$A \subset [0,1]$ which contains half of each subinterval, i.e.
it is impossible to have $\mu(A \cap [a,b]) = \mu( [a,b])/2$ for all
$0 < a < b < 1.$  There will always be small intervals with 
a ``high concentration'' of points of $A$ and other subintervals
with a low concentration.  Put another way, it
asserts that given any $p < 1$ there is an interval $U$ such
that a point in $U$ has probability $> p$ of being in $A$.

\begin{thm}\label{thm: easy lebesgue density}
If $A$ is a Lebesgue measurable set and $\mu(A) > 0$ and if
$0 < p < 1$, then there is an open interval $U = (a,b)$ such that
$\mu(A \cap U) \ge p \mu(U) = p(b-a).$
\end{thm}
\begin{proof}
Let $p \in (0,1)$ be given.  We know from the definition of outer
measure and the fact that $\ms(A) = \mu(A),$ that for any $\eps
>0$ there is a countable open cover $\{U_n\}_{n=1}^\infty$ of $A$ such
that
\[
\mu(A) \le \sum_{n=1}^\infty \len(U_n) \le \mu(A) + \eps.
\]
Choosing $\eps = (1 - p) \mu(A)$ we get 
\begin{align*}
\sum_{n=1}^\infty \len(U_n) &\le \mu(A) + (1 -p)\mu(A)\\
&\le \mu(A) + (1 - p)\sum_{n=1}^\infty \len(U_n)
\end{align*}
so 
\begin{equation}\label{eqn:density}
p\sum_{n=1}^\infty \len(U_n) \le \mu(A) \le \sum_{n=1}^\infty \mu(A \cap U_n).
\end{equation}
where the last inequality follows from subadditivity.  
Since these infinite series have finite sums, 
there is at least one $n_0$ such that
$p \mu( U_{n_0}) \le \mu( A \cap U_{n_0})$.  This is because if it were the
case that $p \mu( U_{n}) > \mu( A \cap U_{n})$ for all $n$, then
it would follow that
$p\sum_{n=1}^\infty \len(U_n) > \sum_{n=1}^\infty \mu(A \cap U_n)$
contradicting equation (\ref{eqn:density}).  The interval $U_{n_0}$ is
the $U$ we want.
\end{proof}

There is a much stronger result than the theorem above which we now state,
but do not prove. A proof can be found in Section 9.2 of \cite{T}.

\begin{defn}
If $A$ is a Lebesgue measurable set and $x \in A$, then $x$
is called a \myindex{Lebesgue density point} if
\[
\lim_{\eps \to 0} \frac{\mu( A \cap [x-\eps, x+\eps])}{\mu([x-\eps, x+\eps])} = 1.
\]
\end{defn}

\begin{thm}[Lebesgue Density Theorem]
\index{Lebesgue Density Theorem}
If $A$ is a Lebesgue measurable set, then there is a subset $E \subset A$
with $\mu(E) = 0$ such that every point of $A \setminus E$ is a Lebesgue
density point.
\end{thm}

\section{Lebesgue Measurable Sets -- Summary}

In this section we provide a summary outline of the key 
properties of collection $\M$ of Lebesgue measurable sets
which have been developed in this chapter.
Recall $I$ is a closed interval and $\M(I)$ denotes the
subsets of $I$ which are in $I$.

\begin{enumerate}
\item The collection of Lebesgue measurable sets $\M$ is a 
\myindex{$\sigma$-algebra}, which  means
\subitem $\bullet$ If $A \in \M$, then $A^c \in M.$
\subitem $\bullet$ If $A_n \in \M$ for $n \in \N$, then $\bigcup_{n=1}^\infty A_n \in M.$
\subitem $\bullet$ If $A_n \in \M$ for $n \in \N$, then $\bigcap_{n=1}^\infty A_n \in M.$
\item All open sets  and  all closed sets are in $\M$.
Any null set is in $\M$.
\item If $A \in \M(I)$, then there is a real number $\mu(A)$ called
its Lebesgue measure which satisfies:
\subitem $\bullet$ The Lebesgue measure of an interval is its length.
\subitem $\bullet$ Lebesgue measure is translation invariant.
\subitem $\bullet$ If $A \in \M$, then $\mu(A^c) = 1 - \mu(A).$
\subitem $\bullet$ If $A \in \M$ is a null set if and only if $\mu(A) = 0.$
\subitem $\bullet$ \myindex{Countable Subadditivity}:
If $A_n \in \M$ for $n \in \N$, then 
\[
\mu\big( \bigcup_{n=1}^\infty A_n\big) \le \sum_{n=1}^\infty \mu(A_n).
\]
\subitem $\bullet$ \myindex{Countable Additivity}:
If $A_n \in \M$ for $n \in \N$ are pairwise disjoint, then 
\[
\mu\big( \bigcup_{n=1}^\infty A_n\big) = \sum_{n=1}^\infty \mu(A_n).
\]
\subitem $\bullet$ {\em Increasing sequences}: If $A_n \in \M$ for $n
\in \N$ satisfy $A_n \subset A_{n+1}$, then
\[
\mu\big( \bigcup_{n=1}^\infty A_n\big) = \lim_{n \to \infty} \mu(A_n).
\]
\subitem $\bullet$ {\em Decreasing sequences}: If $A_n \in \M$ for $n
\in \N$ satisfy $A_n \supset A_{n+1}$, then
\[
\mu\big( \bigcap_{n=1}^\infty A_n\big) = \lim_{n \to \infty} \mu(A_n).
\]
\end{enumerate}

\vfill\eject

\begin{exer}\label{exer: measure}\ 
\begin{enumerate}
\item Prove for $a,b \in I$ that $\mu([a,b]) = \mu([a,b[) = b-a.$
\item Let $X$ be the subset of irrational numbers in $I$.  Prove
$\mu(X) = 1.$  Prove that if $Y \subset I$ is a closed set and
$\mu(Y) = 1$, then $Y = I.$

\item ({\em The Cantor middle third set})\index{Cantor middle third set} 
We first recursively define a nested sequence $\{J_n\}_{n=0}^\infty$
of closed subsets of $I$.  Each $J_n$ consists of a finite union of
closed intervals.  We define $J_0$ to be $I = [0,1]$ and let $J_n$ be
the union of the closed intervals obtained by deleting the open middle
third interval from each of the intervals in $J_{n-1}.$ Thus $J_0 =
[0,1],\ J_1 = [0,1/3] \cup [2/3,1]$ and $J_2 = [0,1/9] \cup [2/9,1/3]
\cup [2/3,7/9] \cup [8/9,1]$ etc.

Let $C = \cap_{n=0}^\infty J_n.$ It is called the {\em Cantor Middle Third
set.}

\subitem {(a)} When the open middle thirds of the intervals in
$J_{n-1}$ are removed we are left with two sets of closed intervals:
the left thirds of the intervals in $J_{n-1}$ and the right thirds of
these intervals.  We denote the union of the left thirds by $L_n$ and
the right thirds by $R_n$, and we note note $J_n = L_n \cup R_n.$
Prove that $L_n$ and $R_n$ each consist of $2^{n-1}$ intervals of
length $1/3^n$ and hence $J_n$ contains $2^n$ intervals of length
$1/3^n$.  

\subitem{(b)} Let $\D$ be the uncountable set
set of all infinite sequences $d_1 d_2 d_3 \dots d_n
\dots$ where each $d_n$ is either $0$ or $1$ (see part 4. of Exercise 
(\ref{exer: uncountable}))
and define a function $\psi: C \to \D$ by $\psi(x) =
d_1 d_2 d_3 \dots d_n \dots$ where each $d_n = 0$ if $x \in L_n$ and
$d_n = 1$ if $x \in R_n$.  Prove that $\psi$ is surjective and hence
by Corollary~(\ref{cor: uncountable}) the set
$C$ is uncountable. {\em Hint: } You will need to use 
Theorem~(\ref{thm: nested}).

\subitem{(c)}  Prove that $C$ is Lebesgue measurable and that $\mu(C) = 0.$ 
{\em Hint: } Consider $C^c$, the complement of $C$ in $I$. Show it is 
measurable and calculate $\mu(C^c).$  {\em Alternative hint: } Show
directly that $C$ is a null set by finding for each $\eps >0$
a collection of open intervals $\{U_n\}_{n=1}^\infty$ such that
\[
\sum_{n=1}^\infty \len(U_n) < \eps
\text{ and }
C \subset \bigcup_{n=1}^\infty U_n.
\]
\end{enumerate}
\end{exer}
\vfill\eject

\chapter{The Lebesgue Integral}
\section{Measurable Functions}

In this chapter we want to define the Lebesgue integral in
a fashion which is analogous to our definitions of regulated
integral and Riemann integral from Chapter 1. The difference
is that we will no longer use step functions to approximate
a function we want to integrate, but instead will use a
much more general class called simple functions. 

\begin{defn}[Characteristic Function]
If $A \subset [0,1],$ its \myindex{characteristic function} $\X_A(x)$ is 
defined by
\[
\X_A(x) = 
\begin{cases}
	1, &\text{if $x \in A$;} \\
	0, &\text{otherwise.} 
\end{cases}
\]
\end{defn}

\begin{defn}[Measurable partition]
A finite \myindex{measurable partition} of $[0,1]$ is a collection
$\{A_i\}_{i=1}^n$ of measurable subsets which are pairwise
disjoint and whose union is $[0,1]$. 
\end{defn}

We can now define \myindex{simple functions}.  Like step
functions these functions have only finitely many values,
but unlike step functions the set on which a simple function
assumes a given value is no longer an interval. Instead
a simple function is constant on each subset of a finite
measurable partition of $[0,1].$

\begin{defn}[Simple Function]\label{def: simple}
A function $f: [0,1] \to \R$ is called 
\myindex{Lebesgue simple} or \myindex{simple}, for short, provided
there exist a finite measurable partition
$\{A_i\}_{i=1}^n$ and real numbers $r_i$ such that
$f(x) = \sum_{i=1}^n r_i \X_{A_i}.$
The \myindex{Lebesgue integral} of a simple function is
defined by $\int f\ d\mu = \sum_{i=1}^n r_i \mu(A_i).$
\end{defn}

The definition of integral of a simple function should come
as no surprise. The fact that $\int \X_A(x)\ d\mu$ is defined
to be $\mu(A)$ is the generalization of the fact that
the Riemann integral $\int_a^b 1\ dx = (b-a)$.  The value
of $\int f\ d\mu$ for a step function $f$ is, then forced if
we want our integral to have the linearity property.

\begin{lemma}[Properties of simple functions]\label{lem: simple functions}
The set of simple functions is a vector space and the Lebesgue
integral of simple functions satisfies the following properties:
\begin{enumerate}
\item \myindex{Linearity}: If $f$ and $g$ is simple
functions and $c_1,c_2 \in \R$, then
\[
\int c_1f + c_2g\ d\mu = c_1 \int f\ d\mu + c_2 \int g\ d\mu.
\]
\item \myindex{Monotonicity}: If $f$ and $g$ are simple and $f(x) \le g(x)$
for all $x$, then $\int f\ d\mu \le \int g\ d\mu.$
\item {\em Absolute value}: If $f$ is simple 
then $|f|$ is simple and $|\int f\ d\mu| \le \int |f|\ d\mu.$
\end{enumerate}
\end{lemma}

\begin{proof}
If $f$ is simple, then clearly $c_1f$ is simple.  Hence to show that
simple functions form a vector space it suffices to show that the
sum of two simple functions are simple. 

Suppose
$\{A_i\}_{i=1}^n$ and $\{B_j\}_{j=1}^m$ are measurable partitions
of $[0,1]$ and that $f(x) = \sum_{i=1}^n r_i \X_{A_i}$
and $g(x) = \sum_{j=1}^m s_j \X_{B_j}$ are simple functions.
We consider the measurable partition $\{C_{i,j}\}$ with
$C_{i,j} = A_i \cap B_j.$ Then $A_i = \bigcup_{j=1}^m C_{i,j}$
and $B_j = \bigcup_{i=1}^n C_{i,j},$ so
\[
f(x) = \sum_{i=1}^n r_i \X_{A_i} =  \sum_{i=1}^n r_i \sum_{j=1}^m \X_{C_{i,j}}(x)
=  \sum_{i,j} r_i \X_{C_{i,j}}.
\]
Likewise 
\[
g(x) = \sum_{j=1}^m s_j \X_{B_j} =  \sum_{j=1}^m s_j \sum_{i=1}^n \X_{C_{i,j}}(x)
=  \sum_{i,j} s_j \X_{C_{i,j}}.
\]
Hence $f(x) + g(x) =   \sum_{i,j} (r_i + s_j) \X_{C_{i,j}}(x)$ is simple
and the set of simple functions forms a vector space.

It follows immediately from the definition that if $f$ is simple
and $a \in \R$, then $\int af \ d\mu = a \int f \ d\mu.$  So to prove linearity
we need only show that if $f$ and $g$ are simple functions as above
then $\int (f + g) \ d\mu = \int f\ d\mu  +\int g \ d\mu.$
But this follows because
\begin{align*}
\int (f + g) \ d\mu  &=   \sum_{i,j} (r_i + s_j) \mu(C_{i,j})\\
&=   \sum_{i,j} r_i \mu(C_{i,j}) + \sum_{i,j} s_j \mu(C_{i,j})\\
&= \sum_{i=1}^n r_i \sum_{j=1}^m \mu(C_{i,j}) 
+ \sum_{j=1}^m s_j \sum_{i=1}^n \mu(C_{i,j})\\
&= \sum_{i=1}^n r_i \mu(A_i)  + \sum_{j=1}^m s_j \mu(B_j)\\
&= \int f \ d\mu + \int g \ d\mu.
\end{align*}

Monotonicity follows from the fact that if $f$ and $g$ are simple functions
with $f(x) \le g(x)$, then $g(x) - f(x)$ is a non-negative simple function.
Clearly from the definition of the integral of a simple function, if
the function is non-negative, then its integral is $\ge 0.$  Thus
$\int g \ d\mu - \int f \ d\mu = \int g - f \ d\mu \ge 0.$

If $f(x) = \sum r_i \X_{A_i}$, the absolute value property follows 
from the fact that 
\[
\Big |\int f \ d\mu\Big | = \Big |\sum r_i\mu(A_i)\Big | 
\le \sum |r_i| \mu(A_i).
\]
\end{proof}

\begin{exer}\ 
\begin{enumerate}
\item Prove that if $f$ and $g$ are simple functions, then
so is $fg.$  In particular, if $E \subset [0,1]$ is measurable
then $f \X_E$ is a simple function.
\end{enumerate}
\end{exer}

A function $f: [0,1] \to \R \cup \{\infty\} \cup \{-\infty\}$ will
be called an \index{extended real valued function} extended real
valued function. For $a \in \R$
we will denote the set $(-\infty,a] \cup \{-\infty\}$ by
$[-\infty,a]$ and the set $[a,\infty) \cup \{\infty\}$ by
$[a,\infty].$

\begin{prop}\label{prop: meas func}
If $f:[0,1] \to \R$ is an extended real valued function, then the following
are equivalent:
\begin{enumerate}
\item For any $a \in [-\infty, \infty]$ the set $f^{-1}([-\infty,a])$
is Lebesgue measurable.
\item For any $a \in [-\infty, \infty]$ the set $f^{-1}([-\infty,a))$
is Lebesgue measurable.
\item For any $a \in [-\infty, \infty]$ the set $f^{-1}([a,\infty])$ is
Lebesgue measurable.
\item For any $a \in [-\infty, \infty]$ the set $f^{-1}((a,\infty])$ is
Lebesgue measurable.
\end{enumerate}
\end{prop}

\begin{proof}
We will show 1) $\Rightarrow$ 2) $\Rightarrow$ 3) 
$\Rightarrow$ 4) $\Rightarrow$ 1).

First assume 1), then $[-\infty,a[ = \bigcup_{n=1}^\infty [-\infty, a
- 2^{-n}].$ So
\[
f^{-1}([-\infty,a[) = \bigcup_{n=1}^\infty f^{-1}([-\infty, a - 2^{-n}])
\]
which is measurable by Theorem (\ref{thm: countable union}).  Hence 2)
holds.

Now assume 2), then $[a,\infty] = [-\infty,a)^c$ so  
\[
f^{-1}([a,\infty]) = f^{-1}([-\infty,a)^c) = (f^{-1}([-\infty,a)))^c.
\]
Hence 3) holds.

Assume 3), $(a,\infty] = \bigcup_{n=1}^\infty [a + 2^{-n}, \infty].$ So
\[
f^{-1}((a,\infty]) = \bigcup_{n=1}^\infty f^{-1}([a - 2^{-n},\infty])
\]
which is measurable by Theorem (\ref{thm: countable union}).  Hence 4) holds.

Finally, assume 4), then $[-\infty,a] = (a,\infty]^c$ so
\[
f^{-1}([-\infty,a]) 
= f^{-1}((a,\infty]^c) = (f^{-1}((a,\infty]))^c.
\]
  Hence 1) holds.
\end{proof}

\begin{defn}[Measurable Function]
\index{measurable}\index{Lebesgue measurable}
An extended real valued function $f$ is called Lebesgue measurable if
it satisfies one (and hence all) of the properties of Proposition
(\ref{prop: meas func}).
\end{defn}

\begin{prop}\label{prop: meas zero func}
If $f(x)$ is a function which has the value
$0$ except on a set of measure $0$, then $f(x)$ is measurable.
\end{prop}
\begin{proof}
Suppose $f(x) = 0$ for all $x \notin A$ where $A \subset [0,1]$ has
measure $0$.  That is, if $A = f^{-1}([-\infty, 0[) \cup f^{-1}((0, \infty])$
then $A$ is a null set.
For $a < 0$ the set
$U_a = f^{-1}([-\infty, a])$ is a subset of $A$ so $U_a$ is a null
set and hence measurable.
For $a \ge 0$ the set
$U_a = f^{-1}([-\infty, a])$ is the complement of the null set
$f^{-1}((a, \infty])$ and hence measurable.
In either case $U_a$ is measurable so $f$ is a measurable function.
\end{proof}

\begin{thm}\label{thm: lim}
Let $\{ f_n \}_{n=1}^\infty$ be a sequence of measurable functions.
Then the extended real valued functions 
\begin{align*}
	g_1(x) &= \sup_{n \in \N} f_n(x)\\
	g_2(x) &= \inf_{n \in \N} f_n(x)\\
	g_3(x) &= \limsup_{n \to \infty} f_n(x)\\
	g_4(x) &= \liminf_{n \to \infty} f_n(x)
\end{align*}
are all measurable.
\end{thm}

\begin{proof}
If $a \in [-\infty, \infty]$, then 
\[
\{ x\ |\ g_1(x) > a \} = \bigcup_{n=1}^\infty \{ x\ |\ f_n(x) > a \}.
\]
Each of the sets on the right is measurable so 
$\{ x\ |\ g_1(x) > a \}$ is also by Theorem (\ref{thm: countable union}).
Hence $g_1$ is measurable. 

Since $g_2(x) = \inf_{n \in \N} f_n(x) = - \sup_{n \in \N} -f_n(x)$
it follows that $g_2$ is also measurable.

Since the limit of a decreasing sequence is the $\inf$ of the terms,
\[
g_3(x) = \limsup_{n \to \infty} f_n(x) = \inf_{m \in \N}
\sup_{n \ge m} f_n(x).
\]
It follows that $g_3$ is measurable. And since 
\[
g_4(x)= \liminf_{n \to \infty} f_n(x) = - \limsup_{n \to \infty} -f_n(x)
\]
it follows that $g_4$ is measurable.
\end{proof}

For the following result we need to use honest real valued functions,
i.e., not extended.  The reason for this is that there is no way to
define the sum of two extended real valued functions if one has the
value $+\infty$ at a point and the other has the value 
$-\infty$ at the same point.

\begin{thm}\label{thm: vec space}
The set of Lebesgue measurable functions from $[0,1]$ to $\R$ is
a vector space.  The set of bounded Lebesgue measurable functions
is a vector subspace.
\end{thm}
\begin{proof}
It is immediate from the definition that for $c \in \R$ the function
$cf$ is measurable when $f$ is.  Suppose $f$ and $g$ are measurable.
We need to show that $f+g$ is also measurable, i.e., that for any $a \in \R$
the set $U_a = \{x\ |\ f(x) + g(x) > a\}$ is measurable.  

Let $\{r_n\}_{n=1}^\infty$ be an enumeration of the rationals. If
$x_0 \in U_a$, i.e., if $f(x_0) + g(x_0) > a$, then $f(x_0)  > a - g(x_0).$
Since the rationals are dense there is an $r_m$ such that 
$f(x_0)> r_m > a -  g(x_0).$  Hence if we define
\[
V_m = \{x\ |\ f(x) > r_m\} \cap \{x\ |\ g(x) > a - r_m\}
\]
then $x_0 \in V_m$.  So every point of $U_a$ is in some $V_m.$ 
Conversely if $y_0 \in V_m$ for some $m$, then
$f(y_0) > r_m > a -  g(y_0)$, so $f(y_0) + g(y_0) > a$ and $y_0 \in U_a.$
Thus $U_a = \bigcup_{m=1}^\infty V_m$ and since each
$V_m$ is measurable, we conclude that $U_a$ is measurable.
This shows that $f+g$ is a measurable function and hence
the measurable functions form a vector space.

Clearly if $f$ and $g$ are bounded measurable functions and $c\in \R$
then $cf$ and $f+g$ are bounded.  We just showed they are also 
measurable, so the bounded measurable functions are a vector subspace.
\end{proof}

\begin{exer}\ 
\begin{enumerate}
\item Prove that if $f$ is a measurable function, then
so is $f^2.$  
\item Prove that if $f$ and $g$ are measurable functions, then
so is $fg.$  {\em Hint: } $2fg = (f+g)^2 - f^2 -g^2.$
\end{enumerate}
\end{exer}

\section{The Lebesgue Integral of Bounded Functions}

In this section we want to define the Lebesgue integral and
characterize the bounded integrable functions.  In the case of the
regulated integral, the integrable functions are the uniform limits of
step functions.  In the case of the Riemann integral a function $f$ is
integrable if the {\it infimum} of the integrals of step function
bigger than $f$ equals the {\it supremum} of the integrals of step
function less than $f.$ It is natural to alter both these definitions,
replacing step function with simple function.  It turns out that when
we do this for bounded functions we get the {\em same class of
integrable functions} whether we use the analog of regulated integral
or the analog of Riemann integral.  Moreover, this class is precisely
the bounded measurable functions!

\begin{thm}\label{thm: meas func}
If $f:[0,1] \to \R$ is a bounded function, then the following
are equivalent:
\begin{enumerate}
\item The function $f$ is Lebesgue measurable.
\item There is a sequence of simple functions $\{f_n\}_{n=1}^\infty$ 
which converges uniformly to $f$.
\item If $\U_{\mu}(f)$ denotes the set of all simple functions $u(x)$ such
that $f(x) \le u(x)$ for all $x$ and if $\yL_{\mu}(f)$ denotes 
the set of all simple
functions $v(x)$ such that $v(x) \le f(x)$ for all $x,$, then
\[
 \sup_{v \in \yL_{\mu}(f)} \Big\{\int v\ d\mu\Big\} =
\inf_{u \in \U_{\mu}(f)} \Big\{\int u\ d\mu\Big\}. 
\]
\end{enumerate}
\end{thm}

\begin{proof}
We will show 1) $\Rightarrow$ 2) $\Rightarrow$ 3) $\Rightarrow$ 1).
To show 1) $\Rightarrow$ 2), assume $f$ is a bounded measurable
function, say $a \le f(x) \le b$ for all $x \in [0,1].$

Let $\eps_n = (b-a)/n.$
We will partition the {\em range} $[a,b]$ of $f$ by intervals follows:
Let $c_i = a + i\eps_n$ so $a = c_0 < c_1 < \dots < c_n = b.$
Now define a measurable partition of $[0,1]$ by
$A_i = f^{-1}([c_{i-1}, c_i))$ for $i < n$ and $A_n =f^{-1}([c_{n-1}, b])$.
Then clearly $f_n(x) = \sum_{i=1}^n c_i \X_{A_i}$ is a simple function.
Moreover we note that for any $x \in [0,1]$ we have $|f(x) - f_n(x)| \le 
\eps_n.$  This is because $x$ must lie in one of the $A$'s, say
$x \in A_j.$ So $f_n(x) = c_j$ and  $f(x) \in [c_{j-1}, c_j[$. 
Hence $|f(x) - f_n(x)| \le c_j - c_{j-1} = \eps_n.$
This implies that the sequence of simple functions
$\{ f_n \}_{n=1}^\infty$ converges uniformly to $f$.

To show 2) $\Rightarrow$ 3), assume 
$f$ is the uniform limit of the sequence of simple functions 
$\{ f_n \}_{n=1}^\infty$.  This means if
$\delta_n = \sup_{x \in [0,1]} |f(x) - f_n(x)|$, then $\lim \delta_n = 0.$
We define simple functions 
$v_n(x) = f_n(x) -\delta_n$ and $u_n(x) = f_n(x) +\delta_n$ so
$v_n(x) \le f(x) \le u_n(x).$  

Then 
\begin{align}\label{eqn:ineq}
\inf_{u \in \U_\mu(f)} \Big\{\int u \ d\mu\Big\} &\le \liminf_{n \to \infty} 
\int u_n \ d\mu \notag\\
&= \liminf_{n \to \infty} \int (f_n + \delta_n)\ d\mu \notag\\
&= \liminf_{n \to \infty} \int f_n \ d\mu \notag\\
&\le \limsup_{n \to \infty} \int f_n \ d\mu \notag\\
&= \limsup_{n \to \infty} \int (f_n - \delta_n)\ d\mu \notag\\
&\le \limsup_{n \to \infty} \int v_n \ d\mu \notag\\
&\le \sup_{v \in \yL_\mu(f)} \Big\{\int v \ d\mu\Big\}.
\end{align}

For any $v \in \yL_\mu(f)$ and any $u \in \U_\mu(f)$ we have
$\int v \ d\mu \le \int u \ d\mu$ so
\[
\sup_{v \in \yL_\mu(f)} \Big\{\int v \ d\mu \Big\} 
\le \inf_{u \in \U_\mu(f)} \Big\{\int u \ d\mu \Big\}.
\]
Combining this with the inequality above we conclude that 
\begin{equation}\label{eqn:inf-sup}
\sup_{v \in \yL_\mu(f)} \Big\{\int v \ d\mu \Big\}
= \inf_{u \in \U_\mu(f)} \Big\{\int u \ d\mu\Big\}.
\end{equation}

All that remains is to show that 3) $\Rightarrow$ 1).  For this we
note that if 3) holds, then for any $n >0$ there are simple functions
$v_n$ and $u_n$ such that $v_n(x) \le f(x) \le u_n(x)$ for all $x$ and
such that 
\begin{equation}\label{eqn:$u_n$}
\int u_n\ d\mu - \int v_n\ d\mu < 2^{-n}.
\end{equation}

By Theorem (\ref{thm: lim}) the functions 
\[
g_1(x) = \sup_{n \in \N} \Big\{ v_n(x)\Big\} \text{ and }
g_2(x) = \inf_{n \in \N} \Big\{ u_n(x)\Big\}
\]
are measurable.  They are also bounded and satisfy
$g_1(x) \le f(x) \le g_2(x).$  We want to show that
$g_1(x) =  g_2(x)$ except on a set of measure zero,
which we do by contradiction.  Let $B = \{ x\ |\ g_1(x) < g_2(x) \}$
and suppose $\mu(B) > 0.$   Then since $B = \bigcup_{i=1}^\infty B_m$
where $B_m = \{ x\ |\ g_1(x) < g_2(x) -\frac{1}{m} \}$ we
conclude that $\mu( B_{m_0}) > 0$ for some $m_0.$
This implies that for every $n$ and every $x \in B_{m_0}$ we have 
$v_n(x) \le g_1(x) < g_2(x) -\frac{1}{m_0} \le u_n(x)-\frac{1}{m_0}.$
So $u_n(x) - v_n(x) > \frac{1}{m_0}$ for all $x \in B_{m_0}$ and
hence $u_n(x) - v_n(x) > \frac{1}{m_0}\X_{B_{m_0}}(x)$ for all $x$.
But this would mean that
$\int u_n\ d\mu - \int v_n\ d\mu = \int u_n - v_n\ d\mu 
\ge \int \frac{1}{m_0}\X_{B_{m_0}} \ d\mu = 
\frac{1}{m_0}\mu(B_{m_0})$ for all $n$ which contradicts
equation (\ref{eqn:$u_n$}) above.

Hence it must be the case that $\mu(B) = 0$  so $g_1(x) = g_2(x)$
except on a set of measure zero.  But since
$g_1(x) \le f(x) \le g_2(x)$ this means if we define 
$h(x) = f(x) - g_1(x)$, then $h(x)$ is zero except on a subset of 
$B$ which is a set of measure $0$.  It, then follows from 
Proposition (\ref{prop: meas zero func}) that $h$ is a measurable
function.  Consequently, $f(x) = g_1(x) + h(x)$ is also measurable
and we have completed the proof that 3) $\Rightarrow$ 1).
\end{proof}

\begin{defn}[Lebesgue Integral of a bounded function]
\label{def: lebesgue integral}\index{Lebesgue integral of a bounded function}
If $f : [0,1] \to \R$ is a bounded measurable function, then
we define its Lebesgue integral by
\[
\int f \ d\mu = \inf_{u \in \U_\mu(f)} \Big\{\int u \ d\mu\Big\},
\]
or equivalently (by Theorem (\ref{thm: meas func})),
\[
\int f \ d\mu = \sup_{v \in \yL_\mu(f)} \Big\{\int v \ d\mu\Big\}.
\]
\end{defn}

Alternatively, as the following proposition shows, we could have
defined it to be the limit of the integrals of a sequence of simple
functions converging uniformly to $f.$

\begin{prop}
If $\{g_n\}_n^\infty$
is any sequence of simple functions converging uniformly to a bounded
measurable function $f$, then
$\lim_{n \to \infty} \int g_n \ d\mu$ exists and is equal to
$\int f \ d\mu.$
\end{prop}

\begin{proof}
If we let
$\delta_n = \sup_{x \in [0,1]} |f(x) - g_n(x)|$, then $\lim \delta_n = 0$
and 
\[
g_n(x)-\delta_n \le f(x) \le g_n(x) + \delta_n.
\]
So $g_n -\delta_n  \in \yL_\mu(f)$ and $g_n + \delta_n \in \U_\mu(f)$.
Hence
\begin{align*}
\int f \ d\mu & = \inf_{u \in \U_\mu(f)} \int u \ d\mu\\
&\le \liminf_{n \to \infty} \int (g_n +\delta_n) \ d\mu \\
&=  \liminf_{n \to \infty} \int g_n \ d\mu \\
&\le \limsup_{n \to \infty} \int g_n \ d\mu\\
&\le \limsup_{n \to \infty} \int (g_n -\delta_n) \ d\mu\\
&\le \sup_{v \in \yL_\mu(f)} \Big\{\int v \ d\mu\Big\}\\
&= \int f \ d\mu.
\end{align*}

Hence these inequalities must be equalities and $\lim_{n \to \infty}
\int g_n \ d\mu = \int f\ d\mu.$
\end{proof}

\begin{exer}\ 
\begin{enumerate}
\item Since simple functions are themselves bounded measurable
functions, we have actually given two definitions of their Lebesgue
integral: the one in Definition (\ref{def: simple}) and the one above
in Definition (\ref{def: lebesgue integral}).  Prove that these
definitions give the same value.
\end{enumerate}
\end{exer}

\begin{thm}\label{thm: Lebesgue}
The Lebesgue integral, defined on
the vector space of bounded Lebesgue measurable functions on $[0,1]$,
satisfies the following properties:
\begin{description}
\item[I. Linearity:]
If $f$ and $g$ are Lebesgue measurable
functions and $c_1,c_2 \in \R$, then
\[
\int c_1f + c_2g\ d\mu = c_1 \int f\ d\mu + c_2 \int g\ d\mu.
\]
\item [II. Monotonicity:] If $f$ and $g$ are Lebesgue measurable
 and $f(x) \le g(x)$ for all $x$, then $\int f\ d\mu \le \int g\ d\mu.$
\item [III. Absolute value:] If $f$ is Lebesgue measurable 
then $|f|$ is also and $|\int f\ d\mu| \le \int |f|\ d\mu.$
\item [IV. Null Sets:] If $f$ and $g$ are bounded functions
and $f(x) = g(x)$ except on a set of measure zero, then $f$
is measurable if and only if $g$ is measurable. If they
are measurable, then $\int f\ d\mu = \int g\ d\mu.$
\end{description}
\end{thm}
\begin{proof}
If $f$ and $g$ are measurable there exist sequences of simple functions
$\{ f_n\}_{n=1}^\infty$ and $\{ g_n\}_{n=1}^\infty$ converging uniformly
to $f$ and $g$ respectively.  This implies that the sequence
$\{ c_1f_n + c_2g_n\}_{n=1}^\infty$ converges uniformly to 
the bounded measurable function $c_1f + c_2g.$ 
The fact that
\begin{align*}
\int c_1f + c_2g \ d\mu &= \lim_{n \to \infty} \int (c_1f_n + c_2g_n) \ d\mu \\
&= c_1\lim_{n \to \infty} \int f_n \ d\mu \ +
\ c_2\lim_{n \to \infty} \int g_n \ d\mu \\
&= c_1\int f \ d\mu\ +\ c_2\int g \ d\mu
\end{align*}
implies the linearity property.

Similarly the absolute value property follows from 
Lemma (\ref{lem: simple functions}) because
\[
\Big |\int f\ d\mu \Big | = \lim_{n \to \infty}\Big |\int f_n\ d\mu \Big |
\le \lim_{n \to \infty}\int |f_n|\ d\mu = \int |f|\ d\mu.
\]

To show monotonicity we use the definition of the Lebesgue integral.
If $f(x) \le g(x)$, then
\[
\int f \ d\mu = \sup_{v \in \yL_\mu(f)} \Big\{\int v \ d\mu\Big\}
\le \inf_{u \in \U_\mu(g)} \int u \ d\mu = \Big\{\int g \ d\mu\Big\}.
\]

If $f$ and $g$ are bounded functions which are equal except on a set
$E$ with $\mu(E) = 0$, then $h(x) = f(x) - g(x)$ is non-zero only on
the set $E$.  The function $h$ is measurable by Proposition
(\ref{prop: meas zero func}).  Clearly, since $f = g + h$ the function
$f$ is measurable if and only if $g$ is.

In case they are both measurable 
$|\int f\ d\mu - \int g\ d\mu| =|\int h\ d\mu| \le  \int |h|\ d\mu.$
But the function $h$ is bounded; say $|h(x)| \le M$.  Then
$|h(x)| \le M \X_E(x)$ so by monotonicity
$\int |h|\ d\mu \le \int M \X_E\ d\mu = M \mu(E) = 0.$
It follows that $|\int f\ d\mu - \int g\ d\mu| = 0$ so
$\int f\ d\mu = \int g\ d\mu.$
\end{proof}

\begin{defn} If $E \subset [0,1]$ is a measurable set and
$f$ is a bounded measurable function we define the Lebesgue integral
of $f$ over $E$ by
\[
\int_E f \ d\mu = \int f\X_E \ d\mu.
\]
\end{defn}

\begin{prop}[Additivity] \index{additivity}
If $E$ and $F$ are disjoint measurable subsets
of $[0,1]$, then 
\[
\int_{E \cup F} f \ d\mu = \int_E f \ d\mu + \int_F f \ d\mu.
\]
\end{prop}
\begin{proof}
If $E$ and $F$ are disjoint measurable subsets
of $[0,1]$, then $\X_{E \cup F} = \X_{E} + \X_{F}$ so 
\[
\int_{E \cup F} f \ d\mu = \int f\X_{E \cup F}\ d\mu
= \int f(\X_E + \X_F) \ d\mu = \int_E f \ d\mu + \int_F f \ d\mu.
\]
\end{proof}

\begin{prop}[Riemann integrable functions are Lebesgue integrable]
Every bounded Riemann integrable function $f:[0,1] \to \R$ is 
measurable and hence Lebesgue integrable.  The values of the
Riemann and Lebesgue integrals coincide.
\end{prop}

\begin{proof}
The set $\U(f)$ of step functions greater than $f$ is a subset of
the set $\U_\mu(f)$ of simple functions greater than $f$.  Likewise
the set $\yL(f) \subset \yL_\mu(f).$  Hence
\[
\sup_{v \in \yL(f)} \Big\{\int_0^1 v(t) \ dt\Big\} \le
\sup_{v \in \yL_\mu(f)} \Big\{ \int v \ d\mu \Big\}
\le \inf_{u \in \U_\mu(f)} \Big\{ \int u \ d\mu\Big\}
\le \inf_{u \in \U(f)} \Big\{ \int_0^1 u(t) \ dt\Big\}.
\]
The fact that $f$ is Riemann integrable asserts the first and
last of these values are equal.  Hence they are all equal and
$f$ is measurable and the Riemann and Lebesgue integrals
coincide.
\end{proof}

\section{The Bounded Convergence Theorem}

We want to investigate when the fact that a sequence of 
functions
$\{f_n\}_{n=1}^\infty$ converges pointwise to a function $f$ implies that
their Lebesgue integrals converge to the integral of $f$.
It is straightforward to prove that if a sequence of bounded
measurable functions converges uniformly to $f$, then their
integrals converge to the integral of $f$.  We will not do
this, because we prove a stronger result below.  But first
we consider an example which shows what can go wrong.

\begin{ex}
Let
\[
f_n(x) = 
\begin{cases}
	n, &\text{if $x \in [\frac{1}{n},\frac{2}{n}]$;} \\
	0, &\text{otherwise.} 
\end{cases}
\]
Then $f_n$ is a step function equal to $n$ on an interval of length
$\frac{1}{n}$ and $0$ elsewhere.  Thus $\int f_n \ d\mu = n\frac{1}{n}
= 1.$ But, for any $x \in [0,1]$ we have $f_n(x) = 0$ for all sufficiently
large $n$.  Thus the sequence
$\{f_n\}_{n=1}^\infty$ converges \myindex{pointwise} to the constant function
$0$.  Hence
\[
\int (\lim_{n \to \infty} f_n(x)) \ d\mu = 0
\text{ and } \lim_{n \to \infty} \int  f_n \ d\mu = 1.
\]
\end{ex}

In this example each $f_n$ is a bounded step function, but there is no
single bound which works for all $f_n$ since the maximum value of
$f_n$ is $n.$  It turns out that any example of this sort must be
a sequence of functions which is not uniformly bounded.

\begin{thm}[The Bounded Convergence Theorem]\label{thm: bounded convergence}
\index{Bounded Convergence Theorem}
Suppose $\{f_n\}_{n=1}^\infty$ is a sequence of measurable
functions which converges pointwise to a function $f$ and there is a
constant $M>0$ such that $|f_n(x)| \le M$ for all $n$ and all $x \in [0,1].$
Then $f$ is a bounded measurable function and 
\[
\lim_{n \to \infty} \int  f_n \ d\mu = \int  f \ d\mu.
\]
\end{thm}

\begin{proof}
For each $x\in [0,1]$ we know that $\lim_{m \to \infty} f_m(x) = f(x).$
This implies that $|f(x)| \le M$ and by Theorem (\ref{thm: lim}) that
$f(x)$ is measurable.

We must show that
\[
\lim_{n \to \infty} \Big | \int  f_n \ d\mu - \int  f \ d\mu \Big | = 0.
\]
but 
\begin{equation}\label{eqn:limit $f_n$}
\lim_{n \to \infty} \Big | \int  f_n \ d\mu - \int  f \ d\mu \Big | 
= \lim_{n \to \infty} \Big | \int  (f_n - f) \ d\mu \Big | 
\le \lim_{n \to \infty} \int  |f_n - f| \ d\mu.
\end{equation}
So we need to estimate the integral of $|f_n - f|.$

Given $\eps >0$ define $E_n = \{x\ |\ |f_m(x) - f(x)| <\eps/2 \text{
for all } m \ge n\}.$  Notice that if for some $n$ the set $E_n$ were all of
$[0,1]$ we would be able to estimate 
$\int |f_m-f| \ d\mu \le \int \eps/2 \ d\mu = \eps/2$
for all $m \ge n.$  But we don't know that.  Instead
we know that for any $x$ the limit $\lim_{m \to \infty} f_m(x) = f(x)$ which
means that each $x$ is in some $E_n$ (where $n$ depends on $x$).  In other
words $\bigcup_{n=1}^\infty E_n = [0,1].$

Since $E_n \subset E_{n+1}$ by Proposition
(\ref{prop: increasing union}) we know $\lim_{n \to \infty}\mu(E_n) =
\mu([0,1]) = 1.$ Thus there is an $n_0$ such that $\mu(E_{n_0}) > 1 -
\frac{\eps}{4M}$, so $\mu(E_{n_0}^c) < \frac{\eps}{4M}$.

Now for any $n > n_0$ we have
\begin{align*}
\int  |f_n - f| \ d\mu 
&= \int_{E_{n_0}}  |f_n - f| \ d\mu + \int_{E_{n_0}^c}  |f_n - f| \ d\mu\\
&\le \int_{E_{n_0}}  \frac{\eps}{2} \ d\mu + \int_{E_{n_0}^c}  2M \ d\mu\\
&\le \frac{\eps}{2}\mu(E_{n_0}) +  2M \mu(E_{n_0}^c)\\
&\le \frac{\eps}{2} +  2M \frac{\eps}{4M} = \eps.\\
\end{align*}
Thus we have shown $\lim_{n\to \infty} \int  |f_n - f| \ d\mu = 0.$
Putting this together with equation (\ref{eqn:limit $f_n$}) we see that
\[
\lim_{n \to \infty} \Big | \int  f_n \ d\mu - \int  f \ d\mu \Big | = 0
\]
as desired.
\end{proof}

\begin{defn}[Almost everywhere]
If a property holds for all $x$ except for a set of
measure zero, we say that it holds \myindex{almost everywhere} or
for \myindex{almost all} values of $x$.
\end{defn}
For example, we say that two functions $f$ and $g$ defined on $[0,1]$
are equal {\em almost everywhere} if the set of $x$ with $f(x) \ne
g(x)$ has measure zero.  The last part of Theorem (\ref{thm: Lebesgue})
asserted that if $f(x) = g(x)$ almost everywhere, then $\int f\ d\mu =
\int g\ d\mu.$ As another example, we say $\lim_{n \to \infty}f_n(x) =
f(x)$ for {\em almost all} $x$ if the set of $x$ where the limit does
not exist or is not equal to $f(x)$ is a set of measure zero.

\begin{thm}[Better Bounded Convergence Theorem]\label{thm: bbounded convergence}
Suppose $\{f_n\}_{n=1}^\infty$ is a sequence of bounded measurable functions
and $f$ is a bounded function such that
\[
\lim_{n \to \infty}f_n(x) = f(x)
\]
for almost all $x$. Suppose also there is a constant $M>0$ such that
for each $n >0$, $|f_n(x)| \le M$ for almost all $x \in [0,1]$.
Then $f$ is a measurable function, satisfying $|f(x)| \le M$ 
for almost all $x \in [0,1]$ and 
\[
\lim_{n \to \infty} \int  f_n \ d\mu = \int  f \ d\mu.
\]
\end{thm}

\begin{proof}
Let $A = \{x\ |\ \lim_{n \to \infty} f_n(x) \ne f(x)\}$, then
then $\mu(A) = 0$.  Define the set $D_n = \{x\ |\ |f_n(x)| > M\}$, then
then $\mu(D_n) = 0$  so if $E = A \cup \bigcup_{n=1}^\infty D_n$, then
$\mu(E) = 0.$ 
Let
\[
g_n(x) = f_n(x)\X_{E^c}(x) = 
\begin{cases}
	f_n(x), &\text{ if $x \notin E$;} \\
	0, &\text{ if $x \in E$.}
\end{cases}
\]
Then $|g_n(x)| \le M$ for {\em all} $x \in [0,1]$
and for any $x \notin E$ we have
$\lim_{n \to \infty} g_n(x) = \lim_{n \to \infty} f_n(x) = f(x).$
Also for $x \in E,\ g_n(x) = 0$ so so for all $x \in [0,1]$ we have
$\lim_{n \to \infty} g_n(x) = f(x)\X_{E^c}(x).$

Define the function $g$ by $g(x) = f(x)\X_{E^c}(x).$
For any $x \notin E$ we have
$g(x) = \lim_{n \to \infty} g_n(x) = \lim_{n \to \infty} f_n(x) = f(x)$
so $g(x) = f(x)$ almost everywhere.  We know from its definition
that $|g(x)| \le M$ since $|g_n(x)| \le M$.
And by Theorem (\ref{thm: lim}) $g$ is
measurable.  Since $f(x) - g(x)$ is zero almost everywhere it
is measurable by Proposition (\ref{prop: meas zero func}).  It
follows that $f$ is measurable and by  Theorem (\ref{thm: Lebesgue})
\[
\int f\ d\mu = \int g\ d\mu.
\]
Since $f_n = g_n$ almost everywhere we also know that
\[
\int f_n\ d\mu = \int g_n\ d\mu.
\]
Hence it will suffice to show that 
\[
\lim_{n \to \infty} \int g_n\ d\mu = \int g\ d\mu.
\]
But this is true by Theorem (\ref{thm: bounded convergence})
\end{proof}

\vfill\eject

\chapter{The Integral of Unbounded Functions}

In this section we wish to define and investigate the Lebesgue
integral of functions which are not necessarily bounded and even
extended real valued functions.  In fact, henceforth we will use the
term ``measurable function'' to refer to extended real valued
measurable functions.  If a function is unbounded both above and below
it is more complicated than if it is only unbounded above.  Hence we
first focus our attention on this case.

\section{Non-negative Functions}
\begin{defn}[Integrable Function]
\index{integrable function}
If $f: [0,1] \to \R$ is a non-negative Lebesgue measurable function
we let $f_n(x) = \min \{f(x), n\}.$  Then $f_n$ is a bounded
measurable function and we  define 
\[
\int f\ d\mu = \lim_{n \to \infty} \int f_n \ d\mu.
\]
If $\int f\ d\mu < \infty$ we say $f$ is {\em integrable.}
\end{defn}

Notice that the sequence $\{\int f_n\ d\mu\}_{n=1}^\infty$ is a monotonic
increasing sequence of numbers so the limit
$\lim_{n \to \infty} \int f_n \ d\mu$ either exists or is  $+\infty.$

\begin{prop}
If $f$ is a non-negative integrable function and
$A = \{x\ |\ f(x) = +\infty\}$, then $\mu(A) = 0.$
\end{prop}
\begin{proof}
For $x \in A$ we observe that $f_n(x) = n$ and hence 
$f_n(x) \ge n \X_A(x)$ for all $x$.  Thus
$\int f_n\ d\mu \ge \int n\X_A\ d\mu = n \mu(A).$
If $\mu(A) > 0$, then
$\int f\ d\mu = \lim \int f_n\ d\mu \ge \lim n \mu(A) = +\infty.$
\end{proof}

\begin{ex}
Let $f(x) = 1/\sqrt{x}$ for $x \in (0,1]$ and let $f(0) = +\infty.$  Then
$f:[0,1] \to \R$ is a non-negative measurable function.
Then the function
\[
f_n(x) = 
\begin{cases}
	n, &\text{ if $0 \le x < \frac{1}{n^2}$;} \\
	\frac{1}{\sqrt{x}}, &\text{ if $ \frac{1}{n^2} \le x \le 1$.}
\end{cases}
\]
Hence if $E_n = [0,1/n^2[$, then 
\begin{align*}
\int f_n\ d\mu &= \int_{E_n} f_n\ d\mu + \int_{E_n^c} f_n\ d\mu\\
&= \int n\X_{E_n} \ d\mu + \int_{\frac{1}{n^2}}^1 \frac{1}{\sqrt{x}}\ dx\\
&= n\mu(E_n) + \Big (2 - \frac{2}{n} \Big )\\
&= \frac{n}{n^2} + 2 - \frac{2}{n} = 2-\frac{1}{n}.
\end{align*}
Hence 
\[
\int f\ d\mu = \lim_{n \to \infty} \int f_n \ d\mu = 2.
\]
So $f$ is integrable.
\end{ex}

\begin{prop}\label{prop: integral monotonicity}
Suppose $f$ and $g$ are non-negative measurable functions with
$g(x) \le f(x)$ for almost all $x$.  
If $f$ is integrable, then $g$ is integrable
and $\int g\ d\mu \le \int f\ d\mu.$  In particular if $g = 0$
almost everywhere, then $\int g\ d\mu = 0.$
\end{prop}

\begin{proof}
If $f_n(x) = \min\{f(x),n\}$ and $g_n(x) = \min\{g(x),n\}$, then
$f_n$ and $g_n$ are bounded measurable functions and satisfy
$g_n(x) \le f_n(x)$ for almost all $x$.  It follows that
$\int g_n\ d\mu \le \int f_n\ d\mu \le \int f\ d\mu.$
Since the sequence of numbers $\{\int g_n\ d\mu\}_{n=1}^\infty$
is monotonic increasing and bounded above by $\int f\ d\mu$ it has 
a finite limit.  By definition this limit is $\int g\ d\mu.$
Since for each $n$ we have $\int g_n\ d\mu \le \int f\ d\mu$,
the limit is also bounded by $\int f\ d\mu$.  That is,
\[
\int g\ d\mu =\lim_{n\to \infty}\int g_n\ d\mu \le \int f\ d\mu.
\]
If $g = 0$ almost everywhere, then $0 \le g(x) \le 0$ for almost
all $x$ so we have $\int g\ d\mu = 0.$
\end{proof}

\begin{cor}\label{cor: zero ae}
If $f: [0,1] \to \R$ is a
non-negative integrable function and $\int f\ d\mu = 0$
then $f(x) = 0$ for almost all $x.$
\end{cor}
\begin{proof}
Let $E_n = \{x\ |\ f(x) \ge 1/n\}$.  Then 
$f(x) \ge \frac{1}{n}\X_{E_n}(x)$ so
\[
\frac{1}{n}\mu(E_n) = \int \frac{1}{n}\X_{E_n}\ d\mu \le \int f\ d\mu = 0.
\]
Hence $\mu(E_n) = 0.$  But if $E = \{x\ |\ f(x) > 0 \}$, then
$E = \bigcup_{n=1}^\infty E_n$ so $\mu(E) = 0.$
\end{proof}

\begin{thm}[Absolute Continuity]\label{thm: abs continuity}
\index{absolute continuity of measure}
Suppose $f$ is a non-negative integrable function.  Then for
any $\eps >0$ there exists a $\delta >0$ such that
$\int_A f\ d\mu < \eps$ for every measurable $A \subset [0,1]$
with $\mu(A) < \delta.$
\end{thm}

\begin{proof}
Let $f_n(x) = \min\{ f(x), n\}$ so $\lim \int f_n \ d\mu = \int f\ d\mu.$
Let 
\[
E_n = \{ x \in [0,1]\ |\ f(x) \ge n\}
\]
so 
\[
f_n(x) =
\begin{cases}
	n, &\text{ if $x \in E_n$;} \\
	f(x), &\text{ if $x \in E_n^c$.}
\end{cases}
\]
Consequently
\[
\int f_n\ d\mu =  \int_{E_n} n\ d\mu + \int_{E_n^c} f\ d\mu.
\]

Hence we have
\begin{align*}
\int f\ d\mu &= \int_{E_n} f\ d\mu + \int_{E_n^c} f\ d\mu\\
&= \int_{E_n} (f - n)\ d\mu + \int_{E_n} n\ d\mu + \int_{E_n^c} f\ d\mu\\
&= \int_{E_n} (f - n)\ d\mu + \int f_n\ d\mu.
\end{align*}
Thus $\int f\ d\mu - \int f_n\ d\mu = \int_{E_n} (f - n)\ d\mu$ and
we conclude from integrability of $f$ that
\[
\lim_{n \to \infty} \int_{E_n} (f - n)\ d\mu = 0.
\]

Hence we may choose $N$ such that 
$\int_{E_N} (f - N)\ d\mu < \eps/2.$  Now pick $\delta < \eps/2N.$
Then if $\mu(A) < \delta$ we have
\begin{align*}
\int_A f\ d\mu &= \int_{A \cap E_N} f\ d\mu + \int_{A \cap E_N^c} f\ d\mu\\
&\le  \int_{A\cap E_N} (f -N)\ d\mu + \int_{A\cap E_N} N\ d\mu + \int_{A \cap E_N^c} N\ d\mu\\
&\le  \int_{E_N} (f - N)\ d\mu  + \int_{A} N\ d\mu\\
&<  \frac{\eps}{2}  + N \mu(A) < \frac{\eps}{2} + N \delta < \eps.
\end{align*}
\end{proof}

Theorem (\ref{thm: abs continuity}) is labeled ``Absolute Continuity''
for reasons that will become clear later in Section \S\ref{sec: other measures}.
But as a nearly immediate consequence we have the following generalization
of a result from Exercise (\ref{exer: integral continuous}).

\begin{cor}[Continuity of the Integral]
If $f:[0,1] \to \R$ is a non-negative integrable function and
we define $F(x) = \int_{[0,x]} f\ d\mu$, then $F(x)$ is continuous.
\end{cor}

\begin{proof}
Given $\eps >0$ let $\delta >0$ be the corresponding value guaranteed by
Theorem (\ref{thm: abs continuity}).  Now suppose $x < y$ and 
$|y - x| < \delta$.  Then $\mu([x,y]) < \delta$ so
\[
\big | F(y) - F(x) \big | 
= \Big | \int_{[0,y]} f\ d\mu - \int_{[0,x]} f\ d\mu\Big | 
= \Big | \int_{[x,y]} f\ d\mu \Big | < \eps
\]
by Theorem (\ref{thm: abs continuity}).
We have in fact proven that $F$ is uniformly continuous.
\end{proof}

\begin{exer}\label{exer: Lebesgue thm}\ 
\begin{enumerate}
\item Define $f(x) = \frac{1}{x^p}$ for $x \in (0,1]$ and $f(0) = +\infty.$ 
Prove that $f$ is integrable if and only if $p < 1$.  Calculate the
value of $\int f\ d\mu$ in this case.
\item Give an example of a non-negative extended function $g: [0,1] \to \R$ which
is integrable and which has the value $+\infty$ at infinitely many points of
$[0,1].$
\end{enumerate}
\end{exer}

\section{Convergence Theorems}\label{sec: convergence theorems}

The following result is very similar to the Bounded Convergence
Theorem (see Theorem (\ref{thm: bounded convergence}) and 
Theorem (\ref{thm: bbounded convergence}).  The difference is that
instead of having a constant  bound on the functions $f_n$ we have
them bounded by an integrable function $g.$  This is enough to make
essentially the same proof work, however, because of 
Theorem (\ref{thm: abs continuity}).

\begin{thm}[Lebesgue Convergence for Non-negative functions]
\index{Lebesgue Convergence Theorem}
\label{thm: lebesgue convergence+}
Suppose $f_n$ is a sequence of non-negative measurable functions and $g$ is
a non-negative integrable function such that $f_n(x) \le g(x)$ for
all $n$ and almost all $x$.  If\ $\lim f_n(x) = f(x)$ for almost all
$x$, then $f$ is integrable and
\[
\int f\ d\mu = \lim_{n \to \infty} \int f_n\ d\mu.
\]
\end{thm}
\begin{proof}
If we let $h_n = f_n \X_E$ and
 $h = f \X_E$ where $E= \{x\ |\ \lim f_n(x) = f(x)\}$, then
$f=h$ almost everywhere and $f_n = h_n$ almost everywhere.
So it suffices to prove 
\[
\int h\ d\mu = \lim_{n \to \infty} \int h_n\ d\mu,
\]
and we now have the stronger property that
$\lim h_n(x) = h(x)$ for {\em all} $x$, instead of almost all.
Since $h_n(x) = f_n(x) \X_E(x) \le g(x)$ for almost all
$x$ we know that $h(x) \le g(x)$ for almost all $x$ and
hence by Proposition (\ref{prop: integral monotonicity})
that $h$ is integrable.

The remainder of the proof is very similar to the proof
of Theorem (\ref{thm: bounded convergence}).
We must show that
\[
\lim_{n \to \infty} \Big | \int  h_n \ d\mu - \int  h \ d\mu \Big | = 0.
\]
but 
\begin{equation}\label{eqn:limit $h_n$}
\lim_{n \to \infty} \Big | \int  h_n \ d\mu - \int  h \ d\mu \Big | 
= \lim_{n \to \infty} \Big | \int  (h_n - h) \ d\mu \Big | 
\le \lim_{n \to \infty} \int  |h_n - h| \ d\mu.
\end{equation}
So we need to estimate the integral of $|h_n - h|.$

Given $\eps > 0$ 
define $E_n = \{x\ |\ |h_m(x) - h(x)| <\eps/2 \text{
for all } m \ge n\}.$
We know by Theorem (\ref{thm: abs continuity}) that
there is a $\delta >0$ such that $\int_A g\ d\mu < \eps/4$ whenever
$\mu(A) < \delta.$ 

We also know that for any $x$ the limit $\lim_{m \to \infty} h_m(x) = h(x)$ which
means that each $x$ is in some $E_n$ (where $n$ depends on $x$).  In other
words $\bigcup_{n=1}^\infty E_n = [0,1].$
Since $E_n \subset E_{n+1}$ by Proposition
(\ref{prop: increasing union}) we know $\lim_{n \to \infty}\mu(E_n) =
\mu([0,1]) = 1.$ Thus there is an $n_0$ such that $\mu(E_{n_0}) > 1 -
\delta$, so $\mu(E_{n_0}^c) < \delta$.

Now  $|h_n(x) - h(x)| \le |h_n(x)| + |h(x)|
\le 2 g(x)$ so for any $n > n_0$ we have
\begin{align*}
\int  |h_n - h| \ d\mu 
&= \int_{E_{n_0}}  |h_n - h| \ d\mu + \int_{E_{n_0}^c}  |h_n - h| \ d\mu\\
&\le \int_{E_{n_0}}  \frac{\eps}{2} \ d\mu + \int_{E_{n_0}^c}  2g \ d\mu\\
&\le \frac{\eps}{2}\mu(E_{n_0}) +  2 \int_{E_{n_0}^c}  g \ d\mu\\
&\le \frac{\eps}{2} +  2 \frac{\eps}{4} = \eps.\\
\end{align*}
Thus we have shown $\displaystyle{\lim_{n\to \infty} \int  |h_n - h| \ d\mu = 0
.}$ Putting this together with equation (\ref{eqn:limit $h_n$}) we see that
\[
\lim_{n \to \infty} \Big | \int  h_n \ d\mu - \int  h \ d\mu \Big | = 0
\]
as desired.

\end{proof}

\begin{thm}[Monotone Convergence Theorem]\label{thm: monotone convergence}
\index{Monotone Convergence Theorem}
Suppose $g_n$ is an increasing sequence of non-negative measurable functions.
If\ $\lim g_n(x) = f(x)$ for almost all
$x$, then
\[
\int f\ d\mu = \lim_{n \to \infty} \int g_n\ d\mu.
\]
In particular $f$ is integrable if and only if $\lim \int g_n\ d\mu < +\infty.$
\end{thm}

\begin{proof}
The function $f$ is measurable by Theorem (\ref{thm: lim}).  If it is
integrable, then the fact that $f(x) \ge g_n(x)$ for almost all $x$
allows us to apply the previous theorem to conclude the desired
result.

Hence we need only show that if $\int f\ d\mu = +\infty$, then
$\displaystyle{\lim_{n \to \infty} \int g_n\ d\mu = +\infty}.$
But if $\int f\ d\mu = +\infty$, then for any $N >0$ there exists
$n_0$ such that $\int \min \{f(x), n_0\}\ d\mu > N$.
And we know that
\[
\lim_{n\to \infty}\min \{g_n(x), n_0\} = \min \{f(x), n_0\}
\]
for almost all $x$.  Since these are bounded measurable functions,
\[
\lim_{n\to \infty}\int g_n\ d\mu \ge
\lim_{n\to \infty}\int \min \{g_n, n_0\}\ d\mu = \int \min \{f, n_0\}\ d\mu > N,
\]
where the equality comes from the bounded convergence 
Theorem~(\ref{thm: bbounded convergence}).
Since $N$ is arbitrary we conclude that
\[
\lim_{n\to \infty}\int g_n\ d\mu = +\infty.
\]
\end{proof}

\begin{cor}[Integral of infinite series]\label{cor: infinite series}
\index{infinite series}
Suppose $u_n$ is a non-negative measurable function and $f$ is
an non-negative function such that $\sum_{n=1}^\infty u_n(x) = f(x)$
for almost all $x$.  Then
\[
\int f\ d\mu = \sum_{n =1}^\infty \int u_n\ d\mu.
\]
\end{cor}
\begin{proof}
Define
\[
f_N(x) = \sum_{n=1}^N u_n(x).
\]
Now the result follows from the previous theorem.
\end{proof}

\section{Other Measures}\label{sec: other measures}

There are other measures besides Lebesgue and indeed measures
on other spaces besides $[0,1]$ or $\R$.  We will limit our attention 
to measures defined on $I = [0,1]$.

Recall that a collection $\A$ of subsets of $I$ is called a $\sigma$-algebra
provided it contains the set $I$ and is closed under taking
complements, countable unions, and countable intersections.

\begin{exs} The following are examples of $\sigma$-algebras on $I=[0,1]$:
\begin{enumerate}
\item {\em The trivial $\sigma$-algebra.} $\A = \{ \emptyset, I\}.$
\item  $\A = \{ A \subset I\ |\ A \text{ is countable, or }
A^c\text{ is countable}\}.$
\item $\A = \M$ the Lebesgue measurable sets
\item $\A$ is Borel sets, the smallest $\sigma$-algebra containing 
the open intervals.
\end{enumerate}
\end{exs}

\begin{defn}[Finite Measure]\label{def: measure}
If $\A$ is a $\sigma$-algebra of subsets of $I$, then a function
$\nu: \A \to \R$ is called a \myindex{finite measure} provided 
\begin{itemize}
\item $\nu(A) \ge 0$ for every $A \in \A,$
\item $\nu(\emptyset) = 0,$ $\nu(I) < \infty,$ and
\item $\nu$ is countably additive, i.e. if $\{A_n\}_{n=1}^\infty$
are pairwise disjoint sets in $\A$, then 
\[
\nu(\bigcup_{n=1}^\infty A_n) = \sum_{n=1}^\infty \nu(A_n).
\]
\end{itemize}
\end{defn}

We will restrict our attention to measures defined on the
$\sigma$-algebra of Lebesgue measurable sets.  The integral of a
measurable function with respect to a measure $\nu$ is defined
analogously to Lebesgue measure.

\begin{defn} Let $\nu$ be a finite measure defined on the $\sigma$-algebra
$\M(I)$.  If $f(x) = \sum_{i=1}^n r_i \X_{A_i}$ is a simple function
then its \myindex{integral with respect to $\nu$} is
defined by $\int f\ d\nu = \sum_{i=1}^n r_i \nu(A_i).$
If $g : [0,1] \to \R$ is a bounded measurable function, then
we define its integral with respect to $\nu$ by  
\[
\int g \ d\nu = \inf_{u \in \U_\mu(g)}  \Big\{ \int u \ d\nu \Big\}.
\]
If $h$ is a non-negative extended measurable function we define
\[
\int h \ d\nu = \lim_{n \to \infty}  \int \min\{h,n\} \ d\nu.
\]
\end{defn}

\begin{defn}[Absolutely Continuous Measure]\label{def: abs cont}
\index{absolutely continuous measure}
If $\nu$ is a measure defined on $\M(I),$ the Lebesgue
measurable subsets of $I,$, then we say $\nu$ is 
\myindex{absolutely continuous} with respect to Lebesgue measure
$\mu$ if $\mu(A) = 0$ implies $\nu(A) = 0.$
\end{defn}

The following result motivates the name ``absolute continuity.''

\begin{thm}\label{thm: abs cont}
If $\nu$ is a measure defined on $\M(I)$
which is absolutely continuous
with respect to Lebesgue measure, then for any $\eps > 0$
there is a $\delta >0$ such that 
$\nu(A) < \eps$ whenever $\mu(A) = \delta.$
\end{thm}

\begin{proof}
We assume there is a counter-example and show this leads to a 
contradiction.  If the measure $\nu$ does not satisfy the 
conclusion of the theorem, then there is an $\eps >0$ for which it
fails, i.e. there is no $\delta > 0$ which works for this $\eps.$
In particular, for any positive integer $m$
there is a set $B_m$ such that $\nu(B_m) \ge \eps$ 
and $\mu(B_m) < 1/2^{m}.$
Hence if we define $A_n = \bigcup_{m=n+1}^\infty B_m $, then
\[
\mu(A_n) \le \sum_{m=n+1}^\infty \mu(B_m) \le \sum_{m=n+1}^\infty \frac{1}{2^{m}}
= \frac{1}{2^{n}.}
\]
The sets $A_n$ are nested, i.e. $A_n \supset A_{n+1}.$  It follows from 
Proposition (\ref{prop: increasing union}) that 
\begin{equation}\label{eqn:mu $A_n$}
\mu(\bigcap_{n=1}^\infty A_n) = \lim_{n \to \infty} \mu(A_n) 
\le \lim_{n \to \infty}\frac{1}{2^{n}} = 0.
\end{equation}

The proof of Proposition (\ref{prop: increasing union}) made use only
of the countable additivity of the measure.  Hence it is also valid
for $\nu$, i.e. 
\[
\nu(\bigcap_{n=1}^\infty A_n) = \lim_{n \to \infty} \nu(A_n).
\]
On the other hand $\nu(A_n) \ge \nu(B_{n+1}) \ge \eps,$ so
\[
\nu(\bigcap_{n=1}^\infty A_n) = \lim_{n \to \infty} \nu(A_n) 
\ge \lim_{n \to \infty} \eps = \eps.
\]
This together with equation (\ref{eqn:mu $A_n$}) 
contradicts the absolute continuity of $\nu$ with respect to $\mu.$
We have proven the contrapositive of the result we desire.
\end{proof}

\begin{exer}
Given a point $x_0 \in [0,1]$ define the function 
$\delta_{x_0} : \M \to \R$ by
$\delta_{x_0}( A) = 1$ if $x_0 \in A$ and $\delta_{x_0}( A) = 0$ 
if $x_0 \notin A$. Let $\nu(A) = \delta_{x_0}( A).$
\begin{enumerate}
\item Prove that $\nu$ is a measure. 
\item Prove that if $f$ is a measurable function
$\int f \ d\nu = f(x_0).$
\item Prove that $\nu$ is not absolutely continuous
with respect to Lebesgue measure $\mu.$
\end{enumerate}
The measure $\nu$  is called the \myindex{Dirac $\delta$-measure.}
\end{exer}

\begin{prop}\label{prop: $nu_f$}
If $f$ is a non-negative integrable function on $I$ and we
define 
\[
\nu_f(A) = \int_A f\ d\mu
\]
then $\nu_f$ is a measure with $\sigma$-algebra $\M(I)$ which is
absolutely continuous with respect to Lebesgue measure $\mu$.
\end{prop}

\begin{proof}
Clearly $\nu_f(A) = \int_A f\ d\mu \ge 0$ for all $A \in \M$ since
$f$ is non-negative.  Also $\nu_f(\emptyset) = 0.$  We need to
check countable additivity.  

Suppose $\{A_n\}_{n=1}^\infty$ is
a sequence of pairwise disjoint measurable subsets of $[0,1].$
and $A$ is their union.  Then for all $x \in [0.1].$
\[
f(x)\X_{A}(x) = \sum_{n=1}^\infty f(x) \X_{A_n}(x).
\]
Hence by Theorem (\ref{cor: infinite series}) 
\[
\int f\X_{A} \ d\mu  = \sum_{n=1}^\infty \int f \X_{A_n}\ d\mu
\]
and so 
\[
\nu_f(A) = \sum_{n=1}^\infty \nu_f (A_n).
\]
Thus $\nu$ is a measure.

If $\mu(A) = 0$, then by Proposition (\ref{prop: integral monotonicity})
\[
\nu_f(A) = \int f\X_A \ d\mu =0 
\]
so $\nu$ is absolutely continuous with respect to $\mu.$
\end{proof}

The converse to Proposition (\ref{prop: $nu_f$}) is called the 
Radon-Nikodym Theorem.  Its proof is beyond the scope of this
text.  A proof can be found in Chapter 11 Section 5 of Royden's
book \cite{R}

\begin{thm}[Radon-Nikodym]\label{Radon-Nikodym}
\index{Radon-Nikodym Theorem} If $\nu$ is a measure with
\myindex{$\sigma$-algebra} $\M(I)$ which is absolutely
continuous with respect to Lebesgue measure $\mu$, then there is a
non-negative integrable function $f$ on $[0,1]$ such that define
\[
\nu(A) = \int_A f\ d\mu.
\]
The function $f$ is unique up to measure $0$, i.e. if $g$ is another
function with these properties, then $f = g$ almost everywhere.
\end{thm}

The function $f$ is called the \myindex{Radon-Nikodym derivative}
of $\nu$ with respect to $\mu.$  In fact the Radon-Nikodym Theorem
is more general than we have stated, since it applies to any
two finite measures $\nu$ and $\mu$ defined on a $\sigma$-algebra
$\A$ with $\nu$ absolutely continuous with respect to $\mu.$

\section{General Measurable Functions}\label{sec: general measurable}

In this section we consider extended measurable functions which may be unbounded 
both above and below. 
we will define
\[
f^+(x) = \max \{ f(x), 0 \} \text{ and } f^-(x) = - \min \{ f(x), 0 \}.
\]
These are both non-negative measurable functions.

\begin{defn}\label{def: general integral}
If $f: [0,1] \to \R$ is a measurable function, then we say $f$
is \myindex{Lebesgue integrable} provided both $f^+$ and $f^-$ are
integrable (as non-negative functions). If $f$ is integrable
we define 
\[
\int f \ d\mu = \int f^+ \ d\mu - \int f^- \ d\mu.
\]
\end{defn}

\begin{prop}\label{prop: almost everywhere}
Suppose $f$ and $g$ are measurable functions on $[0,1]$ and $f = g$
almost everywhere.  Then if $f$ is integrable, so is $g$ and
$\int f\ d\mu = \int g\ d\mu.$  In
particular if $f = 0$ almost everywhere $\int f\ d\mu = 0.$
\end{prop}
\begin{proof}
If $f$ and $g$ are measurable functions on $[0,1]$ and $f = g$
almost everywhere, then  $f^+ = g^+$ almost everywhere,
$f^- = g^-$ almost everywhere, and $f^+$ and $f^-$ are integrable.
It, then follows from Proposition (\ref{prop: integral monotonicity})
that $g^+$ and $g^-$ are integrable and that
$\int f^+\ d\mu \ge \int g^+\ d\mu$ and $\int f^-\ d\mu \ge \int g^-\ d\mu$.
Switching the roles of $f$ and $g$ this same proposition gives the
reverse inequalities so we have
$\int f^+\ d\mu = \int g^+\ d\mu$ and $\int f^-\ d\mu = \int g^-\ d\mu$.
\end{proof}

\begin{prop}\label{prop: integral abs}
The measurable function  $f: [0,1] \to \R$ is integrable if and only
if the the function $|f|$ is integrable.
\end{prop}

\begin{proof}
Notice that $|f(x)| = f^+(x) + f^-(x).$   Thus if $|f|$ is integrable,
since $|f(x)| \ge f^+(x)$ and $|f(x)| \ge f^-(x)$ it follows from
Proposition (\ref{prop: integral monotonicity}) that both 
$f^+$ and $f^-$ are integrable.  Conversely if $f^+$ and $f^-$ are integrable
then so is their sum $|f|.$
\end{proof}

\begin{thm}[Lebesgue Convergence Theorem]
\index{Lebesgue Convergence Theorem}
\label{thm: general lebesgue convergence}
Suppose $f_n$ is a sequence of measurable functions and $g$ is
an non-negative integrable function such that $|f_n(x)| \le g(x)$ for
all $n$ and almost all $x$.  If $\lim f_n(x) = f(x)$ for almost all
$x$, then $f$ is integrable and
\[
\int f\ d\mu = \lim_{n \to \infty} \int f_n\ d\mu.
\]
\end{thm}
\begin{proof}
The functions $f_n^+(x) = \max\{ f_n(x), 0\}$ and 
$f_n^-(x) = -\min\{ f_n(x), 0\}$ satisfy 
\[
\lim_{n \to \infty} f_n^+(x) = f^+(x) \text{ and } 
\lim_{n \to \infty} f_n^-(x) = f^-(x) 
\]
for almost all $x$.  Also $g(x) \ge f_n^+(x)$ and
$g(x) \ge f_n^-(x)$ for almost all $x$.  Hence by Theorem 
(\ref{thm: lebesgue convergence+})
\[
\int f^+\ d\mu = \lim_{n \to \infty} \int f_n^+\ d\mu \text{ and }
\int f^-\ d\mu = \lim_{n \to \infty} \int f_n^-\ d\mu.
\]
Thus $f$ is integrable and
\begin{align*}
\int f\ d\mu &= \int f^+\ d\mu - \int f^-\ d\mu\\
&=\lim_{n \to \infty} \int f_n^+\ d\mu - \lim_{n \to \infty} \int f_n^-\ d\mu\\
&=\lim_{n \to \infty} \int f_n^+ - f_n^-\ d\mu \\
&=\lim_{n \to \infty} \int f_n\ d\mu.
\end{align*}
\end{proof}

The following theorem says that for any $\eps >0$ any integrable
function can be approximated within $\eps$ by a step function if we
are allowed to exclude a set of measure $\eps.$

\begin{thm}\label{thm: step approx}
If $f:[0,1] \to \R$ is an integrable function, then given $\eps >0$
there is a step function $g:[0,1] \to \R$ 
and a  measurable subset $A \subset [0,1]$ such that 
$\mu(A) < \eps$ and 
\[
|f(x) - g(x)| < \eps \text{ for all } x \notin A.
\]
Moreover, if $|f(x)| \le M$ for all $x$, then we may choose $g$
with this same bound.
\end{thm}
\begin{proof}
We first prove the result for the special case of $f(x) = \X_E(x)$
for some measurable set $E.$  This follows because there is a 
countable cover of $E$ by open intervals $\{U_i\}_{i=1}^\infty$
such that 
\[
\mu(E) \le \sum_{i=1}^\infty \len(U_i) \le  \mu(E)+ \frac{\eps}{2}.
\]
and hence
\begin{equation}\label{eqn: step approx 1}
\mu \Big( \big( \bigcup_{i=1}^\infty U_i  \big) \setminus E\Big) 
< \frac{\eps}{2}.
\end{equation}
Also we may choose $N>0$ such that
\begin{equation}\label{eqn: step approx 2}
\mu\Big( \bigcup_{i=N}^\infty U_i \Big) \le
\sum_{i=N}^\infty \len(U_i) < \frac{\eps}{2}.
\end{equation}

Let $V_N = \cup_{i=1}^N U_i$.  It is a finite union of intervals,
so the function $g(x) = \X_{V_N}$ is a step function 
and if $A = \{x\ |\ f(x) \ne g(x)\}$, then 
\[
A \subset \Big( V_N \setminus E \Big)
\cup \Big( E \setminus V_N \Big) \subset
\Big( \big( \bigcup_{i=1}^\infty U_i  \big) \setminus E\Big)
\cup \Big( \bigcup_{i=N}^\infty U_i \Big),
\]
so it follows from equations (\ref{eqn: step approx 1}) and (\ref{eqn:
step approx 2}) that $\mu(A) < \eps.$ This proves the result for $f =
\X_E$.

From this the result follows for simple functions $f= \sum r_i\X_{E_i}$
because if $g_i$ is the approximating step function for $\X_{E_i}$
then $g= \sum r_ig_i$ approximates $f$ (with a suitably adjusted $\eps$).

If $f$ is a bounded measurable function by Theorem (\ref{thm: meas func})
there is a simple function $h$ such that $|f(x) - h(x)| < \eps/2$ for
all $x.$  Let $g$ be a step function such that $|h(x) - g(x)| < \eps/2$ for
all $x \notin A$ with $\mu(A) < \eps.$    Then
\[
|f(x) - g(x)| \le  |f(x) - h(x)| + |h(x) - g(x)| 
< \frac{\eps}{2} + \frac{\eps}{2} = \eps,
\]
for all $x \notin A.$ That is, the result is true if $f$ is a 
bounded measurable function.

Suppose $f$ is a non-negative integrable function.
Let $A_n = \{ x\ |\ f(x) > n\}$.  Then
\[
n \mu(A_n) =  \int n\X_{A_n} \ d\mu \le \int f\ d\mu < \infty.
\]
It follows that $\lim \mu(A_n) = 0.$  Hence
there is an $N>0$ such that $\mu(A_N) < \eps/2.$  

If $f_N = \min \{f,N\}$, then $f_N$ is a bounded measurable function.
So we may
choose a step function $g$ such that  $|f_N(x) - g(x)| < \eps/2$ for
all $x \notin B$ with $\mu(B) < \eps/2.$  It follows that
if $A = A_N \cup B$, then $\mu(A) < \eps.$  Also
if $x \notin A$, then $f(x) = f_N(x)$ so
\[
|f(x) - g(x)| \le  |f(x) - f_N(x)| + |f_N(x) - g(x)| = |f_N(x) - g(x)| <\eps.
\]
Hence the result holds for non-negative $f$.

For a general integrable $f$ we have $f = f^+ - f^-$.  The fact
that the result holds for $f^+$ and  $f^-$ easily implies it holds for
$f.$

Suppose now that $f$ is bounded, say $|f(x)| \le M$ for all $x$ 
and $g$ satisfies the conclusion of our theorem, then we define 
\[
g_1 (x) =
\begin{cases}
M, \text{ if } g(x) > M;\\
g(x), \text{ if } -M \le g(x) \le M;\\
-M \text{ if } g(x) < -M.
\end{cases}
\]
The function $g_1$ is a step function 
with $|g_1(x)| \le M$ and $g_1(x) = g(x)$ except
when $|g(x)| > M.$   Note if $g(x) > M$ and  $x \notin A$
then $f(x) \le M = g_1(x) < g(x)$
so $|g_1(x) - f(x) | < \eps.$   The case $g(x) < -M$ is similar.

\end{proof}

\begin{thm}\label{thm: Lebesgue general}
The Lebesgue integral 
satisfies the following properties:
\begin{description}
\item[I. Linearity:]
If $f$ and $g$ are Lebesgue measurable
functions and $c_1,c_2 \in \R$, then
\[
\int c_1f + c_2g\ d\mu = c_1 \int f\ d\mu + c_2 \int g\ d\mu.
\]
\item [II. Monotonicity:] If $f$ and $g$ are Lebesgue measurable
 and $f(x) \le g(x)$ for all $x$, then $\int f\ d\mu \le \int g\ d\mu.$
\item [III. Absolute value:] If $f$ is Lebesgue measurable 
then $|f|$ is also and $|\int f\ d\mu| \le \int |f|\ d\mu.$
\item [IV. Null Sets:] If $f$ and $g$ are bounded functions
and $f(x) = g(x)$ except on a set of measure zero, then $f$
is measurable if and only if $g$ is measurable. If they
are measurable, then $\int f\ d\mu = \int g\ d\mu.$
\end{description}
\end{thm}

The proof is left as an exercise.

\begin{exer}\label{exer: f=g=h ae}\ 
\begin{enumerate}
\item Prove that if $f,g,h$ are measurable functions and $f=g$ almost
everywhere and $g=h$ almost everywhere, then $f=h$ almost everywhere.
\item Prove that if $f:[0,1] \to \R$ is an integrable function, then
given $\eps > 0$ there exists a continuous function $g:[0,1] \to \R$
and a set $A$ such that $\mu(A) < \eps,\ |f(x) - g(x)| < \eps$
for all $x \notin A,$ and $g(0) = g(1).$
\item Prove that the, not necessarily bounded, integrable functions
from $[0,1]$ to $\R$ form a vector space.
\item Prove Theorem (\ref{thm: Lebesgue general}).
Proposition (\ref{prop: almost everywhere}) proves the null
set property.  Prove the remaining parts of this theorem, 
namely linearity, monotonicity, and the absolute value property.
(You may use Theorem (\ref{thm: Lebesgue})).
\end{enumerate}
\end{exer}
\vfill\eject

\chapter{The Hilbert Space $\xL2$}
\section{Square Integrable Functions}

In this chapter we will develop the beginnings of a theory of 
function spaces with many properties analogous to the basic
properties of $\R^n$.  To motivate these developments we first
take a look at $\R^n$ in a different way.  We let 
$X$ be a finite set with $n$ elements, say, 
$X = \{1,2,3,\dots,n\}$ and we define a measure $\nu$ on $X$
which is called the ``counting measure''.  

More precisely, we take as $\sigma$-algebra the family of all
subsets of $X$ and for any $A \subset X$ we define $\nu(A)$
to be the number of elements in the set $A$.  It is easy to
see that this is a measure and that {\em any} function 
$f: X \to \R$ is measurable.  In fact any function is a simple
function.  This is because there is a partition of $X$ given
by $A_i = \{i\}$ and clearly $f$ is constant on each $A_i$,
so $f = \sum_{i=1}^n r_i \X_{A_i}$ where $r_i = f(i).$

Consequently we have 
\[
\int f \ d\nu = \sum_{i=1}^n r_i \nu(A_i) = \sum_{i=1}^n f(i).
\]

For reasons that will be clear below we will denote the collection
of all functions from $X$ to $\R$ by $L^2(X)$.  The important thing
to note is that this is just another name for $\R^n$.  More formally,
there is a vector space isomorphism of $L^2(X)$ and $\R^n$ given by
$f \longleftrightarrow (x_1,x_2,\dots, x_n)$ where $x_i = f(i).$

Under this isomorphism it is important to note what the inner product
(or ``dot'' product $<x,y> = \sum_{i=1}^n x_i y_i$ becomes.  If $f,g \in L^2(X)$
are the functions corresponding to vectors $x$ and $y$ respectively,
then $x_i = f(i)$ and $y_i = g(i)$ so
\[
<x,y> = \sum_{i=1}^n x_i y_i = \sum_{i=1}^n f(i) g(i)  = \int fg\ d\nu.
\]
Also the norm (or length) of a vector is given by
\[
\| x \|^2 = <x,x> = \sum_{i=1}^n x_i^2 = \sum_{i=1}^n f(i)^2 = \int f^2\ d\nu.
\]
It is this way of viewing the inner product and norm on $\R^n$ which 
generalizes nicely to a space of real valued functions on the interval.

In this chapter it will be convenient (for notational purposes) to
consider functions on the interval $[-1, 1]$ rather than $[0,1]$.
Of course, all of our results about measurable functions and their
integrals remain valid on this different interval.

\begin{defn} A measurable function $f :[-1,1] \to \R$ is called
\myindex{square integrable} if $f(x)^2$ is integrable.  We denote
the set of all square integrable functions by $\xL2.$  We define
the \myindex{norm} of $f \in \xL2$ by
\[
\|f\| = \Big(\int f^2\ d\mu \Big)^{\frac{1}{2}}.
\]
\end{defn}

\begin{prop}
The norm $\|\ \|$ on $\xL2$ satisfies $\|a f\| = |a| \|f\|$ for
all $a \in \R$ and all $f \in \xL2.$  Moreover for all $f, \
\|f\| \ge 0 $ with equality only if $f = 0$ almost everywhere.
\end{prop}
\begin{proof}
We see
\[
\|a f\| = \Big( \int a^2 f^2 \ d\mu\Big)^\frac{1}{2}
= \sqrt{a^2}\Big( \int  f^2 \ d\mu\Big)^\frac{1}{2}
= |a| \|f\|.
\]
Since, $\int f^2\ d\mu \ge  0$ clearly $\|f\| \ge 0.$  
Also if 
\[
\|f\| = 0 \text{, then } \int f^2\ d\mu = 0.
\]
So by Corollary (\ref{cor: zero ae})
$f^2 = 0$ almost everywhere and hence $f = 0$ almost everywhere.
\end{proof}

\begin{lemma}\label{lem: ineq}
If $f,g \in \xL2$, then $fg$ is integrable and 
\[
2\int |fg| \ d\mu \le \|f\|^2 + \|g\|^2.
\]
Equality holds if and only if $|f| = |g|$ almost everywhere.
\end{lemma}
\begin{proof}
Since 
\[
0 \le (|f(x)| - |g(x)|)^2  = f(x)^2 - 2|f(x)g(x)| + g(x)^2
\]
we have $2|f(x)g(x)| \le f(x)^2 + g(x)^2$.  Hence by Proposition
(\ref{prop: integral monotonicity}) we conclude that
$|fg|$ is integrable and that 
\[
2\int |fg| \ d\mu \le \|f\|^2 + \|g\|^2.
\]
Equality holds if and only if $\int (|f(x)| - |g(x)|)^2\ d\mu = 0$ and
we may conclude by Corollary (\ref{cor: zero ae}) that this happens
if and only if $(|f(x)| - |g(x)|)^2 = 0$ almost everywhere and hence that
$|f| = |g|$ almost everywhere.
\end{proof}

\begin{thm}
$\xL2$ is a vector space.
\end{thm}
\begin{proof}
We must show that if $f,g \in \xL2$ and $c \in \R$, then $cf \in \xL2$
and $(f+g) \in \xL2.$  The first of these is clear since $f^2$ integrable
implies that $c^2f^2$ is integrable.

To check the second we observe that 
\[
(f+g)^2 = f^2 +2fg + g^2 \le f^2 +2|fg| + g^2.
\]
Since $f^2,\ g^2$ and $|fg|$ are all integrable, it follows from
Proposition (\ref{prop: integral monotonicity}) that $(f+g)^2$ is
also.  Hence $(f+g) \in \xL2.$
\end{proof}

\begin{thm}[H\"older Inequality]\label{thm: Holder}
\index{H\"older Inequality}
If $f,g \in \xL2$, then
\[
\int |fg| \ d\mu \le \|f\|\ \|g\|.
\]
Equality holds if and  only if there is a constant $c$ such that 
$|f(x)| = c|g(x)|$ or $|g(x)| = c|f(x)|$ almost everywhere.
\end{thm}

\begin{proof}
If either $\|f\|$ or $\|g\|$ is $0$ the result is trivial so assume
they are both non-zero.  In that case
the functions $f_0 = f/\|f\|$ and $g_0 = g/\|g\|$ satisfy 
$\|f_0\| = \|g_0\| = 1.$

Then by Lemma~(\ref{lem: ineq}) 
\[
2\int |f_0g_0| \ d\mu \le \|f_0\|^2 + \|g_0\|^2 = 2,
\]
so 
\[
\int |f_0g_0| \ d\mu \le 1 ,
\]
and equality holds if and only if $|f_0| = |g_0|$ almost everywhere.
So 
\[
\frac{1}{\|f\|\ \|g\|}\int |fg| \ d\mu = \int |f_0g_0| \ d\mu \le 1
\]
and hence
\[
\int |fg| \ d\mu \le \|f\|\ \|g\|.
\]

Equality holds if and only if $|f_0| = |g_0|$ almost everywhere, which
implies there is a constant $c$ with $|f(x)| = c|g(x)|$ almost
everywhere.
\end{proof}

\begin{cor}\label{cor: angles} 
If $f,g \in \xL2$, then
\[
\Big| \int fg \ d\mu \Big | \le \|f\|\ \|g\|.
\]
Equality holds if and only if there is a constant $c$ such that
$f(x) = cg(x)$ or $g(x) = cf(x)$  almost everywhere.
\end{cor}
\begin{proof}
The inequality follows from H\"older's inequality and the absolute
value inequality since
\[
\Big| \int fg \ d\mu \Big | \le \int |fg| \ d\mu \le \|f\|\ \|g\|.
\]
Equality holds when both of these inequalities are equalities.
If this case, suppose first that $\int fg \ d\mu \ge 0$.  Then 
$\int |fg| \ d\mu =\int fg \ d\mu$, so
$\int |fg| -fg \ d\mu = 0$ and hence $|fg| = fg$ almost everywhere.
This says that $f$ and $g$ have the same sign almost everywhere.
Since the second inequality is an equality we know from H\"older
that there is a constant $c$ such that 
$|f(x)| = c|g(x)|$ or $|g(x)| = c|f(x)|$ almost everywhere.
This togther with the fact that $f$ and $g$ have the same sign 
almost everywhere implies $f(x) = cg(x)$ or $g(x) = cf(x)$ almost everywhere.
For the case that that $\int fg \ d\mu \le 0$ we can replace $f$ with $-f$
and conclude that $f(x) = -cg(x)$ or $g(x) = -cf(x)$.
Conversely, it is easy to see that if $f(x) = cg(x)$ or $g(x) = cf(x)$
almost everywhere, then the inequality above is an equality.

\end{proof}

The following result called the Minkowski Inequality, is the
triangle inequality for the vector space $\xL2$.

\begin{thm}[Minkowski's Inequality]\label{thm: triangle ineq}
\index{Minkowski Inequality}
If $f,g \in \xL2$, then
\[
\| f + g\| \le \|f\| + \|g\|.
\].
\end{thm}
\begin{proof}
We observe that
\begin{align*}
\| f + g\|^2 &= \int (f+g)^2\ d\mu\\
&= \int (f^2 +2fg + g^2)\ d\mu\\
&\le \int f^2 +2|fg| + g^2\ d\mu\\
&\le \|f\|^2 + 2\|f\|\ \|g\|+ \|g\|^2 \text{\ \ \ \ by H\"older's inequality}\\
&= (\|f\| + \|g\|)^2.
\end{align*}
Taking square roots of both sides of this equality gives the 
triangle inequality.
\end{proof}

\begin{defn}[Inner Product on $\xL2$] \index{inner product on $\xL2$}
If $f,g \in \xL2$, then we define
their \myindex{inner product} by
\[
\la f,g \ra = \int fg\ d\mu.
\]
\end{defn}

\begin{thm}[The Inner Product on $\xL2$]\label{thm: inner prod}
For any $f_1,f_2,g \in \xL2$ and any $c_1,c_2\in \R$ the inner product 
on $\xL2$ satisfies the following properties:
\begin{enumerate}
\item \myindex{Commutativity}: $\la f,g \ra = \la g,f \ra.$
\item \myindex{Bi-linearity}: $\la c_1f_1 +c_2 f_2,g \ra = c_1 \la f_1,g \ra + 
c_2 \la f_2,g \ra.$
\item \myindex{Positive Definiteness}: $\la g,g \ra = \|g\|^2 \ge 0$ 
with equality if and only if $g = 0$ almost everywhere.
\end{enumerate}
\end{thm}
\begin{proof}
Clearly $\la f,g \ra = \int fg\ d\mu = \int gf\ d\mu = \la g,f \ra.$
Bi-linearity holds because of the linearity of the integral. 
Also $\la g,g \ra = \int g^2\ d\mu \ge 0.$ 
Corollary (\ref{cor: zero ae})  implies that equality holds only if
$g^2 = 0$ almost everywhere.
\end{proof}

Notice that we have {\em almost} proved that $\xL2$ is an inner product
space.  The one point where the
definition is not quite satisfied is that $\|f\| = 0$ implies $f = 0$
{\em almost everywhere} rather than everywhere.  The pedantic way to
overcome this problem is to define $\xL2$ as the vector space of
equivalence classes of square integrable functions, where $f$ and $g$
are considered ``equivalent'' if they are equal almost everywhere.  It
is customary, however, to overlook this infelicity and simply consider
$\xL2$ as a vector space of functions rather than equivalence classes of
functions.  In doing this we should keep in mind that we are generally
considering two functions the same if they agree almost everywhere.

\section{Convergence in $\xL2$}

We have discussed uniform convergence and pointwise convergence and now
we wish to discuss convergence in the $\xL2$ norm $\|\ \|.$  The vector
space $\xL2$ is, of course, a metric space with distance function given
by $\dist( f, g) = \| f - g\|.$  Note that $\dist( f, g) = 0$ if and only
if $f = g$ almost everywhere, so again if we wish to be pedantic this
metric space is really the equivalence classes of functions which are
equal almost everywhere.

\begin{defn}
If $\{f_n\}_{n=1}^\infty$ is a sequence in $\xL2$, then it is said
to converge to in measure of order $2$ or to converge in $\xL2$ if
there is a function $f \in \xL2$ such that
\[
\lim_{n \to \infty} \|f - f_n\| = 0.
\]
\end{defn}

\begin{lemma}[Density of Bounded Functions]\label{prop: bounded fcn density}
If we define 
\[
f_n (x) =
\begin{cases}
n, \text{ if } f(x) > n;\\
f(x), \text{ if } -n \le f(x) \le n;\\
-n \text{ if } f(x) < -n,
\end{cases}
\]
then 
\[
\lim_{n \to \infty} \|f - f_n\| = 0.
\]
\end{lemma}
\begin{proof}
We will show that for any $\eps > 0$ there is an $n$ such
that $\|f - f_n\|^2 < \eps.$  First we note that $|f_n(x)| \le |f(x)|$ so
\[
|f(x) - f_n(x)|^2 \le |f(x)|^2 + 2 |f(x)|\ |f_n(x)| + |f(x)|^2 \le
4|f(x)|^2.
\]

Let $E_n = \{ x\ |\ |f(x)| > n \} = \{ x\ |\ |f(x)|^2 > n^2 \}$
and let $C = \int |f|^2 \ d\mu$.
Then
\[
C = \int |f|^2 \ d\mu \ge \int_{E_n} |f|^2 \ d\mu  
\ge \int_{E_n} n^2 \ d\mu = n^2\mu(E_n)
\]
and we conclude that $\mu(E_n) \le C/n^2.$

We know from absolute continuity, Theorem (\ref{thm: abs continuity}),
that there is a $\delta >0$ such that 
$\int_A |f|^2 \ d\mu < \eps/4$ whenever $\mu(A) < \delta.$  Thus we have
\[
\|f - f_n\|^2 = \int |f - f_n|^2 \ d\mu 
= \int_{E_n} |f - f_n|^2 \ d\mu 
\le \int_{E_n} 4|f|^2 \ d\mu < 4 \frac{\eps}{4} = \eps
\]
whenever $n$ is sufficiently large that $\mu(E_n) \le C/n^2 < \delta.$
\end{proof}

\begin{prop}[Density of Step Functions and Continuous Functions]
\label{prop: fcn density} \index{density, of step and
continuous functions in $L^2$} The step functions are dense in $\xL2$.
That is, for any $\eps >0$ and any $f \in \xL2$ there is a step
function $g: [-1,1] \to \R$ such that $\|f - g\| < \eps.$ Likewise,
there is a continuous function $h: [-1,1] \to \R$ such that $\|f - h\|
< \eps.$  The function $h$ may be chosen so $h(-1) = h(1).$
\end{prop}

\begin{proof}
By the preceding result we may choose $n$ so that $\|f - f_n\| <\frac{\eps}{2}.$
Note that $|f_n(x)| \le n$ for all $x.$ 
Suppose now that $\delta$ is any given small positive number.
According to Theorem (\ref{thm: step approx})
there is a step function $g$ with $|g| \le n$  and a measurable 
set $A$ with $\mu(A) < \delta$
such that $|f_n(x) - g(x)| < \delta$ if $x \notin A.$  Hence
\begin{align*}
\|f_n - g\|^2 &= \int |f_n-g|^2 \ d\mu \\
&= \int_A |f_n - g|^2 \ d\mu + \int_{A^c} |f_n-g|^2 \ d\mu \\
&\le \int_A 4n^2 \ d\mu + \int_{A^c} \delta^2 \ d\mu \\
&\le  4n^2\mu(A) + \delta^2 \mu(A^c) \le 4n^2 \delta + 2\delta^2.
\end{align*}

Clearly if we choose $\delta$ sufficiently small, then
\[
\|f_n - g\| \le \sqrt{4n^2 \delta + 2\delta^2} < \frac{\eps}{2}.
\]
It follows that $\|f - g\| \le \|f - f_n\|  + \|f_n - g\| < \eps.$

The proof for continuous functions is the same, except 
Exercise~(\ref{exer: f=g=h ae}) is used in place of 
Theorem~(\ref{thm: step approx}).  The details are left as an
exercise.
\end{proof}

\begin{defn}
An inner product space 
$(\V, \la\ ,\ \ra)$ which is \myindex{complete}, i.e. in which Cauchy sequences
converge, is called a \myindex{Hilbert space}.
\end{defn}

For example, $\R^n$ with the usual dot product is a Hilbert space.

We want to prove that $\xL2$ is a Hilbert space. 

\begin{thm}\label{thm: L2 Hilbert}
$\xL2$ is a Hilbert space.
\end{thm}
\begin{proof}
We have already shown that $\xL2$ is an inner product space.  All that
remains is to prove that the norm $\|\ \|$ is complete, i.e. that 
Cauchy sequences converge.

Let $\{f_n\}_{n=1}^\infty$ be a Cauchy sequence.  Then we may choose numbers
$n_i$ such that $\| f_m - f_n\| < 1/2^i$ whenever $m,n \ge n_i.$
Hence if we define $g_0 = 0$ and $g_i = f_{n_i}$ for $i>0$, then
$\| g_{i+1} - g_i\| < 1/2^i$ so, in particular 
$\sum_{i=0}^\infty \| g_{i+1} - g_i\|$ converges, say to $S.$

Consider the function $h_n(x)$ defined by
\[
h_n(x) = \sum_{i=0}^{n-1} | g_{i+1}(x) - g_i(x)|.
\]
For any fixed $x$ the sequence $\{h_n(x)\}$ is monotone increasing
so we may define the extended real valued function $h$ by
$\displaystyle{h(x) = \lim_{n \to \infty} h_n(x).}$   
Note that by the Minkowski inequality
\[
\|h_n\| \le \sum_{i=0}^{n-1} \| g_{i+1} - g_i\| < S.
\]

Hence $\int h_n^2\ d\mu = \|h_n\|^2 < S^2.$  Since $h_n(x)^2$ is a monotonic
increasing sequence of non-negative measurable functions converging to
$h^2$ we conclude from the Monotone Convergence Theorem 
(\ref{thm: monotone convergence}) that 
$\displaystyle{\int h^2 \ d\mu = \lim_{n \to \infty} \int h_n^2\ d\mu < S^2}$
so $h^2$ is integrable.

Since $h^2$ is integrable, $h(x)$ is finite almost everywhere.  For each
$x$ with finite $h(x)$ the series of real numbers
$\sum_{i=0}^\infty (g_{i+1}(x) - g_i(x))$ converges absolutely and
hence converges by Theorem (\ref{thm: absolute convergence}).  We denote
its sum by $g(x).$  For
$x$ in the set of measure $0$ where $h(x) = +\infty$ we define $g(x) = 0.$
Notice that 
\[
g_n(x) = \sum_{i=0}^{n-1} (g_{i+1}(x) - g_i(x)) 
\]
because it is a telescoping series.   Hence 
\[
\lim_{n \to \infty} g_n(x) 
= \lim_{n \to \infty} \sum_{i=0}^{n-1} (g_{i+1}(x) - g_i(x)) 
= g(x)
\]
for almost all $x$.  
Moreover 
\[
|g(x)| = \lim_{n \to \infty} |g_n(x)| 
\le \lim_{n \to \infty} \sum_{i=0}^{n-1} |g_{i+1}(x) - g_i(x)| 
= \lim_{n \to \infty} h_n(x) = h(x)
\]
for almost all $x$ so $|g(x)|^2 \le h(x)^2$ and hence $|g(x)|^2$ is
integrable and $g \in \xL2.$

We also observe that 
\[
|g(x) - g_n(x)|^2 \le
(|g(x)| + |g_n(x)|)^2 \le (2 h(x))^2.
\]
Since $\displaystyle{\lim_{n\to \infty} |g(x) - g_n(x)|^2 = 0}$ for
almost all $x$ the Lebesgue Convergence Theorem 
(\ref{thm: general lebesgue convergence})
tells us
$\displaystyle{\lim_{n\to \infty} \int |g(x) - g_n(x)|^2\ d\mu = 0}$.
This implies
$\displaystyle{\lim_{n\to \infty} \|g - g_n\| = 0}$.

Hence given $\eps >0$ there is an $i$ such that $\|g - g_i\| < \eps/2$
and $1/2^i < \eps/2.$  Recalling that $g_i = f_{n_i}$ we see that
whenever  $m \ge n_i$ we have
$\|g - f_m\| \le \|g - g_i\| + \|g_i - f_m\| < \eps/2 + \eps/2 = \eps.$
Hence $\displaystyle{\lim_{m\to \infty} \|g - f_m\| = 0}$.
\end{proof}

\section{Hilbert Space}

In any Hilbert space we can, of course, talk about convergent sequences
and series.  The meaning is precisely what you would expect.  In particular,
if $\HH$ is a Hilbert space and $\{x_n\}$ is a sequence, then
\[
\lim_{n \to \infty} x_n = x
\]
means that for any $\eps >0$ there is an $N>0$ such that $\|x - x_n\| <\eps$
whenever $n \ge N.$  This is exactly the usual definition in $\R$ except
we use the norm $\|\ \|$ in place of absolute value.
Also if $\{u_n\}$ is a sequence in $\HH$, then 
\[
\sum_{m = 1}^\infty u_m = s
\]
means $\lim s_n = s$ where
\[
s_n = \sum_{m = 1}^n u_m.
\]
We will say a series $\sum_{m = 1}^\infty u_m$ 
{\em converges absolutely} \index{absolute convergence}
provided $\sum_{m = 1}^\infty \|u_m\|$ converges.

\begin{prop} If a series in a Hilbert space converges absolutely
then it converges.
\end{prop}
\begin{proof}
Given $\eps > 0$ there is an $N >0$ such that
whenever $n> m \ge N,$
\[
\sum_{i = m}^n \|u_m\| \le \sum_{i = m}^\infty \|u_m\| <\eps.
\]
Let $s_n = \sum_{i = 1}^n u_i$, then 
$\|s_n -s_m \| \le \sum_{i = m}^n \|u_m\| < \eps.$
It follows that $\{s_n\}$ is a Cauchy sequence. Hence it
converges.
\end{proof}

We will also talk about perpendicularity in $\HH$.  We say 
$x,y \in \HH$ are \myindex{perpendicular} (written $x \perp y$) if
$\la x,y\ra = 0.$

\begin{thm}[Pythagorean Theorem]\label{thm: pythag}
\index{Pythagorean Theorem}
If $x_1, x_2, \dots x_n$ are mutually perpendicular elements of
a Hilbert space, then
\[
\Big \| \sum_{i=1}^n x_i\Big \|^2 = \sum_{i=1}^n \|x_i\|^2.
\]
\end{thm}
\begin{proof}
Consider the case $n=2$.  If $x \perp y$, then
\[
\| x + y\|^2 = \la x+ y, x+y \ra = \la x, x \ra + 2\la x, y \ra  
+ \la y, y \ra =  \|x\|^2 + \|y\|^2
\]
since $\la x,y \ra = 0.$
The general case follows by induction on $n$.
\end{proof}

\begin{defn} If $\HH$ is a Hilbert space, a {\em bounded linear functional}
\index{linear functional}
on $\HH$ is a function $L: \HH \to \R$ such that for all $v,w \in \HH$ and 
$c_1,c_2 \in \R,\ L(c_1v +c_2w) = c_1L(u) + c_2L(w)$ and such that there is
a constant $M$ satisfying $|L(v)| \le M\|v\|$ for all $v \in \HH.$
\end{defn}

The following result was proved in Proposition~(\ref{prop: cauchy schwartz}).
In the case of the Hilbert space $\xL2$ it is just the corollary
to H\"older's inequality, 
Corollary~(\ref{cor: angles}).

\begin{prop}[Cauchy-Schwarz Inequality]\label{prop: Hilbert cauchy schwartz}
\index{Cauchy-Schwarz Inequality}
If $(\HH, \la\ ,\ \ra)$ is a Hilbert space and $v,w \in \HH$, then 
\[
|\la v, w \ra| \le \|v\|\ \|w\|,
\]
with equality if and only if $v$ and $w$ are multiples of a single
vector.
\end{prop}

For any fixed $x \in \HH$
we may define $L: \HH \to \R$ by $L(v) = \la v,x\ra$.  Then $L$
is a linear function and as a consequence of the Cauchy-Schwarz
inequality it is bounded.  Indeed $\|L(v)\| \le M \|v\|$ where
$M = \|x\|$.   Our next goal is to prove that these are the 
{\em only} bounded linear functions from $\HH$ to $\R.$

\begin{lemma}\label{lem: convex}
Suppose $\HH$ is a Hilbert space and $L: \HH \to \R$ is a bounded
linear functional which is not identically $0$. If
$\V = L^{-1}(1)$.
Then there is a unique $x \in \V$ such that
\[
\|x\| = \inf_{v \in \V} \|v\|.
\]
That is, there is a unique vector in $\V$ closest to $0$. 
Moreover, the vector $x$
is perpendicular to every element of $L^{-1}(0)$, i.e. if $v \in \HH$
and $L(v) = 0$, then $\la x,v \ra = 0.$
\end{lemma}
\begin{proof}
We first observe that $\V$ is closed, i.e. that any convergent sequence
in $\V$ has a limit in $\V$.  To see this suppose $\lim x_n = x$ and
$x_n \in \V.$  Then $|L(x) - L(x_n)| =  |L(x-x_n)| \le M \|x - x_n\|$
for some $M$.  Hence since $L(x_n) = 1$ for all $n$, we have
$|L(x) - 1|  \le \lim M \|x - x_n\| = 0.$  Therefore $L(x) = 1$ and
$x \in \V.$

Now let $\displaystyle{d = \inf_{v \in \V} \|v\|}$ and choose a sequence 
$\{x_n\}_{n=1}^\infty$ in $\V$ such that $\lim \|x_n\| = d.$
We will show that this sequence is Cauchy and hence converges.

Notice that $(x_n + x_m)/2$ is in $\V$ so $\|(x_n + x_m)/2\| \ge d$ 
or $\|x_n + x_m\| \ge 2d$.
By the parallelogram law (Proposition~(\ref{prop: norm props})
\[
\| x_n - x_m\|^2 + \| x_n + x_m\|^2 = 2\| x_n\|^2 + 2\| x_m\|^2.
\]
Hence
\[
\| x_n - x_m\|^2  = 2\| x_n\|^2 + 2\| x_m\|^2 - \| x_n + x_m\|^2
\le 2\| x_n\|^2 + 2\| x_m\|^2 - 4d^2.
\]
As $m$ and $n$ tend to infinity the right side of this equation goes
to $0$.  Hence the left side does also and $\lim \| x_n - x_m\| = 0.$
That is, the sequence $\{x_n\}_{n=1}^\infty$ is Cauchy.  Let $x \in
\V$ be limit limit of this sequence.  Since $\|x\| \le \|x - x_n\| +
\|x_n\|$ for all $n$, we have
\[
\|x\| \le \lim_{n \to \infty}\|x - x_n\| + \lim_{n \to \infty}\|x_n\|
= d.
\]
But $x \in \V$ implies $\|x\| \ge d$ so $\|x\| = d.$

To see that $x$ is unique suppose that $y$ is another element
of $\V$ and $\|y\| = d.$  Then $(x+y)/2$ is in $\V$
so $\|x + y\| \ge 2d.$ 
Hence using the parallelogram law again
\[
\| x - y\|^2  = 2\| x\|^2 + 2\| y\|^2 - \| x + y\|^2 \le 4d^2 - 4d^2 = 0.
\]
We conclude that $x = y.$

Suppose that $v \in L^{-1}(0)$.  We wish to show it is perpendicular
to $x$.  Note that for all $t \in \R$ the vector $x + tv \in L^{-1}(1)$ 
so $\|x + tv\|^2 \ge \|x\|^2.$  Hence
\[
\|x\|^2 +2t<x,v> +t^2 \|v\|^2 \ge \|x\|^2, \text{ so }
\]
$2t<x,v> +t^2 \|v\|^2 \ge 0$ for all $t \in \R.$  This is possible
only if $<x,v> = 0$.
\end{proof}

In the following theorem we characterize all the bounded linear
functionals on  a Hilbert space.  Each of them is obtained by
taking the inner product with some fixed vector.

\begin{thm}
If $\HH$ is a Hilbert space and $L: \HH \to \R$ is a bounded linear
functional, then there is a unique $x \in \HH$ such that
$L(v) = \la v,x\ra.$
\end{thm}
\begin{proof}
If $L(v) =0$ for all $v$, then $x = 0$ has the property we want,
so suppose $L$ is not identically $0.$
Let $x_0 \in \HH$ be the unique point in $L^{-1}(1)$ with smallest
norm, guaranteed by Lemma~(\ref{lem: convex}).

Suppose first that $v \in \HH$ and $L(v) = 1$  Then 
$L(v - x_0) = L(v) - L(x_0) = 1 - 1 = 0$  so by Lemma~(\ref{lem: convex})
$\la v-x_0, x_0\ra = 0.$  It follows that the vector $x = x_0/\|x_0\|^2$ is also
perpendicular to $v-x_0$ so
\[
\la v,x\ra =
\la v,\frac{x_0}{\|x_0\|^2}\ra  = \la v-x_0 ,\frac{x_0}{\|x_0\|^2}\ra  
+ \la x_0,\frac{x_0}{\|x_0\|^2}\ra  = 1 = L(v).
\]
Hence for any $v$ with $L(v) = 1$ we have $L(v) = \la v,x\ra.$
Also for any $v$ with $L(v) = 0$ we have $L(v) = 0 = \la v,x\ra$
by Lemma~(\ref{lem: convex}).

Finally for an arbitrary $w \in \HH$ with $L(w) = c \ne 0$ 
we define $v = w/c$ so $L(v) = L(w)/c = 1.$
Hence 
\[
L(w) = L(cv) = cL(v) = c \la v,x\ra
= \la c v, x\ra = \la w, x\ra.
\]

To see that $x$ is unique, suppose that $y \in \HH$ has the same properties
then for every $v \in \HH$ we have $\la v,x \ra = L(v) = \la v, y\ra.$
Thus $\la v,x-y\ra = 0$ for all $v$ and in particular for $v = x -y.$ We conclude
that $\|x-y\|^2 = \la x-y,x-y\ra = 0$ so $x = y$.

\end{proof}

\section{Fourier Series}

\begin{defn}
A family of vectors $\{u_n\}$ in a Hilbert space $\HH$ is
called {\em orthonormal}\index{orthonormal family} provided
for each $n,\  \|u_n\| = 1$ and $\la u_n, u_m \ra = 0$ if
$n \ne m.$
\end{defn}

\begin{thm}
The family of functions 
$\F = \{ \frac{1}{\sqrt{2}},  \cos(n \pi x), \sin(n\pi x)\}_{n=1}^\infty$
is an orthonormal family in $\xL2.$
\end{thm}

For a proof see Chapter 1 of \cite{Pin}.

\begin{thm}\label{thm: mean approx}
If $\{u_n\}_{n=0}^N$ is a finite orthonormal family of vectors in a
Hilbert space $\HH$ and $w \in \HH$, then the 
minimum value of 
\[
\Big \| w - \sum_{n=0}^N c_n u_n \Big \|
\]
for all choices of $c_n \in \R$ occurs when $c_n = \la w, u_n\ra.$
\end{thm}
\begin{proof}
Let $c_n$ be arbitrary real numbers and define $a_n = \la w, u_n \ra$ 
Let
\[
u = \sum_{n=0}^N a_n u_n, \text{ and } v = \sum_{n=0}^N c_n u_n.
\]
Notice that by Theorem~(\ref{thm: pythag}) 
$\la u, u\ra =  \sum_{n=0}^N a_n^2$ and 
$\la v, v\ra =  \sum_{n=0}^N c_n^2$.  Also 
\[
\la w, v\ra =  \sum_{n=0}^N c_n \la w,u_n\ra =\sum_{n=0}^N a_n c_n
\]
Hence
\begin{align*}
\| w - v\|^2 &= \la w - v, w - v\ra\\
&= \|w\|^2 - 2\la w, v\ra + \|v\|^2\\
&= \|w\|^2 - 2\sum_{n=0}^N a_n c_n +  \sum_{n=0}^N c_n^2\\
&= \|w\|^2 - \sum_{n=0}^N a_n^2 + \sum_{n=0}^N (a_n - c_n)^2\\
&= \|w\|^2  - \|u\|^2  + \sum_{n=0}^N (a_n - c_n)^2.
\end{align*}
It follows that
\[
\| w - v\|^2 \ge \|w\|^2  - \|u\|^2 
\]
for any choices of the $c_n$'s and we have equality if only if
$c_n = a_n = \la w, u_n\ra.$  That is, for all choices of
$v$, the minimum value of
$\| w - v\|^2$ occurs precisely when $v = u.$
\end{proof}

\begin{defn}
If $\{u_n\}_{n=0}^\infty$ is an orthonormal family of vectors in a
Hilbert space $\HH$, it is called {\em complete}
\index{complete orthonormal family} if 
every $w \in \HH$ can be written as an infinite series
\[
w = \sum_{n=0}^\infty c_n u_n
\]
for some choice of the numbers $c_n \in \R.$
\end{defn}

Theorem~(\ref{thm: mean approx}) suggests that the only reasonable
choice for $c_n$ is $c_n = \la w, u_n\ra$ and we will show
that this is the case.  These numbers are sufficiently frequently 
used that they have a name.

\begin{defn}[Fourier Series]
The \myindex{Fourier coefficients} of $w$ with respect 
to an orthonormal family $\{u_n\}_{n=0}^\infty$ are the
numbers $\la w, u_n\ra$.  The infinite series
\[
\sum_{n=0}^\infty \la w, u_n \ra u_n
\]
is called the \myindex{Fourier series}.
\end{defn}

\begin{ex}[Classical Fourier Series]
We will show later that the orthnormal family of functions 
$\F = \{ \frac{1}{\sqrt{2}},
\cos(n \pi x),  \sin(n \pi x),\}_{n=1}^\infty$
is complete.  If $f \in \xL2$, then the Fourier coefficients are
\begin{align*}
A_0 &= \frac{1}{\sqrt{2}}\int f\ d\mu\\
A_n &= \int f \cos(n \pi x)\ d\mu \text{ for $n > 0$}\\
B_n &= \int f \sin(n \pi x)\ d\mu \text{ for $n > 0$ },
\end{align*}
and the Fourier series is
\[
\frac{1}{\sqrt{2}}A_0 + \sum_{n=1}^\infty A_n \cos(n \pi x)
+ \sum_{n=1}^\infty B_n \sin(n \pi x)
\]
\end{ex}

\begin{thm}[Bessel's Inequality]\label{thm: bessel}
\index{Bessel's inequality}
If $\{u_i\}_{i=0}^\infty$ is an orthonormal family of vectors in a Hilbert space 
$\HH$ and $w \in \HH$, then the series
\[
\sum_{i=0}^\infty \la w, u_i \ra^2 \le \|w\|^2.
\]
In particular this series converges.
\end{thm}

\begin{proof}
Let $s_n$ be the partial sum for the Fourier series.  That is,
$s_n = \sum_{i=0}^n \la w, u_n \ra u_n.$
Then since the family is orthogonal, we know by Theorem~(\ref{thm: pythag})
that 
\begin{equation}\label{eqn: fourier}
\|s_n\|^2 = \sum_{i=0}^n \|\la w, u_i \ra u_i\|^2 
= \sum_{i=0}^n \la w, u_i \ra^2.
\end{equation}
This implies  that $s_n \perp (w -s_n)$ because
\[
\la w- s_n, s_n\ra  = \la w, s_n\ra - \la s_n, s_n \ra
= \sum_{i=0}^n \la w, u_n \ra^2 - \|s_n\|^2 = 0.
\]

Since $s_n \perp (w -s_n)$ we know 
\begin{equation}\label{eqn: fourier2}
\|w\|^2 = \|s_n\|^2 +  \|w -s_n\|^2
\end{equation}
by Theorem~(\ref{thm: pythag}) again.  Hence by equation~(\ref{eqn: fourier})
$\sum_{i=0}^n \la w, u_n \ra^2  = \|s_n\|^2 \le \|w\|^2.$
Since $\|s_n\|^2$ is an increasing sequence
it follows that the series
\[
\sum_{i=0}^\infty\la w, u_n \ra^2
= \lim_{n \to \infty}\|s_n\|^2 \le \|w\|^2
\]
converges.
\end{proof}

\begin{cor}
If $\{u_n\}_{n=0}^\infty$ is an orthonormal family of vectors in a
Hilbert space $\HH$ and $w \in \HH$, then the Fourier series 
$\sum_{i=0}^\infty \la w, u_i \ra u_i$
with respect to $\{u_i\}_{i=0}^\infty$ converges.
\end{cor}

\begin{proof}
Let $s_n$ be the partial sum for the Fourier series.  That is,
$s_n = \sum_{i=0}^n \la w, u_i \ra u_i.$ So if $n > m,\ 
s_n - s_m = \sum_{i=m+1}^n \la w, u_i \ra u_i.$

Then since the family is orthogonal, we know by Theorem~(\ref{thm: pythag})
that 
\[
\|s_n - s_m\|^2 = \sum_{i=m+1}^n \|\la w, u_i \ra u_i\|^2 
= \sum_{i=m+1}^n \la w, u_i \ra^2.
\]
Since the series $\sum_{i=0}^\infty \la w, u_i \ra^2$ converges
we conclude that given $\eps >0$ there is an $N>0$ such that
$\|s_n - s_m\|^2 < \eps^2$ whenever $n,m \ge N.$  In other words
the sequence $\{s_n\}$ is Cauchy.
\end{proof}

If Bessel's inequality is actually an equality, then
the Fourier series for $w$ must converge to $w$ in $\HH.$

\begin{thm}[Parseval's Theorem]
If $\{u_n\}_{n=0}^\infty$ is an orthonormal family of vectors in a
Hilbert space $\HH$ and $w \in \HH$, then
\[
\sum_{i=0}^\infty \la w, u_i \ra^2 = \|w\|^2
\]
if and only if the Fourier series with respect 
to $\{u_n\}_{n=0}^\infty$ converges
to $w$, i.e. 
\[
\sum_{i=0}^\infty \la w, u_i \ra u_i = w.
\]
\end{thm}

\begin{proof}
As above let $s_n$ be the partial sum for the Fourier series.
We showed in equation~(\ref{eqn: fourier2}) that
$\|w\|^2 = \|s_n\|^2 +  \|w -s_n\|^2$.  Clearly, then,
$\lim \|w -s_n\| = 0$ if and only if  $\lim \|s_n\|^2 = \|w\|^2.$ 
Equivalently (using equation (\ref{eqn: fourier})) 
$\sum_{n=0}^\infty \la w, u_n \ra u_n = w$
if and only if
$\sum_{n=0}^\infty \la w, u_n \ra^2 = \|w\|^2.$
\end{proof}

Recall that an \myindex{algebra of functions} is a vector space
$\A$ of real valued functions with the additional property that 
if $f,g \in \A$, then $fg \in \A.$  If $X = [a,b]$ is a 
closed interval in $\R$
we will denote by $C(X)$ the algebra of all continuous functions
from $X$ to $\R$  and by $C_{end}(X) = \{f\ |\ f(a) = f(b)\}$,
the subalgebra of functions which agree at the endpoints.
The following theorem is a special case of
a much more general theorem called the Stone-Weierstrass Theorem.

\begin{thm}
Suppose that $X = [-1,1]$ and $\A \subset C_{end}(X)$ 
is an algebra satisfying 
\begin{enumerate}
\item The constant function $1$ is in $\A,$ and
\item $\A$ separates points (except endpoints): for any distinct 
$x,y \in X$ with $\{x,y\} \ne \{-1,1\}$ there is $f \in \A$
such that $f(x) \ne f(y).$
\end{enumerate}
Then $\A$ is dense in $C(X)$, i.e. given any $\eps >0$ and any $g \in C(X)$
there is $f \in \A$ such that $|f(x) - g(x)| < \eps$ for all $x \in X.$
\end{thm}

A proof can be found in 5.8.2 of \cite{Mar} or in \cite{R}.
This result is usually stated in greater generality than we
do here.  For example the set $X$ need only be a compact
metric space, but since we have not defined these concepts
we state only the special case above.

\begin{cor}
If $\eps >0$ and $g: [-1,1] \to \R$ is a continuous function
satisfying $g(-1) = g(1)$, then there are $a_n, b_n \in \R$
such that $|g(x) -p(x)| < \eps,$ for all $x$, where
\[
p(x) = a_0 + \sum_{n=1}^N a_n \cos(n \pi x)
+ \sum_{n=1}^N b_n \sin(n \pi x).
\]
\end{cor}
\begin{proof}
Let $X$ be the unit circle in the plane $\R^2$, i.e. 
$X = \{ (\cos(\pi x), \sin( \pi x))\ |\ x \in [-1,1] \}.$
So if $\phi: [-1,1] \to \R^2$ is given by 
$\phi(x) = (\cos(\pi x), \sin( \pi x))$, then $X = \phi([-1,1]).$
For any function $f: [-1,1] \to \R$, with $f(-1) = f(1)$ we define
$\hat f : X \to \R$ to be the continuous function such that
$\hat f(\phi(x)) = f(x).$  We need the fact that
$f(-1) = f(1)$ because $\phi(-1) = \phi(1).$  Conversely given
any function $\hat h \in C(X)$ we can define $h: [-1,1] \to \R$
by $h(x) = \hat h(\phi(x))$ and we will have $h(-1) = h(1).$

Let $\A$ be the collection of all functions on $[-1,1]$
of the form
\[
q(x) = a_0 + \sum_{n=1}^N a_n \cos(n \pi x)
+ \sum_{n=1}^N b_n \sin(n \pi x).
\]
for some choices of $N, a_n,$ and $b_n.$  Then $\A$ is a vector
space and contains the constant function $1$.  
It is an algebra as a consequence of the trigonometric
identities
\begin{align*}
\sin(x) \cos(y) &= \frac{1}{2} \big (\sin(x+y) + \sin(x-y) \big)\\
\cos(x) \cos(y) &= \frac{1}{2} \big (\cos(x+y) + \cos(x-y) \big)\\
\sin(x) \sin(y) &= \frac{1}{2} \big (\cos(x+y) - \cos(x-y) \big)
\end{align*}

It is also the case that $\A$ separates points with the exception
of the one pair of points $x = -1,\ y= 1.$  To see this note that
if $x$ and $y$ are not this pair and if one is positive and one
negative, then $\sin(\pi x) \ne \sin(\pi y)$.  On the other hand if both are
$\ge 0$ or both $\le 0$, then $\cos(\pi x) \ne \cos(\pi y)$.

It follows that if $\hat \A = \{\hat q\ |\  q \in \A\}$, then
$\hat \A$ is an algebra which separates points of $X$ and contains the
constaint function $1$.  Note that the points  $x = -1,\ y= 1$ correspond to a
single point of $X$, namely $(-1,0) = \phi(-1) = \phi(1).$  So they cause
no problem.  Thus $\hat \A$ satisfies the hypothesis of the Stone-Weierstrass
theorem.

Hence given $\eps >0$ and $g: [-1,1] \to \R$, a continuous function
satisfying $g(-1) = g(1),$ we consider $\hat g$.  By Stone-Weierstrass
there is a $\hat p \in \hat \A$ such that the value of $\hat p$ differs
from the value of $\hat g$ by less than $\eps$ for all points of $X$.
Thus if $p(x) = \hat p(\phi(x))$ we have $|g(x) - p(x)| = 
|\hat g(\phi(x)) - \hat p(\phi(x))| < \eps$ and $p \in \A.$
\end{proof}

\begin{thm} If $f \in \xL2$, then the Fourier series for $f$
with respect to the orthonormal family $\F$ converges to $f$ in $\xL2$.
In particular the orthonormal family $\F$ is complete.
\end{thm}
\begin{proof}
Given $\eps >0$, we know
by Proposition~(\ref{prop: fcn density})  there is a continuous function
$g \in \xL2$ such that $g(-1) = g(1)$ and $\|f - g\| < \eps/2.$

By the corollary to the Stone-Weierstrass theorem there is a function 
\[
p(x) = a_0 + \sum_{n=1}^N a_n \cos(n \pi x)
+ \sum_{n=1}^N b_n \sin(n \pi x).
\]
with $|g(x) -p(x)| < \eps/4$ for all $x$. So
\[
\|g - p\|^2 = \frac{1}{\pi}\int (g - p)^2\ d\mu 
\le \frac{1}{\pi}\int \frac{\eps^2}{16}\ d\mu = \frac{\eps^2}{8}.
\]
Hence, 
$\|f - p\| \le \|f - g\| + \|g - p\| <  \eps/2 + \eps/\sqrt{8} < \eps.$

Let
\[
S_N(x) = \frac{1}{\sqrt{2}}A_0 + \sum_{n=1}^N A_n \cos(n \pi x)
+ \sum_{n=1}^N B_n \sin(n \pi x)
\]
where $A_n$ and $B_n$ are the Fourier coefficients for $f$ with respect
to $\F.$  Then $S_N(x)$ is the partial sum of the Fourier series of 
$f$.  According to Theorem~(\ref{thm: mean approx}) 
for every $m \ge N,\ \| f - S_m\| \le \| f - p\|$ so $\| f - S_m\| < \eps.$
This proves $\displaystyle{\lim \| f - S_m\| = 0.}$
\end{proof}

\begin{exer} 
Suppose $X = [-1,1].$
\begin{enumerate}
\item Prove that $C_{end}(X)$ is a subalgebra of $C(X)$, i.e.
it is a vector subspace closed under multiplication.
\item Let $\A_p$ be the polynomials 
\end{enumerate}
\end{exer}

$\setcounter{section}{0}
$\renewcommand{\thechapter}{\Alph{chapter}}
\appendix
\chapter[Lebesgue Measure]{Lebesgue Measure}
\section{Introduction}
We want to define a generalization of length called {\em
measure} for bounded subsets of the real line or subsets of the
interval $[a,b]$.  
There are several properties which we want it to have.
For each bounded subset $A$ of $\R$ we would like to be able
to assign a non-negative
real number $\mu(A)$ that satisfies the following:
\begin{description}
\item[I. Length.] If $A = (a,b)$ or $[a,b],$ 
then $\mu(A) = \len (A) = b-a,$ i.e., the measure of an open or closed interval
is its length
\item [II. Translation Invariance.]
If $A \subset \R$ is a bounded subset of $\R$ and $c\in \R$, then 
$\mu( A + c ) = \mu(A),$ where $A + c$ denotes the set
$\{ x +c\ |\ x \in A\}.$

\item [III. Countable Additivity.]
If $\{A_n\}_{n=1}^\infty$ is a countable collection of
bounded subsets of $\R$, then
\[
\mu( \bigcup_{n=1}^\infty A_n) \le \sum_{n=1}^\infty \mu(A_n)
\]
and if the sets are  {\em pairwise disjoint}, then
\[
\mu( \bigcup_{n=1}^\infty A_n) = \sum_{n=1}^\infty \mu(A_n)
\]
Note the same conclusion applies to finite collections
$\{A_n\}_{n=1}^m$ of bounded sets (just let $A_i = \emptyset$  for $i > m$).

\item [IV. Monotonicity] If $A \subset B$, then $\mu(A) \le \mu(B).$
Actually, this property is a consequence of additivity since $A$ and
$B \setminus A$ are disjoint and their union is $B$.
\end{description}
\bigskip

It turns out that it is not possible to find a $\mu$ which
satisfies I--IV  and which is defined for {\em all}
bounded  subsets of the reals.  But we can do it for a very large
collection including the open sets and the closed sets.

\section{Outer Measure}\label{sec: outer}

We first describe the notion of ``outer measure'' which comes close to
what we want.  It is defined for {\em all} bounded sets of the reals
and satisfies properties I and II above.  It also satisfies the
inequality part of the additivity condition, III, which is called
\myindex{subadditivity.}  But it fails to be additive for some choices of
disjoint sets.  The resolution of this difficulty will be to restrict
its definition to a certain large collection of nice sets (called
measurable) on which the additivity condition holds.  Our task is to
develop the definition of \myindex{measurable set}, to define the notion
of \myindex{Lebesgue measure} for such a set and, then to prove that 
properties I-IV hold, if we restrict our attention to measurable
sets.

Suppose $A \subset \R$ is a bounded set and $\{ U_n\}$ is a countable
covering of $A$ by open intervals, i.e. $A \subset \bigcup_n U_n$ where
$U_n = (a_n,b_n)$.  Then if we were able to define a function 
$\mu$ satisfying the properties I-IV above we would expect that
\[
\mu(A) \le \mu\big ( \bigcup_{n=1}^\infty U_n \big)
\le \sum_{n=1}^\infty \mu( U_n) = \sum_{n=1}^\infty \len( U_n)
\]
and hence that $\mu(A)$ is less than or equal to the {\it infimum}
of all such sums where we consider all possible coverings of $A$
by a countable collection of open intervals.   This turns out to
be a very useful definition.

\begin{defn}[Lebesgue Outer Measure]\label{def: outer}
\index{outer measure}
Suppose $A \subset \R$ is a bounded set and $\U(A)$ is the 
collection of all countable coverings of $A$ by open intervals.
We define the \myindex{Lebesgue outer measure} $\ms(A)$ by
\[
\ms(A) = \inf_{\{U_n\} \in \U(A)} \Big\{\sum_{n=1}^\infty \len ( U_n)\Big\},
\]
where the {\it infimum} is taken over all possible countable coverings
of $A$ by open intervals.
\end{defn}

Notice that this definition together with the definition of a null 
set, Definition~(\ref{def: null}),
says that a set $A \subset I$ is a null set if and only if
$\ms(A) = 0.$

We can immediately show that property I, the length property, holds
for Lebesgue outer measure.

\begin{prop}\label{prop: outer length}
For any $a,b \in \R$ with $a \le b$ we have $\ms( [a,b]) = \ms((a,b)) = b-a.$
\end{prop}
\begin{proof}
First consider the closed interval $[a,b]$. It is covered by the
single interval $U_1 = (a-\eps, b+\eps)$ so
$\ms([a,b]) \le \len( U_1) = b-a + 2\eps.$  Since 
$2\eps$  is arbitrary we conclude that 
$\ms([a,b]) \le b-a.$

On the other hand by the \myindex{Heine-Borel Theorem}
any open covering of $[a,b]$ has
a finite subcovering so it suffices to prove that for any finite
cover $\{U_i\}_{i=1}^n$ we have $\sum \len(U_i) \ge b-a$ as this will 
imply $\ms([a,b]) \ge b-a.$  We prove this by induction on $n$ the
number of elements in the cover by open intervals.  Clearly the result
holds if $n=1.$  If $n >1$ we note that two of the open intervals
must intersect.  This is because one of the intervals (say $(c,d)$) contains
$b$ and another interval contains $c$ and hence these two intersect.  
By renumbering the intervals we can assume that $U_{n-1}$ and $U_n$ intersect. 

Now define $V_{n-1} = U_{n-1} \cup U_n$ and $V_i = U_i$ for
$i < n-1.$  Then $\{V_i\}$ is an open cover of
$[a,b]$ containing $n-1$ intervals.  By the induction hypothesis
\[
\sum_{i=1}^{n-1} \len(V_i) \ge b-a.
\]
But $\len( U_{n-1}) + \len( U_n) > \len( V_{n-1})$ and
$\len( U_i) = \len( V_{i-1})$ for $i>2.$  Hence
\[
\sum_{i=1}^{n} \len(U_i) > \sum_{i=1}^{n-1} \len(V_i) \ge b-a.
\]
This completes the proof that $\ms([a,b]) \ge b-a$  and
hence that $\ms([a,b]) = b-a.$

For the open interval $(a,b)$ we note that $U =(a,b)$
covers itself so $\ms((a,b)) \le b - a.$  On the other hand
any cover $\{U_i\}_{i=1}^\infty$ of $(a,b)$ by open intervals
is also a cover of the closed interval $[a+\eps, b-\eps]$
so, as we just showed, 
\[
\sum_{i=1}^\infty \len(U_i) \ge b-a - 2\eps.
\]
As $\eps$  is arbitrary $\sum \len(U_i) \ge b-a$ and hence
$\ms((a,b)) \ge b-a$ which completes our proof.
\end{proof}

Two special cases are worthy of note:

\begin{cor}\label{cor: point}
The outer measure of a set consisting of single point is $0$.  The
outer measure of the empty set is also $0$.
\end{cor}

Lebesgue outer measure satisfies a monotonicity property with
respect to inclusion.

\begin{prop}\label{prop: outer monotonicity}
If $A$ and $B$ are bounded subsets of $\R$ and $A \subset B$
then $\ms(A) \le \ms(B).$
\end{prop}

\begin{proof}
Since $A \subset B$, every countable cover
$\{U_n\} \in \U(B)$ of $B$ by open intervals is also
in $\U(A)$ since it also covers $A$.  Thus
\[
\inf_{\{U_n\} \in \U(A)} \Big\{\sum_{n=1}^\infty \len ( U_n) \Big\}\le
\inf_{\{U_n\} \in \U(B)} \Big\{\sum_{n=1}^\infty \len ( U_n)\Big\},
\]
so $\ms(A) \le \ms(B).$
\end{proof}

We can now prove the first part of the countable additivity property
we want.  It turns out that this is the best we can do if we want our
measure defined on all bounded sets.  Note that the following result
is stated in terms of a countably infinite collection
$\{A_n\}_{n=1}^\infty$ of sets, but it is perfectly valid for a 
finite collection also.

\begin{thm}[Countable Subadditivity]\label{thm: outer subadditivity}
\index{countable subadditivity}
If $\{A_n\}_{n=1}^\infty$ is a countable collection of
bounded subsets of $\R$, then
\[
\ms( \bigcup_{n=1}^\infty A_n) \le \sum_{n=1}^\infty \ms(A_n)
\]
\end{thm}
\begin{proof}
By the definition of outer measure we know that each $A_n$ has 
a countable cover by open intervals $\{U_i^n\}$ such that
\[
\sum_{i=1}^\infty \len(U_i^n) \le \ms(A_n) + 2^{-n}\eps.
\]
But the union of all these covers $\{U_i^n\}$ is a countable cover of
$\bigcup_{n=1}^\infty A_n.$  So
\[
\ms( \bigcup_{n=1}^\infty A_n) \le \sum_{n=1}^\infty \sum_{i=1}^\infty \len(U_i^n)
\le \sum_{n=1}^\infty \ms(A_n) + \sum_{n=1}^\infty 2^{-n}\eps
= \sum_{n=1}^\infty \ms(A_n) + \eps.
\]
Since this is true for every $\eps$ the result follows.  The result
for a finite collection $\{A_n\}_{n=1}^m$ follows from this by
letting $A_i =\emptyset$ for $i > m.$
\end{proof}

\begin{cor}
If $A$ is countable, then $\ms(A) = 0$
\end{cor}

\begin{proof}
Suppose $A = \bigcup_{i=1}^\infty \{ x_i \}.$  
We saw in Corollary (\ref{cor: point}) that
$\ms(\{x_i\}) = 0$ so
\[
\ms( A) = \ms( \bigcup_{i=1}^\infty \{x_i\}) \le \sum_{i=1}^\infty \ms(\{x_i\}) = 0.
\]
which implies $\ms( A) = 0.$
\end{proof}

Since countable sets have outer measure $0$ and $\ms([a,b]) = 
b-a$ we also immediately obtain the following non-trivial result
(cf. part 4. of Exercise~(\ref{exer: uncountable})).

\begin{cor}\label{cor: reals uncountable}
If $a < b$, then $[a,b]$ is not countable.
\end{cor}

Outer Lebesgue measure satisfies property II of those we
enumerated at the beginning, namely it is translation invariant.
\begin{thm}\label{thm: trans inv}
If $c \in \R$ and $A$ is a bounded subset of $\R$, 
then $\ms(A) = \ms(A +c)$
where $A +c = \{ x + c \ | \ x \in A\}.$
\end{thm}
We leave the (easy) proof as an exercise.

\begin{exer}\ 
\begin{enumerate}
\item Prove Theorem (\ref{thm: trans inv}).
\item Prove that given $\eps >0$ 
there exist a countable collection of open intervals
$U_1,U_2,\dots, U_n, \dots$ such that $\bigcup_n U_n$ contains all rational
numbers in $\R$ and such that 
$\sum_{n=1}^\infty \len(U_n) = \eps.$
\item Give an example of a subset $A$ of $I$ such that
$\ms(A) = 0$, but with the property that if $U_1, U_2, \dots, U_n$ is
a {\em finite} cover by open intervals, then $\sum_{i=1}^n \len(U_i) \ge 1.$
\end{enumerate}
\end{exer}

\section{Lebesgue Measurable Sets}\label{sec: Leb measurable sets}

In Definition~(\ref{def: leb sigma}) we defined the $\sigma$-algebra
$\M$ to be the $\sigma$-algebra of subsets of $\R$ generated by open
intervals and null sets (it is also the $\sigma$-algebra of subsets of
$\R$ generated by Borel sets and null sets).  We defined a set to
be \myindex{Lebesgue measurable} if it is in this
$\sigma$-algebra. However now, in order to prove the existence of
Lebesgue measure, we want to use a different, but equivalent
definition.

Our program is roughly as follows:
\begin{itemize}
\item We will define a collection $\M_0$ of subsets of $I$
The criterion used to define $\M_0$ is often given as the 
definition of Lebesgue measurable sets.
\item We will define the Lebesgue measure $\mu(A)$ of a set $A$ in $\M_0$
to be the outer measure of $A$.
\item We will show that
the collection $\M_0$  is 
a $\sigma$-algebra of subsets of $I$ and in fact precisely 
the $\sigma$-algebra $\M(I)$ so we have defined $\mu(A)$
for all $A \in \M(I).$
\item We will prove that $\mu$ defined in this way satisfies
the properties promised in Chapter~\ref{chap: leb meas}, 
namely properties I-V of Theorem~(\ref{thm: lebesgue measure}).
Several of these properties follow from the corresponding properties
for outer measure $\ms$, which we proved in Section~(\ref{sec: outer}).
\end{itemize}

Henceforth for definiteness we will consider subsets of
the unit interval $I = [0,1].$  We could, of course, use any
other closed interval or even, with some extra work, the whole
real line.  Lebesgue outer measure as in Definition (\ref{def: outer})
has most of the properties we want. There is one serious problem,
however; namely, there exist subsets $A$ and $B$ of $I$ such that
$A \cap B = \emptyset$ and  $A \cup B = I$ but 
$\ms(A) + \ms(B) \ne \ms(I).$  That is, the additivity property
fails even with two sets whose union is an interval.

Fortunately, the sets for which it fails are rather exotic and
not too frequently encountered.  Our strategy is to restrict our
attention to only certain subsets of $I$ which we will call
``measurable'' and to show that on these sets $\ms$ has all the
properties we want.

If $A \subset I$ we will denote the \myindex{complement} of $A$
by $A^c$, that is, 
\[
A^c = I \setminus A = \{x \in I \ |\ x \notin A\}.
\]

\begin{defn}[Alternate Definition of Lebesgue Measurable]
\label{altdef: measurable}
\index{Lebesgue measurable} \index{measurable}
Let $\M_0$ denote the collection of all subsets of $I$
defined as follows: A subset $A$ of $I$ is in $\M_0$
provided for {\em any} subset $X \subset I$
\[
\ms(A \cap X) + \ms(A^c \cap X)  = \ms(X).
\]
For any set $A \in \M_0$ we define
$\mu(A)$ to be $\ms(A).$ 
\end{defn}

The goal of the remainder of this section is to prove
that in fact $\M_0$ is nothing other than the $\sigma$-algebra $\M(I)$
of Lebesgue measurable subsets of $I$ and the function 
$\mu: \M_0 \to \R$ satisfies the properties for Lebesgue
measure we claimed in Chapter~\ref{chap: leb meas}.
The defining condition above for a set $A$ to be in $\M_0$
is often taken as the definition of a
\myindex{Lebesgue measurable} subset of $I$ because it is what is
needed to prove the properties we want for Lebesgue 
measure.  Since we have already given a different definition
of Lebesgue measurable sets in Definition~(\ref{def: leb sigma})
we will instead prove the properties of $\M_0$ and $\mu$
which we want and, then show $\M_0 = \M(I)$ so the two
definitions coincide.
Indeed, we will prove in Corollary~(\ref{cor: equal sigma}) that the
sets in $\M_0$ are precisely the sets in $\M(I)$ 
the $\sigma$-algebra generated by Borel subsets and null subsets
of $I$.

As a first step we show that $\M_0 \subset \M(I).$

\begin{prop}\label{prop: M0}
Every set $A \in \M_0$ can be written as 
\[
A = B \setminus N
\]
where $B$ is the intersection of a countable nested family
of open sets (and,
in particular, is a Borel set) and $N = A^c \cap B$ is a null set. It follows
that $\M_0 \subset \M(I).$
\end{prop}

\begin{proof}
Since $\mu(A) = \ms(A)$ it follows from the definition
of outer measure that for any $\eps > 0$ there is a
cover $\V_\eps$ of $A$ by open intervals $U_n$ such
that 
\[
\sum_{n=1}^\infty \ms(U_n) = \sum_{n=1}^\infty \len(U_n) < \ms(A) + \eps.
\]
The monotonicity of outer measure, then implies that
if $V_\eps = \cup_{n=1}^\infty U_n$ we have
$\ms(A) \le \ms(V_\eps) \le \ms(A) + \eps.$
If we let 
\[
W_k = \bigcap_{i=1}^k V_{\frac{1}{k}}
\], then each $\{W_k\}$ is a nested family of open
sets and $\ms(A) \le \ms(W_k) \le \ms(A) + 1/k.$
open.  

Let $B = \cap_{k=1}^\infty W_k$.
By monotonicity again we have
\[
\ms(A) \le \ms(B) \le \ms(V_k) < \ms(A) + \frac{1}{k}.
\]
Since this holds for all $k >0$ we conclude $\ms(B) = \ms(A).$

In the defining equation of $\M_0$ (see Definition~(\ref{altdef: measurable})
we take $X = B$ and
conclude
\[
\ms(A \cap B) + \ms(A^c \cap B)  = \ms(B).
\]
Since $A \subset B$ we have  $A \cap B = A$
and hence $\ms(A) + \ms(A^c \cap B)  = \ms(B).$
From the fact that $\ms(B) = \ms(A)$ it follows that $\ms(A^c \cap B) = 0$.
Therefore if $N = A^c \cap B$, then $N$ is a null set.
Finally $A = B \setminus (B \cap A^c)$ so
$A = B \setminus N.$
\end{proof}

The definition of $\M_0$ is relatively
simple, but to show it has the properties we want
requires some work.  If we were to replace the $=$ sign
in this definition with $\ge$
we would obtain a statement which is true for 
all subsets of $I$.  So to prove a set is in $\M_0$ we need only check
the reverse inequality.  More precisely, 

\begin{prop}\label{prop: meas}\
Suppose $A \subset I$, then
\begin{enumerate}
\item[(1)] The set $A$ is in $\M_0$
provided for {\em any} subset $X \subset I$
\[
\ms(A \cap X) + \ms(A^c \cap X)  \le \ms(X).
\]
\item[(2)] The set $A$ is in $\M_0$ if and only if  $A^c$ is
in $\M_0$.  In this case $\mu( A^c) = 1 - \mu(A).$
\end{enumerate}
\end{prop}

\begin{proof}
For part (1) observe $X = (A \cap X) \cup (A^c \cap X)$ 
so the subadditivity property
of outer measure in Theorem (\ref{thm: outer subadditivity}) tells us that
\[
\ms(A \cap X) + \ms(A^c \cap X)  \ge \ms(X)
\]
always holds. This plus the inequality of our hypothesis gives
the equality of the definition of $\M_0.$

For part (2) suppose $A$ is an arbitrary subset of $I$.
The fact that $(A^c)^c = A$ implies immediately from 
Definition~(\ref{altdef: measurable}) that $A$ is in $\M_0$
if and only if $A^c$ is.  Also
taking $X = I$ in this definition we conclude
\[
\ms(A \cap I) + \ms(A^c \cap I)  = \ms(I) = 1.
\]
So $\mu(A) + \mu(A^c) = 1.$
\end{proof}

\begin{prop}\label{prop: null-complements}
A set $A \subset I$ is a null set if and only if $A \in \M_0$
and $\mu(A) = 0.$ 
\end{prop}
\begin{proof}
By definition  a set $A$ is a null set if and only if
$\ms(A) = 0.$
If $A$ is a null set, then since $A \cap X \subset A$  we know
by the monotonicity of outer measure 
(Proposition~(\ref{prop: outer monotonicity}))
that $\ms(A \cap X) = 0.$
Similarly, since $A^c \cap X \subset X$  we know
that $\ms(A^c \cap X) \le \ms(X).$ 
Hence again using monotonicity of outer measure 
from Proposition (\ref{prop: outer monotonicity}) we
know that 
\[
\ms(A \cap X) + \ms(A^c \cap X) = \ms(A^c \cap X) \le \ms(X)
\]
and the fact that $A \in \M_0$ follows from
part (1) Proposition (\ref{prop: meas}.
\end{proof}

\begin{prop}\label{prop: union}
If $A$ and $B$ are in $\M_0$, then $A \cup B$ and $A \cap B$ are
in $\M_0$.
\end{prop}
\begin{proof}
To prove if two sets, $A$ and $B$, then there union 
is in $\M_0$ requires some work.  Suppose $X \subset I.$  
And since $(A \cup B) \cap X = (B \cap X) \cup (A \cap B^c \cap X),$
the subadditivity of Theorem (\ref{thm: outer subadditivity}) tells us 
\begin{equation}\label{eqn:union1}
\ms( (A \cup B) \cap X) \le \ms(B \cap X) + \ms(A \cap B^c \cap X).
\end{equation}

Also the definition of $\M_0$ tells us
\begin{equation}\label{eqn:union2}
\ms(B^c \cap X) = \ms( A \cap B^c \cap X) + \ms( A^c \cap B^c \cap X).
\end{equation}

Notice that $(A \cup B)^c = A^c \cap B^c.$ 
So we get
\begin{align*}
\ms( (A &\cup B) \cap X) + \ms( (A \cup B)^c \cap X) & \\
&= \ms( (A \cup B) \cap X) + \ms( A^c \cap B^c \cap X) & \\
&\le \ms(B \cap X) + \ms(A \cap B^c \cap X) + \ms( A^c \cap B^c \cap X)
\hskip 0.5cm &\text{by equation (\ref{eqn:union1}),}\\
&= \ms(B \cap X) + \ms(B^c \cap X)
&\text{by equation (\ref{eqn:union2}),}\\
&= \ms(X).&
\end{align*}
According to part (1) of Proposition (\ref{prop: meas}) this
implies that $A \cup B$ is in $\M_0.$

The intersection now follows easily using what we know about the union
and complement.  More precisely, $A \cap B = (A^c \cup B^c)^c$ so if
$A$ and $B$ are in $\M_0$, then so is $(A^c \cup B^c)$ and hence its
complement $(A^c \cup B^c)^c$ is also.
\end{proof}

Next we wish to show intervals are in $\M_0.$

\begin{prop}\label{prop: interval meas}
Any subinterval of $I$, open, closed or half open,in $\M_0.$
\end{prop}

\begin{proof}
First consider $[0,a]$ with complement $(a,1].$ 
If $X$ is an arbitrary  subset of $I$ we must show
$\ms([0,a] \cap X) + \ms((a,1] \cap X) = \ms(X).$
Let $X^- = [0,a] \cap X$ and $X^+ = (a,1] \cap X$. 
Given $\eps >0$, the definition of outer measure tells us
we can find a countable cover of $X$ by open
intervals $\{U_n\}_{n=1}^\infty$ such that 
\begin{equation}\label{eqn:interval}
\sum_{n=1}^\infty \len(U_n) \le \ms(X) + \eps.
\end{equation}

Let $U_n^- = U_n \cap [0,a]$ and $U_n^+ = U_n \cap (a,1]$.  Then
$X^- \subset \bigcup_{n=1}^\infty U_n^-$ and
$X^+ \subset \bigcup_{n=1}^\infty U_n^+$.  Subadditivity of 
outer measure implies
\[
\ms(X^-) \le \ms(\bigcup_{n=1}^\infty U_n^-) \le \sum_{n=1}^\infty \len(U_n^-)
\]
and
\[
\ms(X^+) \le \ms(\bigcup_{n=1}^\infty U_n^+) \le \sum_{n=1}^\infty \len(U_n^+).
\]

Adding these inequalities and using equation (\ref{eqn:interval}) we get
\begin{align*}
\ms(X^-) + \ms(X^+) &\le \sum_{n=1}^\infty \len(U_n^-) + \len(U_n^+)\\
&= \sum_{n=1}^\infty \len(U_n)\\
&\le \ms(X) + \eps.
\end{align*}
Since $\eps$ is arbitrary we conclude that
$\ms(X^-) + \ms(X^+) \le \ms(X)$ which by Proposition (\ref{prop: meas})
implies that $[0,a]$ is in $\M_0$ for any $0 \le a \le 1.$  A
similar argument implies that $[a,1]$ is in $\M_0$.  Taking complements,
unions and intersections it is clear that any interval, open closed or
half open, is in $\M_0$.
\end{proof}

\begin{lemma}\label{lem: additive}
Suppose $A$ and $B$ are disjoint sets in $\M_0$ and
$X \subset I$ is arbitrary.  Then 
\[
\ms( (A \cup B) \cap X) = \ms( A \cap X) + \ms( B \cap X).
\]
The analogous result for a finite union of disjoint measurable
sets is also valid.
\end{lemma}
\begin{proof}
It is always true that 
\[
A \cap (A \cup B) \cap X = A \cap X.
\]
Since $A$ and $B$ are disjoint
\[
A^c \cap (A \cup B) \cap X = B \cap X.
\]
Hence the fact that $A$ is in $\M_0$ tells us
\begin{align*}
\ms((A \cup B) \cap X) &= \ms(A \cap (A \cup B) \cap X) + 
\ms(A^c \cap (A \cup B) \cap X)\\
&= \ms(A \cap X) + \ms(B \cap X).
\end{align*}

The result for a finite collection $A_1, A_2, \dots , A_n$ follows
immediately by induction on $n$. 
\end{proof}

\begin{thm}\label{thm: countable union}
The collection $\M_0$ of subsets of $I$
is closed under countable unions and countable
intersections.
Hence $\M_0$ is a $\sigma$-algebra.
\end{thm}

\begin{proof}
We have already shown that 
the complement of a set in $\M_0$ is a set in $\M_0.$

We have also shown that the union or intersection of a finite
collection of 
sets in $\M_0$ is a set in $\M_0$.

Suppose $\{ A_n\}_{n=1}^\infty$ is a countable collection of 
sets in $\M_0.$  We want to construct a countable collection
of pairwise disjoint sets $\{ B_n\}_{n=1}^\infty$ which
are in $\M_0$ and have the same union.

To do this we define $B_1 = A_1$ and 
\[
B_{n+1} = A_{n+1} \setminus \bigcup_{i=1}^n A_n 
= A_{n+1} \cap \Big (\bigcup_{i=1}^n A_n \Big)^c.
\]

Since finite unions, intersections and complements of sets
in $\M_0$ are
sets in $\M_0$, it is clear that $B_n$ is measurable.  Also
it follows easily by induction that $\bigcup_{i=1}^n B_i = \bigcup_{i=1}^n A_i$
for any $n$.  Thus $\bigcup_{i=1}^\infty B_i = \bigcup_{i=1}^\infty A_i$

Hence to prove $\bigcup_{i=1}^\infty A_i$ is in $\M_0$ we will prove
that $\bigcup_{i=1}^\infty B_i$ is in $\M_0$.  Let 
$F_n = \bigcup_{i=1}^n B_i$ and $F = \bigcup_{i=1}^\infty B_i$. 
If $X$ is an arbitrary subset of $I$, then since $F_n$ is in $\M_0$
\[
\ms(X) = \ms( F_n \cap X) + \ms( F_n^c \cap X) 
\ge \ms( F_n \cap X) + \ms( F^c \cap X)
\]
since $F^c \subset F_n^c.$
By Lemma (\ref{lem: additive})
\[
\ms(F_n \cap X) = \sum_{i=1}^n \ms(B_i \cap X).
\]

Putting these together we have 
\[
\ms(X) \ge \sum_{i=1}^n \ms(B_i \cap X) + \ms( F^c \cap X)
\]
for all $n >0.$ Hence
\[
\ms(X) \ge \sum_{i=1}^\infty \ms(B_i \cap X) + \ms( F^c \cap X).
\]
But subadditivity of $\ms$ implies
\[
\sum_{i=1}^\infty \ms(B_i \cap X) \ge \ms(\bigcup_{i=1}^\infty(B_i \cap X)) 
= \ms(F \cap X).
\]
Hence 
\[
\ms(X) \ge \ms(F \cap X)) + \ms( F^c \cap X)
\]
and $F$ is in $\M_0$ by Proposition (\ref{prop: meas}).

To see that a countable intersection of sets in $\M_0$ is in $\M_0$
we observe that 
\[
\bigcap_{n=1}^\infty A_n = \Big ( \bigcup_{n=1}^\infty A_n^c \Big )^c
\]
so the desired result follows from the result on unions together 
with the fact that  $\M_0$ is closed under taking complements.
\end{proof}

\begin{cor}\label{cor: equal sigma}
The $\sigma$-algebra $\M_0$ of 
subsets of $I$ equals $\M(I)$ the $\sigma$-algebra of subsets
of $I$ generated by 
Borel sets and null sets.
\end{cor}

\begin{proof}
The $\sigma$-algebra $\M_0$ contains open intervals and
closed intervals in $I$ by 
Proposition~(\ref{prop: interval meas}) and hence contains
the $\sigma$-algebra they generate, the Borel subsets
of $I$.
Also $\M_0$ contains null sets by
Proposition~(\ref{prop: null-complements}).
Therefore $\M_0$ contains $\M(I)$ the 
$\sigma$-algebra $\M_0$ generated by Borel sets and
null sets.

On the other hand by Proposition~(\ref{prop: M0}) $\M_0 \subset \M(I).$
Hence $\M_0 = \M(I).$
\end{proof}

Since we now know the sets in $\M_0$, i.e.  the sets which satisfy
Definition~(\ref{altdef: measurable}), coincide with the sets in
$\M[I]$, we will refer to them sets as \myindex{Lebesgue measurable
sets}, or simply \myindex{measurable sets} for short.  We also no
longer need to use outer measure, but can refer to the Lebesgue
measure $\mu(A)$ of a measurable set $A$ (which, of course, has the
same value as the outer measure $\ms(A)).$

\begin{thm}[Countable Additivity]\label{thm: countable additivity}
\index{countable additivity}
If $\{ A_n\}_{n=1}^\infty$ is a countable collection of 
measurable subsets of $I$, then 
\[
\mu( \bigcup_{n=1}^\infty A_n) \le \sum_{n=1}^\infty \mu(A_n).
\]
If the sets are pairwise disjoint, then
\[
\mu( \bigcup_{n=1}^\infty A_n) = \sum_{n=1}^\infty \mu(A_n).
\]
The same equality and inequality are valid for a finite collection of
measurable subsets $\{A_n\}_{n=1}^m.$
\end{thm}
\begin{proof}
The first inequality is simply a special case of the 
subadditivity from Theorem (\ref{thm: outer subadditivity}).
If the sets $A_i$ are pairwise disjoint, then by Lemma (\ref{lem: additive})
we know that for each $n$
\[
\mu( \bigcup_{i=1}^\infty A_i) \ge \mu( \bigcup_{i=1}^n A_i) 
= \sum_{i=1}^n \mu(A_i).
\]
Hence
\[
\mu( \bigcup_{i=1}^\infty A_i) \ge \sum_{i=1}^\infty \mu(A_i).
\]
Since the reverse inequality follows from subadditivity we have
equality.
\end{proof}

We can now prove the main result of this Appendix,
which was presented as Theorem~(\ref{thm: lebesgue measure}) 
in Chapter~2.

\begin{thm}[Existence of Lebesgue Measure]
\index{Lebesgue measure}
There exists a unique  function $\mu$, called {\it Lebesgue measure},
from $\M(I)$ to the non-negative real numbers satisfying:
\begin{description}
\item[I. Length.] If $A = (a,b)$ 
then $\mu(A) = \len (A) = b-a,$ i.e. the measure of an open interval
is its length
\item [II. Translation Invariance.]
Suppose $A \subset I,\ c\in \R$ and $A +c \subset I$
where $A + c$ denotes the set $\{ x +c\ |\ x \in A\}.$
Then  $\mu( A + c ) = \mu(A)$

\item [III. Countable Additivity.]
If $\{A_n\}_{n=1}^\infty$ is a countable collection of
subsets of $I$, then
\[
\mu( \bigcup_{n=1}^\infty A_n) \le \sum_{n=1}^\infty \mu(A_n)
\]
and if the sets are  {\em pairwise disjoint}, then
\[
\mu( \bigcup_{n=1}^\infty A_n) = \sum_{n=1}^\infty \mu(A_n)
\]
\item [IV. Monotonicity] If $A,B \in \M(I)$ and $A \subset B$
then $\mu(A) \le \mu(B)$
\item [V. Null Sets] A subset $A \subset I$ is a null set
set if and only if $A \in \M(I)$ and $\mu(A) = 0.$
\end{description}
\end{thm}

\begin{proof}
The Lebesgue measure $\mu(A)$ of any set $A \in \M(I)$ is
defined to be its outer measure $\ms(A).$  Hence
properties I, II, and IV for $\mu$ follow from the corresponding properties
of $\ms.$  These were established in
Propostion~(\ref{prop: outer length}),
Theorem (\ref{thm: trans inv}), and 
Proposition (\ref{prop: outer monotonicity}) respectively.

Property III, countable additivity,  was proved in
Theorem~(\ref{thm: countable additivity}).
And Property V is a consequence of 
Proposition~(\ref{prop: null-complements}).

We are left with the task of showing that $\mu$ is 
unique.  Suppose $\mu_1$ and $\mu_2$ are two functions
defined on $\M(I)$ and satisfying properties I-V.
They must agree on any open interval by property I.
By Theorem~(\ref{thm: characterize open}) any open
set is a countable union of pairwise disjoint open
intervals, so countable additivity implies 
$\mu_1$ and $\mu_2$ agree on open sets.

Suppose that $B$ is the intersection of a countable nested
family of open sets $U_1 \supset U_2 \supset \dots \supset U_n \dots$.
Then Proposition~(\ref{prop: increasing union}) implies
\[
\mu_1(B) = \mu_1 \big(\bigcap_{n=1}^\infty U_n \big) = 
\lim_{n \to \infty} \mu_1(U_n)
= \lim_{n \to \infty} \mu_2(U_n)
=  \mu_2 \big(\bigcap_{n=1}^\infty U_n \big) = \mu_2(B).
\]
Also if $N \in \M(I)$ is a null set, then
$\mu_1(N) = 0 = \mu_2(N)$ by property V.

Finally if $A$ is an arbitrary set in $\M(I)$ by 
Proposition~(\ref{prop: M0})
$A = B \setminus N$ where $B$ is the intersection
of a countable nested family of open sets and
$N = A^c \cap B$ is a null set. Since $B$ is the disjoint
union of $A$ and $N$
It follows that
\[
\mu_1(A) = \mu_1(B) - \mu_1(N) = \mu_2(B) - \mu_2(N) = \mu_2(A).
\]

\end{proof}

\chapter{A Non-measurable Set}\label{sec: non-measurable}
\setcounter{section}{1}
\setcounter{thm}{0}

We are now prepared to prove the existence of a 
\myindex{non-measurable set}.  
The proof (necessarily) depends on the \myindex{Axiom of Choice}
(see Section (\ref{sec: sets}))
and is highly non-constructive.

\begin{lemma}\label{lem: difference}
Let $A$ be a measurable set with $\mu(A)> 0$  and let
$\Delta = \{ x_1 - x_2 \ | \ x_1,x_2 \in A\}$ be the set
of differences of elements of $A$.  Then for some $\eps> 0$
the set $\Delta$ contains the interval $(-\eps, \eps).$
\end{lemma}
\begin{proof}
By Theorem (\ref{thm: easy lebesgue density}) there is an open interval
$U$ such that $\mu( A \cap U) > \frac{3}{4}\len(U).$  Let $\eps
= \len(U)/2$, so $\len(U) = 2\eps$.
Suppose $y \in (-\eps, \eps)$ and let
$U + y = \{ x +y\ |\ x \in U\},$, then $U \cup (U + y)$ is an interval
of length at most $3\eps$.

Now let $B = A \cap U$ and $B' = B + y$.
Then $\mu(B') = \mu(B) > \frac{3}{4}\len(U) = \frac{3}{2}\eps$
so $\mu(B') + \mu(B) > 3\eps.$  On the other hand the fact that
$B \cup B' \subset U \cup (U + y)$ implies $\mu( B \cup B') \le 3\eps.$
It follows that $B$ and $B'$ cannot be disjoint since otherwise we would
contradict additivity.

If $x_1 \in B \cap B'$, then $x_1 = x_2 + y$ for some $x_2 \in B$.
Hence $y = x_1 - x_2 \in \Delta.$  We have shown that any
$y \in (-\eps, \eps)$ is in $\Delta.$
\end{proof}

\begin{thm}[Non-measurable Set]
\index{non-measurable set}
There exists a subset $E$ of $[0,1]$ which is not Lebesgue measurable.
\end{thm}
\begin{proof}
Let $\Q \subset \R$ denote the rational numbers.  The rationals are
an additive subgroup of $\R$ and we wish to consider the ``cosets''
of this subgroup.  More precisely, we want to consider the sets of
the form $\Q + x$ where $x \in \R$.  

We observe that two such 
sets $\Q + x_1$ and $\Q + x_2$ are either equal or disjoint.  This is
because the existence of one point $z \in (\Q + x_1) \cap (\Q + x_2)$
implies $z = x_1 +r_1 = x_2 +r_2$ with $r_1, r_2 \in \Q$, 
so $x_1 -  x_2 = (r_2 - r_1) \in \Q.$  This, in turn implies that
$\Q + x_2 = \{ x_2 + r\ |\ r\in \Q\} = \{ x_2 + r +(x_1-x_2)\ |\ r\in \Q\}
= \{ x_1 + r\ |\ r\in \Q\} = \Q + x_1.$

Using the Axiom of Choice \index{Axiom of Choice}
we construct a set $E$ which contains one
element from each of the cosets $\Q + x$, that is, for any $x_0 \in \R$
the set $E \cap (\Q + x_0)$ contains exactly one point.  Now let
$\{r_n\}_{n=1}^\infty$ be an enumeration of the rational numbers.
We want to show that $\R = \bigcup_{n=1}^\infty E + r_n$.  To see this
let $x \in \R$ be arbitrary and let $\{x_0\} = E \cap (\Q + x).$  Then
$x_0 = x + r$ for some $r \in \Q$ or $x = x_0 + r_0$ where $r_0 = -r.$
Hence $x \in E - r$ so $x \in \bigcup_{n=1}^\infty E + r_n.$
We have shown $\R = \bigcup_{n=1}^\infty E + r_n$.

We now make the assumption that $E$ is measurable
and show this leads to a contradiction,  We first note
that if we define
$\Delta = \{ x_1 - x_2 \ | \ x_1,x_2 \in E\}$, then $\Delta$
contains no rational points except $0$.  This is because
$x_1 = x_2 + r$ for rational $r$ would imply that
$E \cap (\Q +x_2) \supset \{x_1,x_2\}$ and this intersection 
contains only one point.  Since $\Delta$ contains at most
one rational point it cannot contain an open interval so
by Lemma (\ref{lem: difference}) we must conclude that
$\mu(E) = 0.$

But if we define $V_n = (E + r_n) \cap [0,1]$, then 
$\mu(V_n) \le \mu(E + r_n) = \mu(E) = 0$ so $\mu(V_n)=0.$
The fact that 
$\R = \bigcup_{n=1}^\infty E + r_n$ implies
\[
[0,1] = \bigcup_{n=1}^\infty \big ( (E + r_n) \cap [0,1]\big )
= \bigcup_{n=1}^\infty V_n. 
\]
Subadditivity, then would imply $\mu([0,1]) \le \sum_{n=1}^\infty \mu(V_n) = 0$
which clearly contradicts our assumption that $E$ is measurable.
\end{proof}
\vfill\eject

\vfill\eject

\printindex

\begin{thebibliography}{100}

\bibitem{D}
J. Dieudonn\'e
\newblock {\em Foundations of Modern Analysis,}
\newblock Academic Press, New York. (1960).


\bibitem{Mar}
J.~Marsden and M.~Hoffman
\newblock {\em Elementary Classical Analysis,}
\newblock W.H. Freeman, (1993).

\bibitem{Pin}
M.~Pinsky
\newblock {\em Partial Differential Equations and Boundary Value Problems with 
Applications,}
\newblock McGraw Hill, (1998).

\bibitem{R}
H.~L.~Royden
\newblock {\em Real Analysis,}
\newblock Macmillan, New York. (1963).

\bibitem{T}
S.~J.~Taylor
\newblock {\em Introduction to Measure and Integration,}
\newblock Cambridge University Press, London. (1966).

\bibitem{Z}
A.~C.~Zaanen
\newblock {\em An Introduction to the Theory of Integration,}
\newblock North-Holland, New York. (1961).

\end{thebibliography}
\end{document}